\documentclass[10pt]{amsart}
\usepackage[dvipdfmx]{graphicx}
\usepackage{amsmath,amscd,amsthm,amsxtra,yhmath,stackrel}
\usepackage{epsfig,graphics,color,colortbl}
\usepackage{amssymb,latexsym} 
\usepackage{mathrsfs}
\usepackage{bm}
\usepackage[poly,all]{xy}
\usepackage{hyperref}
\usepackage{tikz-cd}
\usepackage[draft]{todonotes}
\usepackage{upgreek}
\usepackage{ytableau}
\usepackage{stmaryrd}
\usepackage{verbatim}
\usepackage{mathtools}
\usepackage[title]{appendix}
\usepackage{enumerate}
\usepackage[normalem]{ulem}
\usepackage{datetime}
\usepackage[pagewise]{lineno} 
\usepackage{arydshln}
\usetikzlibrary{decorations.pathreplacing}
\usepackage{adjustbox}

\newtheorem{thm}{\bf Theorem}[section]
\newtheorem{df}[thm]{\bf Definition}
\newtheorem{prop}[thm]{\bf Proposition}
\newtheorem{cor}[thm]{\bf Corollary}
\newtheorem{lem}[thm]{\bf Lemma}
\newtheorem{rem}[thm]{\bf Remark}
\newtheorem{ex}[thm]{\bf Example}

\newtheorem{alg}[thm]{\bf Algorithm}

\numberwithin{equation}{section}

\newcommand{\mc}{\mathcal}
\newcommand{\mf}{\mathfrak}
\newcommand{\ms}{\mathscr}

\newcommand{\pf}{\noindent{\bfseries Proof. }}
\newcommand{\ov}{\overline}
\newcommand{\un}{\underline}

\newcommand{\K}{{\mc K}}
\newcommand{\U}{{\mc U}}

\newcommand{\Boson}{B_q}
\newcommand{\cmB}{\Delta_{\Boson}}

\newcommand{\cP}{\mathscr{P}}
\newcommand{\cO}{\mc{O}}
\newcommand{\I}{\mathbb{I}}

\newcommand{\bi}{{\bf i}}
\newcommand{\bj}{{\bf j}}
\newcommand{\ttq}{\texttt{q}}

\newcommand{\N}{\mathbb{N}}
\newcommand{\Z}{\mathbb{Z}}
\newcommand{\Q}{\mathbb{Q}}

\newcommand{\e}{\epsilon}
\newcommand{\de}{\delta}

\newcommand{\te}{\tilde{e}}
\newcommand{\tf}{\tilde{f}}
\newcommand{\tx}{\tilde{x}}
\newcommand{\td}{\widetilde}
\newcommand{\tte}{\td{\texttt e}}
\newcommand{\ttf}{\td{\texttt f}}

\newcommand{\gl}{\mf{gl}}

\newcommand{\g}{\mf{g}}

\newcommand{\La}{\Lambda}

\newcommand{\la}{\lambda}

\newcommand{\hf}{\frac{1}{2}}

\newcommand{\blue}[1]{{\color{blue}#1}}
\newcommand{\red}[1]{{\color{red}#1}}

\newcommand{\ot}{\otimes}

\newcommand{\sa}{(-1)^{\e_a}}

\newcommand{\bq}{{\bf q}}
\newcommand{\bff}{{\bf f}}

\newcommand{\bfF}{{\bf F}}

\newcommand{\bM}{{\bf M}}
\newcommand{\bV}{{\bf V}}
\newcommand{\bc}{{\bf c}}
\newcommand{\burge}{\kappa^\searrow}
\newcommand{\bt}{\btriangle}
\newcommand{\bti}{{\btriangle_i}}
\newcommand{\btic}{{\btriangle_i^c}}
\newcommand{\btiII}{{\btriangle_{i+2}}}
\newcommand{\ti}{{\triangle_i}}

\newcommand{\Bins}{\,{\overset{\,\,{\rm B}}{\longleftarrow}}\,}
\newcommand{\biw}[2]{\binom{\, #1 \,}{\, #2 \,}}
\newcommand{\biws}{{\Omega^{\mf{d}_{m|n}}}}
\newcommand{\ba}{{\bf a}}
\newcommand{\bb}{{\bf b}}
\newcommand{\tP}{\textup{\texttt P}}
\newcommand{\tQ}{\textup{\texttt Q}}

\newcommand{\tT}{\textup{\texttt T}}

\newcommand{\cii}{\mathsf{a}}
\newcommand{\ciI}{\mathsf{b}}
\newcommand{\cII}{\mathsf{c}}
\newcommand{\bsfT}{{\boldsymbol{\mathsf{T}}}}
\newcommand{\bsfX}{{\boldsymbol{\mathsf{X}}}}
\newcommand{\bsfY}{{\boldsymbol{\mathsf{Y}}}}
\newcommand{\bsfx}{{\boldsymbol{\mathsf{x}}}}
\newcommand{\bsfy}{{\boldsymbol{\mathsf{y}}}}

\newcommand{\qo}[1]{\un{#1}}
\newcommand{\re}[1]{\ov{#1}}

\newcommand{\tl}[1]{\substack{\scalebox{0.75}{#1}}}

\newcommand*{\parallelogramm}{%
  \rlap{\rotatebox{-30}{\rule[.04ex]{.4pt}{.6em}}}%
  \kern.03em%
  \rlap{\kern.26em\raisebox{0.5em}{\rule{.5em}{.4pt}}}%
  \rule{.45em}{.4pt}\kern-.03em%
  \rotatebox{-30}{\rule[.04ex]{.4pt}{.6em}}%
}

\makeatletter
\newcommand{\btriangle}{\mathpalette\btriangle@\relax}
\newcommand{\btriangle@}[2]{%
  \begingroup
  \sbox\z@{$\m@th#1\triangle$}%
  \makebox[\wd\z@]{%
    \raisebox{0.04\height}{%
      \resizebox{1.1\wd\z@}{0.96\ht\z@}{%
        $\m@th#1\blacktriangle$%
      }%
    }%
  }%
  \endgroup
}
\makeatother

\def\dotfill#1{\cleaders\hbox to #1{.}\hfill}
\newcommand\dotline[2][.5em]{\leavevmode\hbox to #2{\dotfill{#1}\hfil}}

\newcommand{\hdomino}{\scalebox{0.7}{\raisebox{0.1cm}{\!\!\!\xymatrix@R=-0.82em @C=-0.82em{\boxed{1}\boxed{1}}}}}
\newcommand{\vdomino}{\scalebox{0.7}{\raisebox{0.3cm}{\!\!\xymatrix@R=-0.82em @C=-0.82em{\boxed{1}\\\boxed{2}}}}}
\newcommand{\vd}[2]{\scalebox{0.8}{\raisebox{0.36cm}{\!\!\xymatrix@R=-0.82em @C=-0.82em{\boxed{#1}\\\boxed{#2}}}}}

\makeatletter
\newcommand*{\shifttext}[2]{%
  \settowidth{\@tempdima}{#2}%
  \makebox[\@tempdima]{\hspace*{#1}#2}%
}
\makeatother

\begin{document}
\title
[Crystal base of the quantum orthosymplectic superalgebra]{Crystal base of the negative half of quantum orthosymplectic superalgebra}

\author{IL-SEUNG JANG}
\address{Department of Mathematics, Incheon National University, Incheon 22012, Korea}
\email{ilseungjang@inu.ac.kr}

\author{JAE-HOON KWON}
\address{Department of Mathematical Sciences and RIM, Seoul National University, Seoul 08826, Korea}
\email{jaehoonkw@snu.ac.kr}

\author{AKITO URUNO}
\address{Research institute of Mathematics, Seoul National University, Seoul 08826, Korea}
\email{aki926@snu.ac.kr}

\thanks{I.-S. Jang is supported by Incheon National University Research Grant in 2023 (No. 2023-0205)}
\thanks{J.-H Kwon and A. Uruno are supported by the National Research Foundation of Korea(NRF) grant funded by the Korea government(MSIT) (No.2020R1A5A1016126). J.-H. Kwon is also partially supported by the Global-LAMP Program of the National Research Foundation of Korea (NRF) grant funded by the Ministry of Education (No. RS-2023-00301976)}

\begin{abstract}
We construct a crystal base of the negative half of a quantum orthosymplectic superalgebra. 
It can be viewed as a limit of the crystal bases of $q$-deformed irreducible oscillator representations. 
We also give a combinatorial description of the embedding from the crystal of a  $q$-oscillator representation to that of the negative half subalgebra given in terms of a PBW type basis. 
It is given as a composition of embeddings into the crystals of intermediate parabolic Verma modules, where the most non-trivial one is from an oscillator module to a maximally parabolic Verma module with respect to a quantum subsuperalgebra for $\mf{gl}_{m|n}$. 
A new crystal theoretic realization of Burge correspondence of orthosymplectic type plays an important role for the description of this embedding.
\end{abstract}

\maketitle
\setcounter{tocdepth}{1}

\noindent

\section{Introduction}
\subsection{Quantum orthosymplectic superalgebras and $q$-oscillator modules}%
The crystal base of the negative half of a quantized enveloping algebra associated with a symmetrizable Kac–Moody algebra is an object with a fundamental combinatorial structure, into which the crystal bases of integrable highest weight modules are naturally incorporated \cite{Kas91}. It can be obtained by taking the limit of crystal bases of integrable highest weight modules.

Let $\U(\g)$ be a quantized enveloping algebra over $\Bbbk=\mathbb{Q}(v)$ ($v=q$ or $q^{\hf}$) associated to an orthosymplectic superalgegbra $\g=\mf{b}_{m|n}, \mf{c}_{m|n}, \mf{d}_{m|n}$ ($m\ge 2$) with a standard fundamental system $\Pi=\{\,\alpha_i\,|\,i\in I\,\}$. 
In this paper we construct a crystal base of the negative half $\U(\g)^-$ of $\U(\g)$ as a continuation of our previous works \cite{JKU23,JKU24}.

In case of an orthosymplectic superalgebra $\g$ and even a general linear Lie superalgebra $\gl_{m|n}$, it is not known in general whether a finite-dimensional irreducible $\U(\g)$-module has a crystal base. On the other hand, a family of irreducible highest weight $\U(\g)$-modules was studied in \cite{K15,K16}, which are $q$-analogue of oscillator representations of $\g$ \cite{CZ} naturally arising from the theory of Howe dual pairs \cite{H}. It is shown that this module, which we denote by $V(\la,\ell)$ for $(\la,\ell)\in \cP(\g)$, has a unique crystal base $(\ms{L}(V(\la,\ell)),\ms{B}(V(\la,\ell)))$ together with a combinatorial model of the crystal $\ms{B}(V(\la,\ell))$ called spinor model. 
The classical limit of $V(\la,\ell)$ (after taking a limit of the rank $n$) corresponds to an integrable highest weight module over the Lie algebras of type $BCD$ under the theory of super duality \cite{CLW}. From this point of view, it would be natural to construct a crystal base of $\U(\g)$ by taking a limit of $(\ms{L}(V(\la,\ell)),\ms{B}(V(\la,\ell)))$.

Another major difference from the case of symmetrizable Kac-Moody algebras is that a $q$-oscillator representation does not have an explicit presentation as a quotient of $\U(\g)^-$.
This makes it more complicated to take the limit of crystals of $V(\la,\ell)$ and describe the associated crystal embeddings. 
So in order to overcome these difficulties, we consider a sequence of parabolic Verma modules following the ideas in \cite{JKU23} in case of $\gl_{m|n}$ and also use combinatorial objects and algorithms.
\vspace{-0.2cm}

\subsection{Parabolic Verma modules}
Let $\mf{l}$ be the subalgebra of $\g$ corresponding to $\Pi\setminus\{\alpha_0\}$ which is isomorphic to  $\mf{gl}_{m|n}$, and let $\U(\mf{l})$ be the associated subalgebra of $\U(\g)$. Let $V_{\mf l}(\la)$ be the irreducible polynomial $\U(\mf{l})$-module corresponding to an $(m|n)$-hook partition $\la$, which has a crystal base with its crystal $\ms{B}(V_{\mf l}(\la))$ being the set of $(m|n)$-hook semistandard tableaux of shape $\la$ \cite{BKK}. Let
$P(\la)=\U(\g)\ot_{\U(\mf{p})}V_{\mf{l}}(\la)$,
where $\U(\mf{p})$ is the subalgebra associated to a parabolic subalgebra $\mf{p}\subset \g$ generated by $\mf{l}$ and the Borel subalgebra with respect to $\Pi$. As a $\Bbbk$-space, $P(\la)=\mc{N}\otimes_{\Bbbk} V_{\mf l}(\la)$, where $\mc{N}$ is the subalgebra of $\U(\g)^-$ generated by the root vectors associated to $\mf{g}/\mf{p}$.
In \cite{JKU24}, it is shown that $P(\la)$ has a unique crystal base $(\ms{L}(P(\la)),\ms{B}(P(\la)))$.
The construction of the crystal base of $\mc{N}=P(0)$ is the main part in the proof of the existence.

Let $Q(\la)$ and $R(\la)$ be the parabolic Verma modules induced from polynomial representations of $\U(\gl_{m|0}\oplus\gl_{0|n})$ and $\U(\gl_{m|0})$, respectively. 
Then we have the commuting diagram as follows: 
\begin{equation}\label{eq:main diagram-0}
\xymatrixcolsep{1.7pc}\xymatrixrowsep{2pc}\xymatrix{
  \U(\g)^- \ \ar@{->}[r]\ar@{->}[d]  &\ R(\la) \ \ar@{->}[r] \ar@{->}[d]  &\ Q(\la) \ \ar@{->}[r] \ar@{->}[d]  &\ P(\la) \ \ar@{->}^{\pi_{V(\la,\ell)}}[r] \ar@{->}[d] &\ V(\la,\ell)  \\
   \mc{N}\ot_{\Bbbk} \U(\mf{l})^- \ar@{->}[r] &  \mc{N}\ot_{\Bbbk} X(\la) \ar@{->}[r] & \mc{N}\ot_{\Bbbk} K(\la)  \ar@{->}[r] & \mc{N}\ot_{\Bbbk} V_{\mf{l}}(\la)
},
\end{equation}
where the horizontal maps in the top rows are canonical $\U(\g)^-$-linear projections, the vertical maps are $\Bbbk$-linear isomorphisms, and $K(\la)$ and $X(\la)$ are the corresponding parabolic Verma modules over $\U(\mf{l})$.

One can show by adopting similar arguments as in \cite{JKU23} that the $\U(\g)^-$, $R(\la)$, and $Q(\la)$ have crystal bases (which are unique in case of the last two ones) such that they are isomorphic to the tensor product of crystal bases of $\mc{N}$ and $\U(\mf{l})$-modules in the bottom line of \eqref{eq:main diagram-0} (as modules over $\U(\mf{l})$ or $\U(\mf{l})^-$), and compatible with the canonical projections (or its dual). Then the $\U(\g)$-crystal structures of $\U(\g)^-$, $R(\la)$, and $Q(\la)$ can also be described explicitly from that of $\mc{N}=P(0)$.

The crystal $\ms{B}(\U(\g)^-)$ of $\U(\g)^-$ can be given in terms of a PBW type basis \cite{CHW,Ya94} together with an explicit description of crystal operators, where the roots for $\mc{N}$ precedes the negative roots of $\mf{l}$.
A description of the associated embeddings from $\ms{B}(P(\la))$ to $\ms{B}(\U(\g)^-)$ follows directly from the ones for $\U(\mf{l})$-modules in \eqref{eq:main diagram-0} (see \cite{JKU23}) by tensoring with ${\rm id}_{\ms{B}(\mc{N})}$.

\subsection{Irreducible quotient, crystal embedding and Burge correspondence}
In order to complete a connection between crystal bases of $V(\la,\ell)$ and $\U(\g)^-$, it remains to consider $\pi_{V(\la,\ell)}$ in \eqref{eq:main diagram-0}.

We first show that $\pi_{V(\la,\ell)}$ preserves the crystal bases. It is done by showing that $(\ms{L}(P(\la)),\ms{B}(P(\la)))$ is a limit of $(\ms{L}(V(\la,\ell)),\ms{B}(V(\la,\ell))$ when $\ell\rightarrow \infty$.
Next, we describe the crystal embedding 
\begin{equation}\label{eq:main embedding-0}
\xymatrixcolsep{2.2pc}\xymatrixrowsep{2pc}\xymatrix{
\ms{B}(V(\la,\ell))\ \ar@{->}[r] & \  \ms{B}(P(\la)),} 
\end{equation} 
which is indeed the main result in this paper. 
It consists of two steps: As a first step, we define an embedding
\begin{equation}\label{eq:main embedding-1}
\xymatrixcolsep{2.2pc}\xymatrixrowsep{2pc}\xymatrix{
\ms{B}(V(\la,\ell))\ \ar@{->}[r] &\ {\displaystyle \left(\lim_{\ell\rightarrow \infty}\ms{B}(V(0,\ell))\right)\times \ms{B}(V_{\mf l}(\la))}}. 
\end{equation}
The $\U(\g)$-crystal $\lim_{\ell\rightarrow \infty}\ms{B}(V(0,\ell))$ can be realized as the set of $(m|n)$-hook semistandard tableaux of arbitrary shape $\delta^\pi$ (rotation of $\de$ by $180^\circ$) with even parity condition on its rows or columns of an $(m|n)$-hook partition $\de$ depending on the type of $\g$. The map \eqref{eq:main embedding-1} is given by applying the combinatorial algorithms in \cite{JK21,K18} (called separation algorithm) to spinor model of $\ms{B}(V(\la,\ell))$, which were introduced  to compute Lusztig data of Kashiwara-Nakashima tableaux of type $BCD$.
As a second step, we establish a $\U(\g)$-crystal isomorphism
\begin{equation}\label{eq:main embedding-2}
\xymatrixcolsep{2.2pc}\xymatrixrowsep{2pc}\xymatrix{
{\displaystyle\lim_{\ell\rightarrow \infty}\ms{B}(V(0,\ell))}\ \ar@{->}[r] &\ \ms{B}(\mc{N}).} 
\end{equation}
It is realized by the (super-analogue of) RSK map on symmetric matrices when $\g=\mf{b}_{m|n}, \mf{c}_{m|n}$, and Burge correspondence of type $D$ \cite{B} when $\g=\mf{d}_{m|n}$.
Hence we get \eqref{eq:main embedding-0} by combining \eqref{eq:main embedding-1} and \eqref{eq:main embedding-2} since $\ms{B}(P(\la))=\ms{B}(\mc{N})\times \ms{B}(V_{\mf l}(\la))$. We remark that \eqref{eq:main embedding-0} is given in \cite{JK21,K18} when $n=0$, while its generalization to the case of $\mf{d}_{m|n}$ with $n>0$ is non-trivial and indeed the most technical part in this paper.

\subsection{Organization}
The paper is organized as follows. In Section \ref{sec:prel}, we recall necessary notations and definitions. In Section \ref{sec:osc and parabolic Verma}, we recall the crystal bases of $V(\la,\ell)$ and $P(\la)$, and prove their compatibility under $\pi_{V(\la,\ell)}$. In Section \ref{sec:review on type A}, we recall the crystal base of $\U(\mf{l})^-$ \cite{JKU23}. In Section \ref{sec:Main1}, we construct the crystal base of $\U(\g)^-$. In Section \ref{sec: Main 2}, we describe \eqref{eq:main embedding-0}. The proof for the crystal theoretic interpretation of Burge correspondence providing \eqref{eq:main embedding-2} for $\g=\mf{d}_{m|n}$ is presented in Section \ref{sec:proof of thm:Burge}.

\section{Preliminary}\label{sec:prel}
\subsection{Notations}\label{subsec:notations}
Throughout the paper, we assume that $m$ and $n$ are positive integers with $m\ge 2$. 
We keep the same notations in \cite[Section 2]{JKU24}, which we briefly recall.

\begin{itemize}
\item[$\bullet$] $\I = \{\, 1, \dots ,m , m+1, \dots, m+n \,\}$ with a usual linear order,

\item[$\bullet$] $\I = \I_{\ov 0} \cup \I_{\ov 1}$, where
$\I_{\ov 0} = \{\, 1, \dots, m \,\}$ and $\I_{\ov 1} = \{ \, m+1,\dots, m+n \,\}$,

\item[$\bullet$] $\e=(\e_a)_{a\in \I}$: a sequence with $\e_a=\varepsilon$ for $a\in \I_{\varepsilon}$ ($\varepsilon={\ov 0},{\ov 1}$),

\item[$\bullet$] $P = \bigoplus_{a\in \I}\Z \de_a \oplus \Z\La_0$: a free abelian group,

\item[$\bullet$] $I=\{\,0,1, \ldots,m,m+1,\dots, m+n-1\}$,

\item[$\bullet$] $\g=\mf{b}_{m|n}, \mf{c}_{m|n}, \mf{d}_{m|n}$: an orthosymplectic Lie superalgebra associated to the Dynkin diagram with the set of simple roots $\Pi=\{\,\alpha_i\,|\,i\in I\,\}\subset P$ as follows: \medskip

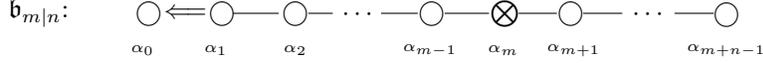
\begin{center}
\hskip -3cm \setlength{\unitlength}{0.16in}
\begin{picture}(24,4)
\put(2,2){\makebox(0,0)[c]{$\mf{b}_{m|n}$:}}
\put(5.6,2){\makebox(0,0)[c]{$\bigcirc$}}
\put(8,2){\makebox(0,0)[c]{$\bigcirc$}}
\put(10.4,2){\makebox(0,0)[c]{$\bigcirc$}}
\put(14.85,2){\makebox(0,0)[c]{$\bigcirc$}}
\put(17.25,2){\makebox(0,0)[c]{$\bigotimes$}}
\put(19.4,2){\makebox(0,0)[c]{$\bigcirc$}}
\put(24.5,1.95){\makebox(0,0)[c]{$\bigcirc$}}
\put(8.35,2){\line(1,0){1.5}}
\put(10.82,2){\line(1,0){0.8}}
\put(13.2,2){\line(1,0){1.2}}
\put(15.28,2){\line(1,0){1.45}}
\put(17.7,2){\line(1,0){1.25}}
\put(19.81,2){\line(1,0){1.28}}
\put(22.8,2){\line(1,0){1.28}}
\put(6.8,2){\makebox(0,0)[c]{$\Longleftarrow$}}
\put(12.5,1.95){\makebox(0,0)[c]{$\cdots$}}
\put(22,1.95){\makebox(0,0)[c]{$\cdots$}}
\put(5.4,0.8){\makebox(0,0)[c]{\tiny $\alpha_{0}$}}
\put(7.8,0.8){\makebox(0,0)[c]{\tiny $\alpha_{1}$}}
\put(10.4,0.8){\makebox(0,0)[c]{\tiny $\alpha_{2}$}}
\put(14.8,0.8){\makebox(0,0)[c]{\tiny $\alpha_{m-1}$}}
\put(17.2,0.8){\makebox(0,0)[c]{\tiny $\alpha_{m}$}}
\put(19.5,0.8){\makebox(0,0)[c]{\tiny $\alpha_{m+1}$}}
\put(24.5,0.8){\makebox(0,0)[c]{\tiny $\alpha_{m+n-1}$}}
\end{picture}
\end{center}
\begin{equation*}
\alpha_i=
\begin{cases}
-\de_{1} & \text{if $i=0$},\\
\de_{i}-\de_{i+1} & \text{if $1 \le i \le m+n-1$},
\end{cases}
\end{equation*}

\begin{center}
\hskip -3cm \setlength{\unitlength}{0.16in}
\begin{picture}(24,4)
\put(2,2){\makebox(0,0)[c]{$\mf{c}_{m|n}$:}}
\put(5.6,2){\makebox(0,0)[c]{$\bigcirc$}}
\put(8,2){\makebox(0,0)[c]{$\bigcirc$}}
\put(10.4,2){\makebox(0,0)[c]{$\bigcirc$}}
\put(14.85,2){\makebox(0,0)[c]{$\bigcirc$}}
\put(17.25,2){\makebox(0,0)[c]{$\bigotimes$}}
\put(19.4,2){\makebox(0,0)[c]{$\bigcirc$}}
\put(24.5,1.95){\makebox(0,0)[c]{$\bigcirc$}}
\put(8.35,2){\line(1,0){1.5}}
\put(10.82,2){\line(1,0){0.8}}
\put(13.2,2){\line(1,0){1.2}}
\put(15.28,2){\line(1,0){1.45}}
\put(17.7,2){\line(1,0){1.25}}
\put(19.81,2){\line(1,0){1.28}}
\put(22.8,2){\line(1,0){1.28}}
\put(6.8,2){\makebox(0,0)[c]{$\Longrightarrow$}}
\put(12.5,1.95){\makebox(0,0)[c]{$\cdots$}}
\put(22,1.95){\makebox(0,0)[c]{$\cdots$}}
\put(5.4,0.8){\makebox(0,0)[c]{\tiny $\alpha_{0}$}}
\put(7.8,0.8){\makebox(0,0)[c]{\tiny $\alpha_{1}$}}
\put(10.4,0.8){\makebox(0,0)[c]{\tiny $\alpha_{2}$}}
\put(14.8,0.8){\makebox(0,0)[c]{\tiny $\alpha_{m-1}$}}
\put(17.2,0.8){\makebox(0,0)[c]{\tiny $\alpha_{m}$}}
\put(19.5,0.8){\makebox(0,0)[c]{\tiny $\alpha_{m+1}$}}
\put(24.5,0.8){\makebox(0,0)[c]{\tiny $\alpha_{m+n-1}$}}
\end{picture}
\end{center}
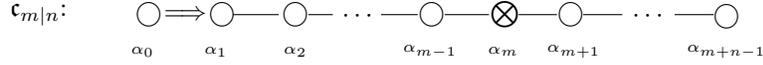
\begin{equation*}
\alpha_i=
\begin{cases}
-2\de_{1} & \text{if $i=0$},\\
\delta_i-\delta_{i+1} & \text{if $1\le i \le m+n-1$},
\end{cases}
\end{equation*}

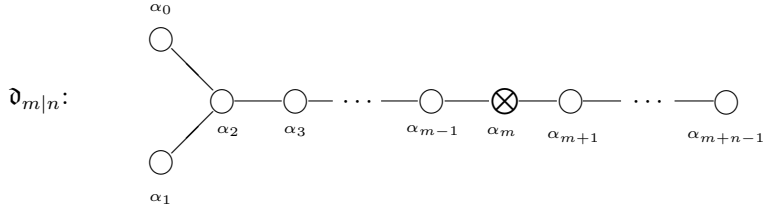
\begin{center}
\hskip -3cm \setlength{\unitlength}{0.16in} \medskip
\begin{picture}(24,5.8)
\put(2,2){\makebox(0,0)[c]{$\mf{d}_{m|n}$:}}
\put(6,0){\makebox(0,0)[c]{$\bigcirc$}}
\put(6,4){\makebox(0,0)[c]{$\bigcirc$}}
\put(8,2){\makebox(0,0)[c]{$\bigcirc$}}
\put(10.4,2){\makebox(0,0)[c]{$\bigcirc$}}
\put(14.85,2){\makebox(0,0)[c]{$\bigcirc$}}
\put(17.25,2){\makebox(0,0)[c]{$\bigotimes$}}
\put(19.4,2){\makebox(0,0)[c]{$\bigcirc$}}
\put(24.5,1.95){\makebox(0,0)[c]{$\bigcirc$}}
\put(6.35,0.3){\line(1,1){1.35}} 
\put(6.35,3.7){\line(1,-1){1.35}}
\put(8.4,2){\line(1,0){1.55}} 
\put(10.82,2){\line(1,0){0.8}}
\put(13.2,2){\line(1,0){1.2}} 
\put(15.28,2){\line(1,0){1.45}}
\put(17.7,2){\line(1,0){1.25}} 
\put(19.8,2){\line(1,0){1.25}}
\put(22.8,2){\line(1,0){1.28}}
\put(12.5,1.95){\makebox(0,0)[c]{$\cdots$}}
\put(22,1.95){\makebox(0,0)[c]{$\cdots$}}
\put(6,5){\makebox(0,0)[c]{\tiny $\alpha_{0}$}}
\put(6,-1.2){\makebox(0,0)[c]{\tiny $\alpha_{1}$}}
\put(8.2,1){\makebox(0,0)[c]{\tiny $\alpha_{2}$}}
\put(10.4,1){\makebox(0,0)[c]{\tiny $\alpha_{3}$}}
\put(14.9,1){\makebox(0,0)[c]{\tiny $\alpha_{m-1}$}}
\put(17.15,1){\makebox(0,0)[c]{\tiny $\alpha_m$}}
\put(19.5,0.8){\makebox(0,0)[c]{\tiny $\alpha_{m+1}$}}
\put(24.5,0.8){\makebox(0,0)[c]{\tiny $\alpha_{m+n-1}$}}
\end{picture}\vskip 8mm
\end{center}
\begin{equation*}
\alpha_i=
\begin{cases}
-\de_{1}-\de_{2} & \text{if $i=0$},\\
\delta_i-\delta_{i+1} & \text{if $1\le i \le m+n-1$}.
\end{cases}
\end{equation*}
\medskip

\item[$\bullet$] $(\, \cdot\, |\, \cdot\, )$ : a symmetric bilinear form on $P$ satisfying
$(\de_a|\de_b)=(-1)^{\e_a}r_{\g}\de_{ab}$,
$(\de_a|\La_0)=-\frac{r_{\g}}{2}$ $(a,b\in \I)$
with {$r_{\mf{b}_{m|n}}=2$, $r_{\mf{c}_{m|n}}=r_{\mf{d}_{m|n}}=1$} (note that $(\La_0|\alpha_i)=\de_{0i}$ for $i\in I$),

\item[$\bullet$] $\Phi^+$ : the set of reduced positive roots,

\item[$\bullet$] $\Phi^+_{\ov{0}}$ (resp.~$\Phi^+_{\ov{1}}$) : the set of even (resp.~odd) roots in $\Phi^+$,

\item[$\bullet$] $\mf{l}=\mf{l}_{m|n}$ : the subalgebra of $\g$ corresponding to $\Pi\setminus\{\alpha_0\}$ isomorphic to $\gl_{m|n}$,

\item[$\bullet$] $\mf{l}_{m|0}$, $\mf{l}_{0|n}$: the subalgebras of $\mf{l}_{\ov 0}$ corresponding to $\gl_m$ and $\gl_n$, respectively,

\item[$\bullet$] $\Phi^+(\mf{l})$ : the set of positive roots of $\mf{l}$,

\item[$\bullet$] $\mf{u}$ : the subalgebra of $\g$ generated by the root vectors corresponding to $\Phi^+\setminus \Phi^+(\mf{l})$,

\item[$\bullet$] $\Phi^+(\mf{u})= \Phi^+\setminus \Phi^+(\mf{l})$.

\end{itemize}

\subsection{Quantum orthosymplectic superalgebra}
We also recall the following \cite{JKU24}:
\begin{itemize}
\item[$\bullet$] $\Bbbk=\Q(v)$ and $\Bbbk^{\times}=\Bbbk \setminus \left\{0\right\}$ where $v=q$ if $\g=\mf{b}, \mf{c} $ and $v=q^{\frac{1}{2}}$ if $\g=\mf{d} $ with $ q^{\frac{1}{2}}$ an indeterminate,

\item[$\bullet$] ${\ttq}_a=\sa q^{\sa r_{\g}}$ $(a\in \I)$, 

\item[$\bullet$] $q_i = q^{d_i}$ $(i\in I\setminus\{m\})$, where $d_i=|(\alpha_i|\alpha_i)|/2$, and $q_m = q^{r_{\mf g}}$,

\item[$\bullet$] $\bq(\mu,\nu) = v^{\sum_{a\in \I} \ell' \mu_a + \ell \nu_a  } \prod_{a\in \I}{\tt q}_a^{\mu_a\nu_a}$ for $\mu=\sum_{a\in\I}\mu_a\de_a +\ell \Lambda_0 $ and $\nu=\sum_{a\in\I}\nu_a\de_a + \ell' \Lambda_0$ (note that $\bq(\Lambda_0,\alpha_i)=q_0^{\delta_{0 i}}$),

\item[$\bullet$] For $s\in\Z_{\ge 0}$ and $z\in q^{\Z_{>0}}$,
$[s]_z=\frac{z^s-z^{-s}}{z-z^{-1}}$, and 
$[s]_z!=[s]_z [s-1]_z \cdots [1]_z$  $(s\geq 1)$ with $[0]_z!=1$, 
where we simply write $[s]$ when $z=q$, and $[s]_i=[s]_{q_i}$ when $z=q_i$ ($i\in I$),

\item[$\bullet$] ${\U}(\g)$: the associative $\Bbbk$-algebra with $1$ 
generated by $k_\mu, e_i, f_i$ for $\mu\in P$ and $i\in I$ 
satisfying
{\allowdisplaybreaks
\begin{gather*}
k_{\mu}=1 \quad (\mu=0), \quad k_{\nu +\nu'}=k_{\nu}k_{\nu'} \quad (\nu, \nu' \in P),\label{eq:Weyl-rel-1} \\ 
k_\mu e_i k_{-\mu}=\bq(\mu,\alpha_i)e_i,\quad 
k_\mu f_i k_{-\mu}=\bq(\mu,\alpha_i)^{-1}f_i\quad (i\in I, \mu\in P), \label{eq:Weyl-rel-2} \\ 
e_if_j - f_je_i =\delta_{ij}\frac{k_{i} - k^{-1}_{i}}{q_i-q_i^{-1}}\quad (i,j\in I),\label{eq:Weyl-rel-3}\\
e_m^2= f_m^2 =0,\label{eq:Weyl-rel-4}
\end{gather*}
where $k_i=k_{\alpha_i}$ for $i\in I$, and 
\begin{gather*}
  x_i x_j -  x_j x_i =0
 \quad \text{($i,j \in I$ and $(\alpha_i|\alpha_j) = 0$)},\\
x_i^2 x_j- (-1)^{\e_i}[2]_i x_i x_j x_i + x_j x_i^2= 0
\quad \text{($i\in I\setminus\{0, m\}$, $j\in I\setminus\{0\}$, $i\neq j$ and $(\alpha_i|\alpha_j)\neq 0$)},    \\
\begin{array}{ll}
		  x_0^3 x_1- [3]_0 x_0^2 x_1 x_0 + [3]_0 x_0 x_1 x_0^2 - x_1 x_0^3 = 0\\ 
		  x_1^2 x_0- [2]_1 x_1 x_0 x_1 + x_0 x_1^2= 0
\end{array}
\quad (\g = \mf{b}_{m|n}),    \\
\begin{array}{ll}
		x_0^2 x_1- [2]_0 x_0 x_1 x_0 + x_1 x_0^2= 0 \\
		x_1^3 x_0- [3]_1 x_1^2 x_0 x_1 + [3]_1 x_1 x_0 x_1^2 - x_0 x_1^3 = 0
\end{array}
\quad (\g = \mf{c}_{m|n}), \\
\begin{array}{ll}
		x_0^2 x_2- [2]_0 x_0 x_2 x_0 + x_2 x_0^2= 0 \\
		x_2^2 x_0- [2]_2 x_2 x_0 x_2 + x_0 x_2^2= 0
\end{array}
	\quad (\g = \mf{d}_{m|n}\ \text{with $m\ge 3$}),\\
\begin{array}{ll}
  x_{m}x_{m'}x_{m}x_{m+1}  
- x_{m}x_{m+1}x_{m}x_{m'} 
+ x_{m+1}x_{m}x_{m'}x_{m} \\
\qquad \qquad \qquad \qquad \qquad - x_{m'}x_{m}x_{m+1}x_{m} 
+ [2]_{m} x_{m}x_{m'}x_{m+1}x_{m} =0 \\ 
\hskip 6cm (m'=0,1 \text{ if $\g={\mf d}_{2|n}$}, \text{ and } m'=m-1 \text{ otherwise} ),
\end{array}
\end{gather*}}
where $x=e,f$,

\item[$\bullet$] $\Delta$: the comultiplication on $\U(\g)$ given by 
$\Delta(k_\mu)=k_\mu\otimes k_\mu$, $\Delta(e_i)= 1\ot e_i + e_i\ot k_i^{-1}$, $\Delta(f_i)= f_i\ot 1 + k_i\ot f_i$ for $\mu\in P$ and $i\in I$,

\item[$\bullet$] $e'_i, e''_i$: a $\Bbbk$-linear map on $\U(\g)^-$ defined by 
\begin{equation}\label{eq: q-derivation}
	\begin{split}
 e'_i(f_j)&=\de_{ij},\quad  
 e'_i(uv)=e'_i(u)v + \bq{(\alpha_i,|u|)}ue'_i(v), \\
 e''_i(f_j)&=\de_{ij},\quad  
 e''_i(uv)=e''_i(u)v + \bq{(\alpha_i,|u|)}^{-1}ue''_i(v), 
	\end{split}
\end{equation}
for $i,j\in I$ and $u,v\in \U(\g)^-$ with $u$ homogeneous.

\end{itemize}
We also consider the following subalgebras of $\U(\g)$:
\begin{itemize}
\item[$\bullet$] $\U(\g)^+= \langle\, e_i \,|\,i\in I\,\rangle$, $\U(\g)^- = \langle\, f_i \,|\,i\in I\,\rangle$, $\U(\g)^0 = \langle\, k_\mu \,|\,\mu \in P\,\rangle$, 

\item[$\bullet$] $\U(\mf{l})=\langle\, k_\mu, e_i, f_i \,|\,\mu\in P,\  i\in I\setminus\{0\}\,\rangle$,

\item[$\bullet$] $\U(\mf{l}_{\ov 0})=\langle\, k_\mu, e_i, f_i \,|\,\mu\in P,\  i\in I\setminus\{0,m\}\,\rangle$,

\item[$\bullet$] $\U(\mf{l}_{m|0})=\langle\, k_\mu, e_i, f_i \,|\,\mu\in P,\  i=1,\dots,m-1 \,\rangle$,

\item[$\bullet$] $\U(\mf{l}_{0|n})=\langle\, k_\mu, e_i, f_i \,|\,\mu\in P,\  i=m+1,\dots,m+n-1 \,\rangle$,

\item[$\bullet$] $\U(\mf{l}_\ast)^\pm = \U(\mf{l}_\ast)\cap \U(\g)^\pm$, where $\mf{l}_\ast = \mf{l}, \mf{l}_{\ov 0}, \mf{l}_{m|0}, \mf{l}_{0|n}$.

\end{itemize}
Note that $\U(\mf{l}_{\ov 0})$ is isomorphic to $U_q(\gl_m)\ot U_{-q^{-1}}(\gl_n)$ as a $\Bbbk$-algebra up to Cartan part, where $U_q(\gl_m)$ is the usual quantized enveloping algebra associated to $\gl_m$.  
Recall that $\U(\g)^\pm$ is naturally graded by $Q^+:=\sum_{i\in I}\Z_{\ge 0}\alpha_i$ and $Q^-:=-Q^+$.

Let $\mf{a}$ denote either $\g$ or $\mf{l}_\ast$ above. 
For a $\U(\mf a)$-module $V$ and $\mu\in P$, let 
$V_\mu = \{\,u\in V\,|\,k_{\nu} u= \bq(\mu,\nu) u \ \ (\nu\in P) \,\}$
be the $\mu$-weight space of $V$.
For $u\in V_\mu\setminus\{0\}$, we call $u$ a weight vector with weight $\mu$ and put ${\rm wt}(u)=\mu$. 
Denote by ${\rm wt}(V)$ the set of weights of $V$. 
For $\La\in P$, let $V_{\mf a}(\La)$ be the irreducible highest weight $\U(\mf{a})$-module generated by a highest weight vector $v^{\mf a}_\La$ or simply $v_\La$ of weight $\La$.

We also remark that the representation theory of $\U(\g)$ is equivalent to that of the quantum superalgebra associated to $\g$ \cite{Ya94} (see \cite[Section 2.2]{JKU24} for more details).

\subsection{Root vectors}\label{subsec:root vectors}
Let $\bff_\beta$ ($\beta\in \Phi^+$) be a root vector in $\U(\g)^-$, which is defined in \cite[Section 3.1]{JKU24} following \cite{CHW,Le}.
Let $\prec$ be the linear order on $\Phi^+$ induced from that on the set of good Lyndon words $\mc{GL}(\g)$ under the bijection with $\Phi^+$. Then we have a PBW type basis $B$ of $\U(\g)^-$ with respect to $\prec$.
Note that $\alpha\prec \beta$ for $\alpha\in \Phi^+(\mf{u})$ and $\beta\in \Phi^+(\mf{l})$. For $\beta\in \Phi^+(\mf{u})$, we simply write  
\begin{equation*}
 \beta= 
\begin{cases}
 (i,i) & \text{if $\beta = -2\de_i/r_{\g}$ for $1\le i\le m+n$},\\
 (i,j) & \text{if $\beta = -\de_i-\de_j$ for $1\le i<j\le m+n$}.
\end{cases}
\end{equation*}
Then
\begin{equation}\label{eq:nilradical roots}
 \Phi^+(\mf{u})=
\begin{cases}
 \{\,(i,i)\,|\,1\le i\le m+n\,\} \sqcup \{\,(i,j)\,|\,1\le i<j\le m+n\,\} & \text{if $\g=\mf{b}_{m|n}$},\\ 
  \{\,(i,i)\,|\,1\le i\le m\,\} \sqcup \{\,(i,j)\,|\,1\le i<j\le m+n\,\} & \text{if $\g=\mf{c}_{m|n}$},\\
   \{\,(i,i)\,|\,m+1\le i\le m+n\,\} \sqcup \{\,(i,j)\,|\,1\le i<j\le m+n\,\} & \text{if $\g=\mf{d}_{m|n}$},
\end{cases}
\end{equation}
and for $(i_1,j_1), (i_2,j_2)\in \Phi^+(\mf{u})$,
	\begin{equation*}
		(i_1,j_1) \prec (i_2,j_2) \ \Longleftrightarrow \
		\begin{cases*}
			i_1<i_2 \text{ or } (i_1=i_2 \text{ and } j_1<j_2) & \text{if $\g=\mf{b}_{m|n}$,} \\
			j_1<j_2 \text{ or } (j_1=j_2 \text{ and } i_1<i_2) & \text{if $\g=\mf{c}_{m|n}, \mf{d}_{m|n}$.}
		\end{cases*}
	\end{equation*}

Let $\U(\mf{u}^-)$ be the $\Bbbk$-space spanned by the vectors in $B$ where the product is over $\beta\in \Phi^+(\mf{u})$.
It is a $\Bbbk$-subalgebra of $\U(\g)^-$ generated by $\bff_\beta$ ($\beta\in \Phi^+(\mf{u})$), which we denote by $\mc{N}$.
As a $\Bbbk$-space 
\begin{equation*}\label{eq:Levi decomp}
 \U(\g)^-\cong \U(\mf{u}^-)\ot \U(\mf{l})^-,
\end{equation*}
where $\ot=\ot_\Bbbk$.

\subsection{Irreducible oscillator modules $V(\la,\ell)$}\label{subsec:osc modules}
Let us recall a semisimple tensor category of $\U(\g)$-modules studied in \cite{K15}.
Let  

\begin{itemize}
 \item[$\bullet$] ${\tt P}=\sum_{a\in \I}\Z\delta _a \subset P$,

 \item[$\bullet$] ${\tt P}^+=\{\,\la=\sum_{i\in \I}\la_i\de_i \in {\tt P}\ |\ \lambda_{1}\geq\ldots\geq\lambda_{m},\ \ \lambda_{m+1}\geq \ldots\geq \lambda_{m+n}\,\}$,

 \item[$\bullet$] $\lambda_{+}=\sum_{i\in \I_0}\la_{i}\de_i$, $\lambda_{-}=\sum_{j\in \I_1}\la_{j}\de_j$ for $\lambda\in {\tt P}$,

 \item[$\bullet$] ${\tt P}_{\geq0}=\sum_{a\in \I}\Z_{\ge 0}\delta _a$.

\end{itemize}

Let $\cO(\mf l)_{\geq0}$ be the category of $\U(\mf{l})$-modules $V$ such that
$V=\bigoplus_{\mu\in {\tt P}_{\ge 0}}V_\mu$ with $\dim V_\mu < \infty$.
Let $\cP$ be the set of all partitions $\la=(\la_i)_{i\ge 1}$, and let $\cP_{m|n}$ be the set of $\la\in\cP$ such that $\la_{m+1}\leq n$.
For $\la\in \cP_{m|n}$, let  
\begin{equation}\label{eq:highest weight correspondence}
\La_\la = \la_1\de_1 +\cdots +\la_m\de_m + \mu_1\de_{m+1}+\cdots+\mu_n\de_{m+n} \in {\tt P}_{\geq 0},
\end{equation}
where $(\mu_1,\ldots,\mu_n)$ is the conjugate of the partition $(\la_{m+1},\la_{m+2},\ldots)$. We often identify $\la\in \cP_{m|n}$ with $\La_\la$ if there is no confusion.
We simply write $V_{\mf{l}}(\la)=V_{\mf{l}}(\La_{\la})$ and $v_\la=v_{\La_\la}$.
Then any irreducible $\U(\mf{l})$-module in $\cO(\mf l)_{\geq0}$ is isomorphic to $V_{\mf{l}}(\la)$ for some $\la\in \cP_{m|n}$ \cite[Propositions 3.4 and 4.5]{BKK} (see also \cite[Remark 4.2]{JKU24}). Moreover, $\cO(\mf l)_{\geq0}$ is semisimple \cite[Proposition 4.3]{JKU24}.

\begin{df}\label{category O^int_q(m|n)}{\rm
Let $\mc{O}(\g)_{\ge 0}$ be  the category of $\U(\mf{g})$-modules $V$ satisfying
\begin{itemize}
\item[(1)] $V=\bigoplus_{\gamma\in P}V_\gamma$ and $\dim V_\gamma <\infty$ for $\gamma\in P$,

\item[(2)] ${\rm wt}(V)\subset \bigcup_{i=1}^r\left(\ell_i\Lambda_{0}+ {\tt P}_{\ge 0}  \right)$ for some $r\geq 1$ and  $\ell_i\in \Z_{\geq 0}$, 

\item[(3)] $f_{0}$ acts locally nilpotently on $V$.
\end{itemize}
}
\end{df}
We may understand that $V\in \mc{O}(\g)_{\ge 0}$ belongs to $\mc{O}(\mf l)_{\ge 0}$ as a $\U(\mf{l})$-module if we ignore $\La_0$, and hence is a direct sum of $V_{\mf l}(\mu)$ for $\mu\in \cP_{m|n}$ with finite multiplicity.
Let $\cP(\mf{g})$ be given by
\begin{equation*}\label{eq:P(g)}
\begin{split}
\cP(\mf{b}_{m|n})&=\{\,(\lambda,\ell)\in \cP_{m|n}\times\Z_{> 0}\,|\,\ell-2\lambda_1\in \Z_{\geq 0}\,\},\\
\cP(\mf{c}_{m|n})&=\{\,(\lambda,\ell)\in \cP_{m|n}\times\Z_{> 0}\,|\,\ell-\lambda_1\in \Z_{\geq 0}\,\},\\
\cP(\mf{d}_{m|n})&=\{\,(\lambda,\ell)\in \cP_{m|n}\times\Z_{> 0}\,|\,\ell-\lambda_1-\lambda_2\in \Z_{\geq 0}\,\},
\end{split}
\end{equation*}
and for $(\la,\ell)\in \cP(\mf{g})$, let
\begin{equation*}\label{eq:hw for osc}
 \La(\la,\ell) = \ell\La_0 + \La_\la. 
\end{equation*}

\begin{thm}{\em (\cite[Theorem 3.8]{K15}) } \label{thm:category O_int}
$\mc{O}(\g)_{\ge 0}$ is a semisimple tensor category with non-trivial irreducible objects $V(\La(\la,\ell))$ for $(\la,\ell)\in \cP(\g)$. 
\end{thm}
For simplicity, let us write $V(\la,\ell)=V(\La(\la,\ell))$ and $v_{(\la,\ell)}=v_{\La(\la,\ell)}$.
For $V\in \mc{O}(\g)_{\ge 0}$, we have
\begin{equation*}
 V \cong \bigoplus_{(\la,\ell)\in \cP(\g)}V(\la,\ell)^{\oplus m(\la,\ell)},
\end{equation*}
for some $m(\la,\ell)\in \Z_{\ge 0}$.

\subsection{Parabolic Verma modules}

Let $\U(\mf{p})$ be the subalgebra of $\U(\g)$ generated by $\U(\mf{l})$ and $\U(\g)^{+}$.
For $\la\in \cP_{m|n}$, let
\begin{equation*}\label{eq:P(la)}
 P(\la)=\U(\g)\ot_{\U(\mf{p})}V_{\mf{l}}(\la),
\end{equation*} 
where we regard $V_{\mf{l}}(\la)$ as a $\U(\mf{p})$-module with $e_0 v_\la=0$.
We have an isomorphism of $\Bbbk$-spaces
\begin{equation}\label{eq:P(la)-2}
\mc{N}\ot V_{\mf l}(\la)\cong  P(\la),
\end{equation}
sending $u\ot v$ to itself since 
$\U(\g)\cong \U(\mf{u}^-)\ot \U(\mf{l})^-\ot \U(\g)^{\ge 0}\cong \mc{N} \ot \U(\mf{p})$,
where $\U(\g)^{\ge 0} = \U(\g)^+\U(\g)^0$. It is a highest weight $\U(\mf{g})$-module with highest weight $\Lambda_\la$ \eqref{eq:highest weight correspondence}. 

Let $\U(\mf{q})$ be the subalgebra of $\U(\g)$ generated by $\U(\mf{l}_{\ov 0})$ and $\U(\g)^{+}$.
For $\la\in {\tt P}^+$, let 
\begin{equation}\label{eq:Q}
 Q(\la)=\U(\g)\ot_{\U(\mf{q})}V_{\mf{l}_{\ov 0}}(\la),
\end{equation}
where we regard $V_{\mf{l}_{\ov 0}}(\la)$ as a $\U(\mf{q})$-module with $e_0 v_{\la}=e_m v_{\la}=0$.

Finally, let $\U(\mf{r})$ be the subalgebra of $\U(\mf{l})$ generated by $\U(\mf{l}_{0|n})$ and $\U(\g)^{+}$.
For $\la\in {\tt P}$ with $\la_-\in {\tt P}^+$, let  
\begin{equation}\label{eq:R}
 R(\la)=\U(\g)\ot_{\U(\mf{r})}V_{\mf{l}_{0|n}}(\la),
\end{equation}
where we regard $V_{\mf{l}_{0|n}}(\la)$ as a $\U(\mf{r})$-module with $e_i v_{\la}=0$ for $i=0,\dots,m$.

For $(\la,\ell)\in \cP(\g)$, we have natural surjective $\U(\g)^-$-linear maps
\begin{equation*}
 \xymatrixcolsep{3pc}\xymatrixrowsep{0pc}\xymatrix{
  \U(\g)^- \ \ar@{->}[r]^{\pi_{R(\la)}} &\ R(\la) \ \ar@{->}[r]^{\pi_{Q(\la)}} &\ Q(\la) \ \ar@{->}[r]^{\pi_{P(\la)}} &\ P(\la) \ \ar@{->}[r]^{\pi_{V(\la,\ell)}} &\ V(\la,\ell) },
\end{equation*}
where each map sends a highest weight vector to a highest weight vector. 

\section{{Crystal bases of $V(\la,\ell)$ and $P(\la)$}}\label{sec:osc and parabolic Verma}
\subsection{Crystal bases of $V(\la,\ell)$ and $P(\la)$}
Let us briefly recall the results on crystal bases of $V(\la,\ell)$ \cite{K15,K16} and $P(\la)$ \cite{JKU24}.

\subsubsection{}\label{subsubsec: crystal base of osc}
Let $V$ be a $\U(\mf{g})$-module in $\mc{O}(\g)_{\ge 0}$.  
For a weight vector $u\in V$ and $i\in I$, we define $\tilde{e}_i u$ and $\tilde{f}_i u$ as follows:\smallskip

\begin{itemize}
 \item[(1)] Suppose that $i\neq m$. We have $u=\sum_{k \geq 0} f_i^{(k)}u_k$, where $f_i^{(k)}=f_i^k/[k]_i!$ and $e_iu_k=0$ for $k\geq 0$. We define
\begin{gather*} 
\te_iu=\sum_{k\geq1}f_i^{(k-1)}u_k,\quad \tf_iu=\sum_{k\geq0}f_i^{(k+1)}u_k \quad (i<m),\label{eq:crystal operator <m}\\
\tilde{e}_iu=\sum_{k\geq1}q_i^{-l_k+2k-1}f_i^{(k-1)}u_k,\quad 
\tilde{f}_iu=\sum_{k\geq 0}q_i^{l_k-2k-1}f_i^{(k+1)}u_k\quad (i>m),\label{eq:crystal operator >m}
\end{gather*}
where $l_k= - ({\rm wt}(u_k)|\alpha_i)/r_{\mf{g}}$. Here we follow the convention in \cite[Section 4.2]{JKU23}.

\item[(2)] Suppose that $i=m$. We define
\begin{equation*}\label{eq:crystal operator m}
\tilde{e}_m u =\eta(f_m) u =q_m^{-1}k_me_m u,\quad \tilde{f}_m u=f_m u,
\end{equation*}
where $\eta$ is the anti-involution on $\U(\mf{g})$ defined by
$\eta(k_\mu)=k_\mu$, $\eta(e_i)=q_if_ik^{-1}_i$, and $\eta(f_i)=q^{-1}_ik_ie_i$ for $\mu \in  P$ and $i\in I$.
\end{itemize}

Let $A_0$ be the subring of rational functions $f(q)\in \Bbbk$ regular at $q=0$. 
We define a crystal base of $V$ to be a pair $(L,B)$ in a usual way with respect to $\te_i$ and $\tf_i$ ($i\in I$), where $L$ is an $A_0$-lattice of $V$, and $B$ is a signed basis of $L/qL$, that is, $B=\mathbf{B}\cup -\mathbf{B}$, where $\mathbf{B}$ is a $\Q$-basis of $L/qL$. We call the set $B/\{\pm 1\}$, which has a colored oriented graph structure, a $\U(\g)$-crystal of $V$ or simply crystal of $V$. We identify $B$ with $B/\{\pm 1\}$ when we consider $B$ as a crystal, if there is no confusion.
The following is obtained in \cite[Theorem 8.8]{K15} for $\g=\mf{b}_{m|n}, \mf{c}_{m|n}$ and \cite[Theorem 5.7]{K16} for $\g=\mf{d}_{m|n}$.  

\begin{thm}\label{thm: crystal base of osc}
 For $(\la,\ell)\in \cP(\g)$, $V(\la,\ell)$ has a crystal base $(\ms{L}(V(\la,\ell)), \ms{B}(V(\la,\ell)))$, where the crystal $\ms{B}(V(\la,\ell))$ is connected.
\end{thm}

Since $\ms{B}(V(\la,\ell))$ is connected, we have a unique crystal base of $V(\la,\ell)$ up to scalar multiplication, that is, 
\begin{equation}\label{eq: L,B by operators}
\begin{split}
\ms{L}(V(\la,\ell))&=\sum_{r \ge 0, i_1,\dots,i_r \in I}A_0\td{x}_{i_1}\dots\td{x}_{i_r}v_{(\la,\ell)} \\ 
\ms{B}(V(\la,\ell))&=\left\{\,\pm\td{x}_{i_1}\dots\td{x}_{i_r}v_{(\la,\ell)} \!\!\!\pmod{q\ms{L}(V(\la,\ell))}\,|\,r\ge 0, i_1,\dots,i_r\in I\,\right\}\setminus\{0\}, 
\end{split}
\end{equation}
where $x = e,f$ for each $i_k$. Note that \eqref{eq: L,B by operators} holds in general for a highest weight module, which has a crystal base (with respect to $\te_i$ and $\tf_i$ defined suitably) such that the associated crystal is connected.
A combinatorial model for $\ms{B}(V(\la,\ell))$ is also given in \cite{K15,K16}, which will be reviewed in Section \ref{subsec:spinor model}.

Regarding $V\in \mc{O}(\g)_{\ge 0}$ as a $\U(\mf{l})$-module, we may consider a crystal base of $V\in \mc{O}(\mf l)_{\ge 0}$ with respect to $\te_i$ and $\tf_i$ ($i\in I\setminus\{0\}$). For $\la\in \cP_{m|n}$,  $V_{\mf{l}}(\la)$ has a crystal base $(\ms{L}(V_{\mf l}(\la)), \ms{B}(V_{\mf l}(\la)))$ \cite{BKK}, where the connected crystal $\ms{B}(V_{\mf l}(\la))$ can be realized as the set $SST_{m|n}(\la)$ of $(m|n)$-hook semistandard tableaux of shape $\la$.
Hence $(\ms{L}(V(\la,\ell)), \ms{B}(V(\la,\ell)))$ is a direct sum of crystal bases of $V_{\mf{l}}(\mu)$'s for $\mu\in \cP_{m|n}$ as a $\U(\mf{l})$-module.

\subsubsection{}\label{subsubsec: crystal base of parabolic P}
Let $\la\in \cP_{m|n}$ be given.
For a weight vector $u\in P(\la)$, we define 
\begin{equation}\label{eq: crystal operator for 0}
 \te_0 u = \sum_{k\ge 1}f_0^{(k-1)}u_k,\quad \tf_0 u = \sum_{k\ge 0}f_0^{(k+1)}u_k,
\end{equation}
where $u=\sum_{k\ge 0}f_0^{(k)}u_k$ for $e'_0 ( u_k ) = 0$.
Then we define a crystal base of $P(\la)$ in the same way as in Section \ref{subsubsec: crystal base of osc}, where $\te_0$ and $\tf_0$ is replaced by \eqref{eq: crystal operator for 0}.
We may define $\te_0$ and $\tf_0$ on $\U(\g)^-$ in the same way, and we have $(\td{x}_0P)(1\ot v_\la) = \td{x}_0(P(1\ot v_\la))$ for $P\in \U(\g)^-$ (cf.~\cite[Section 4.4]{JKU24}).

\begin{thm}{\em (\cite[Section 4.4]{JKU24}) }\label{thm: crystal base of P}
 For $\la \in \cP_{m|n}$, $P(\la)$ has a crystal base $(\ms{L}(P(\la)), \ms{B}(P(\la)))$, where the crystal $\ms{B}(P(\la))$ is connected.
\end{thm}

Let ${\rm ad}_q: \U(\mf{g})  \longrightarrow {\rm End}_{\Bbbk}(\U(\mf{g}))$ be the adjoint representation of $\U(\g)$ with respect to $\Delta$,
where for $\mu \in P$, $i \in I$, and $u \in \U(\g)$,
\begin{equation} \label{eq: quantum adjoint}
\begin{split}
	& {\rm ad}_q(k_\mu)(u) = k_\mu u k_\mu^{-1}, \ \
	{\rm ad}_q(e_i)(u) = \left( e_i u - u e_i \right) k_i, \ \
	{\rm ad}_q(f_i)(u) = f_i u - k_i u k_i^{-1} f_i.
\end{split}
\end{equation}
Let us write ${\rm ad}_q(x)(u)=x\cdot u$ for simplicity.
Then the subalgebra $\mc{N}$ is a $\U(\mf{l})$-submodule of $\U(\mf{g})$ with respect to ${\rm ad}_q$ such that
$x\cdot m_{\mc N}(u_1\ot u_2) = m_{\mc N}(\Delta(x)\cdot(u_1\ot u_2))$
for $x\in \U(\mf{l})$ and $u_1\ot u_2\in \mc{N}\ot \mc{N}$, where $m_{\mc N}$ denotes the multiplication on $\mc{N}$ \cite[Proposition 3.16]{JKU24}.
We have $\mc{N}\in\mc{O}(\mf l)_{\ge 0}$ as a $\U(\mf{l})$-module, and 
the map \eqref{eq:P(la)-2} is an isomorphism of $\U(\mf{l})$-modules.

For $\beta\in \Phi^+(\mf{u})$ and $k\in \Z_{\ge 0}$, let
\begin{equation*}\label{eq:divided power of root vector}
\bfF_\beta^{(k)} = 
\begin{cases}
 \bff_\beta^{(k)} & \text{if $(\beta|\beta)\ge 0$},\\
 q^{-\frac{k(k-1)}{2}r_{\mf g}}\bff_\beta^{(k)} & \text{if $(\beta|\beta)< 0$ and $\beta\ne(i,i)$},\\
 q^{-\frac{k(k-1)}{2}}\bff_\beta^{(k)} & \text{if $\mf{g}=\mf{b}_{m|n}$, $(\beta|\beta)< 0$, and $\beta=(i,i)$},\\
 q^{-k^2}\bff_\beta^{(k)} & \text{if $\mf{g}=\mf{d}_{m|n}$, $(\beta|\beta)< 0$, and $\beta=(i,i)$},
\end{cases}
\end{equation*}
(see \cite[Section 3.1]{JKU24} for the definition of $\bff_\beta^{(k)}$).
Let
\begin{equation*}\label{eq:M(u)}
	M(\mf{u}) = 
	\left\{\, (c_\beta)_{\beta\in \Phi^+(\mf{u})}\in \Z_{\ge 0}^{N} \, | \,
	\text{$c_{\beta}\in\Z_{\ge 0}$ $(\beta\in \Phi^+_{\ov 0}\cup \Phi^+_{\rm{n\text{-}iso}})$,\quad $c_{\beta}=0,1$  $(\beta\in \Phi^+_{{\rm iso}})$} 
	\,\right\},
\end{equation*}
where $N=|\Phi^+(\mf{u})|$ and $\Phi^+_{\text{n-iso}}$ is the set of  non-isotropic roots in $\Phi^+_{\ov 1}$.
For ${\bf c}=(c_\beta)_{\beta\in \Phi^+(\mf{u})}\in M(\mf{u})$, let 
\begin{equation*}
\bfF^{(\bf c)} = \prod^{\rightarrow}_{\beta\in\Phi^+(\mf{u})}\bfF_\beta^{(c_\beta)},
\end{equation*}
where the product is taken in increasing order with respect to $\prec$.

\begin{thm}{\em (\cite[Theorem 4.15]{JKU24}) }\label{thm:crystal base of N}
Let 
\begin{equation*}\label{eq:crystal base of N}
\begin{split}
 \ms{L}\left(\mc{N}\right) &= \bigoplus_{{\bf c}\in M(\mf{u})} A_0 \bfF^{(\bf c)},\ \
 \ms{B}\left(\mc{N}\right) = \left\{\,\pm  \bfF^{(\bf c)}\!\!\!\pmod{q\ms{L}\left(\mc{N}\right)}\,|\,{\bf c}\in M(\mf{u})\,\right\}.
\end{split} 
\end{equation*}
Then $(\ms{L}\left(\mc{N}\right),\ms{B}\left(\mc{N}\right))$ is a crystal base of $\mc{N}$ as a $\U(\mf{l})$-module.
\end{thm}

By \cite[Theorems 4.22 and 4.24]{JKU24}, we have 
\begin{equation} \label{eq: crystal base of Pla}
\ms{L}(P(\la))=\ms{L}(\mc{N})\ot \ms{L}(V_{\mf{l}}(\la)),\quad
\ms{B}(P(\la))=\ms{B}(\mc{N})\ot \ms{B}(V_{\mf{l}}(\la)).
\end{equation}
The crystal operators $\te_i, \tf_i$ for $i\in I\setminus\{0\}$ act on $\ms{B}(\mc{N})\ot \ms{B}(V_{\mf{l}}(\la))$ by tensor product rule \cite[Proposition 2.8]{BKK}, while $\te_0, \tf_0$ acts by decreasing and increasing the multiplicity $c_{\alpha_0}$ in $\bfF^{(\bf c)}$ by $1$, respectively.

\subsection{Compatibility of crystal bases of $V(\la,\ell)$ and $P(\la)$}
Let $(\la,\ell)\in \cP(\g)$ be given.
We show that the crystal base of $P(\la)$ is compatible with that of $V(\la,\ell)$ under $\pi_{V(\la,\ell)}$.
This can be done by following the same idea in \cite[Section 4]{Kas91}.  

Let $\ell'>\ell\ge 0$ be given. 
Note that $(\la,\ell')\in \cP(\g)$.
For $t=\ell'-\ell$, we have a $\U(\g)$-linear embedding 
\begin{equation*}
\xymatrixcolsep{2pc}\xymatrixrowsep{0pc}\xymatrix{
\Phi:V(\la,\ell')  \ \ar@{->}[r] &\ V(\la,\ell)\ot V(0,t) \\
\qquad v_{(\la,\ell')} \ \ar@{|->}[r] &\ v_{(\la,\ell)} \ot v_{(0,t)}
}
\end{equation*}
and a $\U(\g)^-$-linear map
\begin{equation*}
\xymatrixcolsep{2pc}\xymatrixrowsep{0pc}\xymatrix{
S:V(\la,\ell)\ot V(0,t)  \ \ar@{->}[r] &\ V(\la,\ell) \\
 v \ot v_{(0,t)} \ \ar@{|->}[r] &\ v \\
 v \ot v' \ \ar@{|->}[r] &\ 0 
},
\end{equation*}
for $v'$ such that $\mathrm{wt}(v') \ne \La(0,t)=t\Lambda_0$.
Since $\mc{O}(\g)_{\ge 0}$ is a semisimple tensor category, $\Phi$ is a well-defined $\U(\g)$-linear map, and in particular it commutes with $\te_i$ and $\tf_i$ for $i\in I$.

\begin{lem} \label{lem: compatibility lemma1}
Under the above hypothesis, 
\begin{itemize}
 \item[(1)] $S\circ \Phi (\ms{L}(V(\la,\ell'))) =\ms{L}(V(\la,\ell))$,
 
 \item[(2)] the map induced from $S\circ \Phi$ at $q=0$ gives a bijection
\begin{equation*}
\overline{S\circ \Phi} : \{ \, b \in \ms{B}(V(\la,\ell')) \, | \, \ov{S\circ \Phi} (b)   \ne 0 \,  \} \longrightarrow \ms{B}(V(\la,\ell)),
\end{equation*}
which commutes with $\te_i$ and $\tf_i$ for $i\in I$.
\end{itemize}
\end{lem}
\pf Since $\Phi$ commutes with $\te_i$ and $\tf_i$ for $i\in I$, it follows from the connectedness of $\ms{B}(V(\la,\ell'))$ that $\Phi(\ms{L}(V(\la,\ell')))\subset \ms{L}(V(\la,\ell))\ot \ms{L}(V(0,t))$. Hence $S\circ \Phi(\ms{L}(V(\la,\ell')))\subset \ms{L}(V(\la,\ell))$ since $S(\ms{L}(V(\la,\ell))\ot \ms{L}(V(0,t)))=\ms{L}(V(\la,\ell))$ by definition of $S$. 
 
Let $b \in \ms{B}(V(\la,\ell))$ be given, where $b=\td{x}_{i_k}\dots \td{x}_{i_1} v_{(\la,\ell)}$ for some $i_1, \dots, i_k \in I $, where $x=e,f$ for each $i_j$. 
By tensor product rule (cf. \cite[Proposition 2.8]{BKK}), we have
\begin{equation*}
\tx_{i_k}\dots \tx_{i_1} (v_{(\la,\ell)}\ot v_{(0,t)}) \equiv \pm (\tx_{i_k}\dots \tx_{i_1} v_{(\la,\ell)})\ot v_{(0,t)}  \ (\mathrm{mod} \ q\ms{L}(V(\la,\ell))\ot \ms{L}(V(0,t))),
\end{equation*} 
since $\tf_i v_{(0,t)}=0$ for all $i>m$. Hence 
$\ov{S\circ \Phi}(\tx_{i_k}\dots \tx_{i_1} v_{(\la,\ell')})=b$. In particular, $\ov{S\circ \Phi}$ is surjective and hence $S\circ \Phi (\ms{L}(V(\la,\ell'))) =\ms{L}(V(\la,\ell))$. This proves (1).
 
To prove the injectivity of $\ov{S\circ \Phi}$, suppose $b, b' \in \ms{B}(V(\la,\ell'))$ satisfy $\ov{S\circ \Phi}(b) = \ov{S\circ \Phi}(b') \neq 0$.
Write $b=\td{x}_{i_k}\dots \td{x}_{i_1} v_{(\la,\ell')}$, $b'=\td{x}_{j_l}\dots \td{x}_{j_1} v_{(\la,\ell')}$ for some $i_1, \dots, i_k, j_1, \dots j_l \in I $, where each $x$ stand for $e$ or $f$.
Then by assumption, we have $\td{x}_{i_k}\dots \td{x}_{i_1} v_{(\la,\ell)} = \td{x}_{j_l}\dots \td{x}_{j_1} v_{(\la,\ell)}$.
On the other hand, 
$$\ov{\Phi}(b)=\tx_{i_k} \dots \tx_{i_1} (v_{(\la,\ell) }\ot v_{(0,t)}) = (\tx_{i_k} \dots \tx_{i_1} v_{(\la,\ell)}) \ot v_{(0,t)} = \ov{\Phi}(b'). $$
Since $\ov{\Phi}$ is an embedding of crystals, it follows that $b=b'$.
\qed
\medskip

Let $\Phi_{\ell', \ell}=S\circ \Phi$, which is $\U(\g)^-$-linear, and let $\Psi_{\ell}:\U(\g)^- \longrightarrow V(\la,\ell)$ be the $\U(\g)^-$-linear projection sending $1$ to $v_{(\la,\ell)}$.
Since $\Phi_{\ell',\ell} \circ \Psi_{\ell'} = \Psi_{\ell}$, we have $\mathrm{Ker} \, \Psi_{\ell'} \subset \mathrm{Ker} \, \Psi_{\ell}$, especially $(\mathrm{Ker} \, \Psi_{\ell'})_{\xi} \subset (\mathrm{Ker} \, \Psi_{\ell})_{\xi} \subset \U(\g)^-_\xi$ for $\xi\in Q^-$.
Hence for a given $\xi\in Q^-$ there is $N>0$ such that $(\mathrm{Ker} \, \Psi_{\ell+t})_{\xi} = (\mathrm{Ker} \, \Psi_{\ell+N})_{\xi}$ for all $t\ge N$.

This implies that the restriction $\Phi_{\ell+t,\ell+N} : V(\la,\ell+t)_{\Lambda(\la,\ell+t)+\xi} \longrightarrow V(\la,\ell+N)_{\Lambda(\la,\ell+N)+\xi}$ is a $\Bbbk$-linear isomorphism, and 
\begin{equation}\label{eq:stability of crystal lattice}
\Phi_{\ell+t,\ell+N} : \ms{L}(V(\la,\ell+t))_{\Lambda(\la,\ell+t)+\xi} \longrightarrow \ms{L}(V(\la,\ell+N))_{\Lambda(\la,\ell+N)+\xi}
\end{equation}
is an isomorphism of $A_0$-modules for all $t\ge N$. 
	
\begin{lem}\label{lem: compatibility lemma2} {\em (cf. \cite[Lemma 4.4.1]{Kas91})}
Let $\xi \in Q^-$ be given. For $P \in \U(\g)^-_{\xi}$, there exists $N >0$ such that for all $t >N$,	
\begin{equation*}
\begin{split}
(\td{x}_0 P) v_{(\la,\ell+t)} & \equiv \td{x}_0(P v_{(\la,\ell+t)}) \quad (\mathrm{mod} \ q\ms{L}(V(\la,\ell+t))) \quad (x=e,f),
\end{split}
\end{equation*}
where $\td{x}_0$ on $\U(\g)^-$ is defined in the same way as in \eqref{eq: crystal operator for 0}.
\end{lem}
\pf The argument is essentially the same as in \cite[Lemma 4.4.1]{Kas91} after a minor modification.
Assume that $P=f_0^{(k)}Q$ with $e'_0 Q=0$. Then 
$\te_0 P  = f_0^{(k-1)}Q$ and $\tf_0 P  = f_0^{(k+1)}Q$.
We will show that 
\begin{equation*}
\begin{split}
\te_0(P v_{(\la,\ell+t)}) & \equiv  f_0^{(k-1)}Q v_{(\la,\ell+t)},\ \ 
\tf_0(P v_{(\la,\ell+t)})  \equiv  f_0^{(k+1)}Q v_{(\la,\ell+t)} \ \ (\mathrm{mod} \ q\ms{L}(V(\la,\ell+t))),
\end{split}
\end{equation*}
for all sufficiently large $t$.

We remark that \cite[Corollary 3.4.6, Lemma 4.3.9]{Kas91} still hold for $\U(\g)^-$ and $V(\la,\ell)$, where $(L(\la),B(\la))$ and $u$ are replaced by $(\ms{L}(V(\la,\ell)),\ms{B}(V(\la,\ell)))$ and $Qv_{(\la,\ell+t)}$ respectively. 
Thus by \cite[Lemma 4.3.9]{Kas91}, it is enough to show that 
\begin{equation*}
k_0^\nu e_0^{(\nu)} Q v_{(\la,\ell+t)} \in q_0^{\nu(\nu+k+1)}q \ms{L}(V(\la,\ell+t)) \quad \text{for all } 1 \leq \nu \leq k+1,
\end{equation*}
where we have by \cite[Corollary 3.4.6]{Kas91}, 
\begin{equation*}
\begin{split}
k_0^\nu e_0^{(\nu)} Q v_{(\la,\ell+t)} & = \frac{\bq(\alpha_0, \La_\la+(\ell+t)\Lambda_0+\xi-k\alpha_0)^{2\nu} q_0 ^{3\nu^2+\nu} }{(q_0 - q_0^{-1})^\nu} ({e''_0}^{\, \nu} Q) v_{(\la,\ell+t)}, \\
& = \frac{\bq(\alpha_0, \La_\la+\ell \Lambda_0 +\xi-k\alpha_0)^{2\nu} q_0 ^{3\nu^2+\nu +t\nu} }{(q_0 - q_0^{-1})^\nu} ({e''_0}^{\, \nu} Q) v_{(\la,\ell+t)}.
\end{split}
\end{equation*}

Let $N$ be as in \eqref{eq:stability of crystal lattice}. 
Then we have $P v_{(\la,\ell+t)} \in \ms{L}(V(\la,\ell+t))$ if and only if $P v_{(\la,\ell+N)} \in \ms{L}(V(\la,\ell+N))$ for all $t\ge N$.
On the other hand, there exists $M_\nu \in \Z$ such that $q^{M_\nu} ({e''_0}^{\, \nu} Q) v_{(\la,\ell)} \in \ms{L}(V(\la,\ell+t))$ for all $t \geq N$. Then one can choose $N'\geq N$ such that 
\begin{equation*}
\frac{\bq(\alpha_0, \La_\la+\ell \Lambda_0 +\xi-k\alpha_0)^{2\nu} q_0 ^{3\nu^2+\nu +t\nu} }{(q_0 - q_0^{-1})^\nu} \in q_0^{\nu(\nu+k+1)+M_\nu} q A_0
\end{equation*}
for all $t\geq N'$ and $1 \leq \nu \leq k+1$. This completes the proof.
\qed
\vskip 2mm
The following lemma shows that $P(\la)$ can be viewed as the limit of $V(\la,\ell)$ as $\ell \to \infty$.

\begin{lem} \label{lem: character of P(la)}
Let $(\lambda,\ell)\in \cP(\g)$ and $\xi \in Q^-$ be given. 
Then, for sufficiently large integers $t$, we have $P(\la)_{\La(\la,0)+\xi} \cong V(\la,\ell)_{\La(\la,\ell+t)+\xi}$ as vector spaces.
\end{lem}
\pf
This can be seen from the separation algorithm in Section \ref{subsec: separation}. See Remark \ref{rem: proof of ch of P(la)}. 
\qed

\begin{thm}\label{thm: compatibility 1}
 Let {$(\lambda,\ell)\in \cP(\g)$} be given.
Then  
\begin{enumerate}
\item[\em (1)] $\pi_{V(\la,\ell)}(\ms{L}(P(\la)))=\ms{L}(V(\la,\ell))$,

\item[\em (2)] $\ov{\pi}_{V(\la,\ell)}(\ms{B}(P(\la)))= \ms{B}(V(\la,\ell)) \cup\{0\}$, where $$\ov{\pi}_{V(\la,\ell)} : \ms{L}(P(\la))/q\ms{L}(P(\la)) \rightarrow \ms{L}(V(\la,\ell))/q\ms{L}(V(\la,\ell))$$ is the induced $\Q$-linear map at $q=0$,

\item[\em (3)] $\ov{\pi}_{V(\la,\ell)}$ restricts to a bijection 
$$\ov{\pi}_{V(\la,\ell)} : \{\, b\in \ms{B}(P(\la))\,|\,\ov{\pi}_{V(\la,\ell)}(b)\neq 0 \,\} \longrightarrow \ms{B}(V(\la,\ell)),$$ 
which commutes with $\te_i$ and $\tf_i$ for $i\in I$.
\end{enumerate}
\end{thm}
\pf Since $P(\la)$ and $V(\la,\ell)$ are semisimple $\U(\mf{l})$-modules, every $\U(\mf{l})$-irreducible component $L$ of $P(\la)$ is isomorphic to its image under $\pi_{V(\la,\ell)}$ if it is not zero.
Thus $\pi_{V(\la,\ell)}$ commutes with $\te_i$ and $\tf_i$ if $i\in I\setminus \{ 0 \}$.
	
Now we claim that 
\begin{equation}\label{eq: compatibility lemma}
\pi_{V(\la,\ell)} (\tx_{i_k} \dots \tx_{i_1} (1 \ot v_\la) ) \equiv \pm \tx_{i_k} \dots \tx_{i_1} v_{(\la,\ell)} \quad (\mathrm{mod} \ q \ms{L}(V(\la,\ell))),
\end{equation}
for all $(\la,\ell) \in \cP(\g)$ and $i_i,\dots,i_k\in I$.
We use induction on $k$. It is clearly true when $k=0$.
Suppose that \eqref{eq: compatibility lemma} is true for $k-1$.

If $i_k \ne 0$, then \eqref{eq: compatibility lemma} is true since $\tx_i$ commutes with $\pi_{V(\la,\ell)}$.
So we assume $i_k=0$. 
	Let $P\in \U(\g)^-$ such that $\tx_{i_{k-1}}\dots \tx_{i_1}(1 \ot v_\la)=P(1\ot v_\la)$. 
By induction hypothesis, 
\begin{equation}\label{eq:induction k-1}
P v_{(\la,\ell)} \equiv \pm \tx_{i_{k-1}} \dots \tx_{i_1} v_{(\la,\ell)}  \  (\mathrm{mod} \ q \ms{L}(V(\la,\ell)))
\end{equation}
for all $(\la,\ell) \in \cP(\g)$.
By Lemma \ref{lem: compatibility lemma2} and \eqref{eq:induction k-1}, there exists $N$ such that for all $t \geq N$
\begin{equation}\label{eq:induction k}
 (\tx_0 P) v_{(\la,\ell+t)} \equiv \tx_0 (P v_{(\la,\ell+t)}) 
 \equiv \pm \tx_0 \tx_{i_{k-1}} \dots \tx_{i_1} v_{(\la,\ell+t)} \  (\mathrm{mod} \ q \ms{L}(V(\la,\ell+t))).
\end{equation} 
Applying $\Phi_{\ell+t,\ell}$ to \eqref{eq:induction k} and using Lemma \ref{lem: compatibility lemma1} together with induction hypothesis, we obtain
\begin{equation}\label{eq:induction ell}
 (\tx_0 P) v_{(\la,\ell)} \equiv \pm \tx_0 \tx_{i_{k-1}} \dots \tx_{i_1} v_{(\la, \ell)} \equiv \pm \tx_0 (P v_{(\la,\ell)}) \  (\mathrm{mod} \ q \ms{L}(V(\la,\ell))).
\end{equation}
Finally, since $\tx_0 (P(1\ot v_\la)) = (\tx_0 P)(1 \ot v_\la)$ by definition, we have
\begin{equation*}
\begin{split}
			\pi_{V(\la,\ell)} (\tx_{i_k} \dots \tx_{i_1} (1 \ot v_\la) ) & = \pi_{V(\la,\ell)} (\tx_{0} (P (1 \ot v_\la)) ) \\
			& = \pi_{V(\la,\ell)} ((\tx_0 P) (1 \ot v_\la) ) \\
			& = (\tx_0 P) v_{(\la,\ell)} \\
			& \equiv \pm \tx_0 (P v_{(\la,\ell)}) \quad{(\mathrm{mod} \ q\ms{L}(V(\la,\ell)))}  \\
			& \equiv \pm \tx_0 \tx_{i_{k-1}} \dots \tx_{i_1} v_{(\la,\ell)} \quad{(\mathrm{mod} \ q\ms{L}(V(\la,\ell)))}, 
\end{split}
\end{equation*}
where the second-to-last equivalence follows from \eqref{eq:induction ell}.
This completes the proof of \eqref{eq: compatibility lemma}. 
Then (1) and (2) follow from \eqref{eq: compatibility lemma} and the connectedness of $\ms{B}(P(\la))$ and $\ms{B}(V(\la,\ell))$.
We now prove (3). 
Suppose there exist $b,b' \in \ms{B}(P(\la))_{\La(\la,0)+\xi}$ for some $\xi \in Q^-$ with $\ov{\pi}_{V(\la,\ell)}(b)=\ov{\pi}_{V(\la,\ell)}(b') \neq 0$.
By Lemma \ref{lem: character of P(la)}, for a sufficiently large integer $t$ we have an isomorphism $\ms{B}(P(\la))_{\La(\la,0)+\xi} \cong \ms{B}(V(\la,\ell+t))_{\La(\la,t)+\xi}$, 
and therefore $\ov{\pi}_{V(\la,\ell+t)}(b) \neq \ov{\pi}_{V(\la,\ell+t)}(b')$.
Combining this with Lemma \ref{lem: compatibility lemma1} (2) and $\pi_{V(\la,\ell)} = \Phi_{\ell+t,\ell+N} \circ \pi_{V(\la,\ell+t)}$ yields $\ov{\pi}_{V(\la,\ell)}(b) \neq \ov{\pi}_{V(\la,\ell)}(b')$, a contradiction.
Together with \eqref{eq: compatibility lemma}, the proof of (3) is complete.
\qed
\section{Crystal base of $\U(\mf{l})^-$}\label{sec:review on type A}\label{sec:U^- of type A}
Let us briefly recall the results on the crystal base of $\U(\mf{l})^-$ in \cite{JKU23}.

\subsection{Crystal base of a Kac module}
Let us first recall the crystal bases of $q$-deformed Kac modules \cite[Section 4]{K14}.
Let $\Phi^+(\mf{l})_{\kappa}=\Phi^+(\mf{l})\cap \Phi^+_{\kappa}$ ($\kappa={\ov 0}, {\ov 1}$).
Let $\bf{f}'_\beta$ be the root vector for $\beta\in \Phi^+(\mf{l})$ with a partial order $\prec'$ on $\Phi^+(\mf{l})$ defined as in \cite[Section 2.5]{JKU23}, where $\alpha \prec'\beta$ for $\alpha\in \Phi^+(\mf{l})_{\ov 1}$ and $\beta\in \Phi^+(\mf{l})_{\ov 0}$. 
Note that $\bf{f}'_\beta$ is different from $\bf{f}_\beta$ in Section \ref{subsec:root vectors}. 

Let $\mc{K}$ be the $\Bbbk$-span of ${\bf f}'_S$ for $S\subset \Phi^+(\mf{l})_{\ov 1}$, where ${\bf f}'_S$ be the product of $\bf{f}'_\alpha$ ($\alpha\in S$) with respect to $\prec'$. With respect to the associated PBW type basis, we have
\begin{equation}\label{eq:Levi decomp 2}
 \U(\mf{l})^-
 \cong \mc{K} \ot \U(\mf{l}_{\ov 0})^-
 \cong \mc{K} \ot \U(\mf{l}_{m|0})^-\ot \U(\mf{l}_{0|n})^-.
\end{equation}
Moreover $\mc{K}$ is equal to the subalgebra generated by $\bf{f}'_\alpha$ for $\alpha\in \Phi^+(\mf{l})_{\ov 1}$, and it is a $\U(\mf{l}_{\ov 0})$-module with respect to \eqref{eq: quantum adjoint}, which is a direct sum of $V_{\mf{l}_{\ov 0}}(\la)$ for some $\mf{l}_{\ov 0}$-dominant weights $\la\in \texttt{P}^+$.
The pair $(\ms{L}(\K),\ms{B}(\K))$ with
\begin{equation*}\label{eq:LK BK}
 \ms{L}(\K)=\bigoplus_{S\subset \Phi^+(\mf{l})_{\ov 1}}A_0\bff'_S,\quad \ms{B}(\K)=\{\,\pm {\bf f}'_S \pmod{q\ms{L}(\K)}\,|\,S\subset \Phi^+(\mf{l})_{\ov 1}\,\}
\end{equation*}
is a crystal base of $\mc{K}$ as a $\U(\mf{l}_{\ov 0})$-module.

For $\la\in {\tt P}^+$,
let
\begin{equation*}\label{eq:Kac}
K(\la) = \U(\mf{l})\otimes_{\mc{Q}}V_{\mf{l}_{\ov 0}}(\la),
\end{equation*}
where $\mc{Q}$ is the subalgebra of $\U(\mf{l})$ generated by $\U(\mf{l}_{\ov 0})$ and $\U(\mf{l})^+$, and $V_{\mf{l}_{\ov 0}}(\la)$ is extended to a $\mc{Q}$-module trivially. It is a highest weight $\U(\mf{l})$-module generated by $1\ot v_{\la}$ with highest weight $\la$.
As a $\U(\mf{l}_{\ov 0})$-module, we have
\begin{equation}\label{eq:Kac as l-module}
K(\la) \cong \mc{K}\ot V_{\mf{l}_{\ov 0}}(\la).
\end{equation}
We define a crystal base of $K(\la)$ as in Section \ref{subsubsec: crystal base of osc}, where $\te_m, \tf_m$ is replaced by
\begin{equation}\label{eq: odd crystal operator}
 \te_m u = e'_m(u),\quad \tf_m u = f_m u.
\end{equation}
Then $K(\la)$ has a unique crystal base $\left(\ms{L}({K(\la)}),\ms{B}(K(\la))\right)$, where $\ms{B}(K(\la))$ is connected.
Let $(\ms{L}(V_{\mf{l}_{\ov 0}}(\la)), \ms{B}(V_{\mf{l}_{\ov 0}}(\la)))$ be the crystal base of $V_{\mf{l}_{\ov 0}}(\la)$ generated by $v_{\la}$.
Then we have
\begin{equation*}
(\ms{L}({K(\la)}),\ms{B}(K(\la))) \cong (\ms{L}(\mc{K})\ot \ms{L}(V_{\mf{l}_{\ov 0}}(\la)),\ms{B}(\mc{K}) \otimes \ms{B}(V_{\mf{l}_{\ov 0}}(\la))),
\end{equation*}
under \eqref{eq:Kac as l-module}, where $\cong$ means the isomorphism of crystal bases.

For $\la\in \cP_{m|n}$, let  
\begin{equation}\label{eq: Kac to poly}
 \xymatrixcolsep{3pc}\xymatrixrowsep{0pc}\xymatrix{
  K(\la) \ \ar@{->}[r]^{\pi_{V_{\mf{l}}(\la)} } & V_{\mf l}(\la)} 
\end{equation}
be the $\U(\mf{l})$-module homomorphism such that $\pi_{V_{\mf{l}}(\la)}(1\ot v_\la)=v_\la$.
Then $\pi_{V_{\mf{l}}(\la)}$ preserves the crystal bases in the sense of Theorem \ref{thm: compatibility 1}. 

For $\la,\mu\in {\tt P}^+$, we define $\la<\mu$ if and only if $\mu-\la\in {\tt P}^+$.
It is shown in \cite{JKU23} that $\{\,\ms{B}(K(\la))\,|\,\la\in {\tt P}^+\,\}$ forms an inductive system of (abstract) $\U({\mf l}_{\ov 0})$-crystals (not $\U(\mf l)$-crystal), with respect to the injective morphism of $\U({\mf l}_{\ov 0})$-crystals;
\begin{equation*} \label{eq: Theta map}
 \xymatrixcolsep{3pc}\xymatrixrowsep{0pc}\xymatrix{
 \ms{B}({K(\la)})  \ \ar@{->}[r]^{\Theta_{\la,\mu}} &\ \ms{B}({K(\mu)}) },
\end{equation*}
for $\la < \mu$.
Then its limit, say $\ms{B}(K(\infty))$, becomes an (abstract) $\U(\mf l)$-crystal (see \cite[Section 5.3]{JKU23} for more details), where for $b\in \ms{B}(K(\infty))$ 
\begin{equation}\label{eq:ep phi}
\varepsilon_i(b)=\max\{\,k\,|\,\te_i^kb\neq {\bf 0}\,\},\quad
 \varphi_i(b) =
 \begin{cases}
(-1)^{\e_i}({\rm wt}(b),\alpha_i ) + \varepsilon_i(b) & \text{for $i\neq m$},\\
 \max\{\,k\,|\,\tf_i^kb\neq {\bf 0}\,\} & \text{for $i= m$},
 \end{cases}
\end{equation}
where ${\bf 0}$ denotes the formal zero (when we consider an abstract crystal).

\subsection{Crystal base of $\U(\mf{l})^-$}
Let us recall a crystal base of $\U(\mf{l})^-$ \cite[Section 6]{JKU23}. It is given in such a way that the associated crystal is $\ms{B}(K(\infty))$.

We first need a parabolic Verma module and its crystal base, which interpolates those of $K(\la)$ and $\U(\mf{l})^-$ \cite[Section 6.2]{JKU23}.
For $\la\in {\tt P}$ with $\la_-\in {\tt P}^+$,
let
\begin{equation*}\label{eq:X}
X(\la) = \U(\mf{l})\otimes_{\mc{R}}V_{\mf{l}_{0|n}}(\la),
\end{equation*}
where $\mc{R}$ is the subalgebra of $\U(\mf{l})$ generated by $\U(\mf{l}_{0|n})$ and $\U(\mf{l})^+$, and $V_{\mf{l}_{0|n}}(\la)$ is extended to an $\mc{R}$-module trivially.
It is also a highest weight $\U(\mf{l})$-module generated by $1\ot v_{\la}$ with highest weight $\la$.
As a $\Bbbk$-space, we have
\begin{equation}\label{eq:decomp of X(la)}
 X(\la) \cong \K\ot \U(\mf{l}_{m|0})^-\ot V_{\mf{l}_{0|n}}(\la). 
\end{equation}

We define a crystal base of $X(\la)$ as in $K(\la)$, where $\te_i, \tf_i$ for $i<m$ is replaced by those defined as in \eqref{eq: crystal operator for 0} with respect to $e'_i$.
Then $X(\la)$ also has a unique crystal base $\left(\ms{L}({X(\la)}),\ms{B}(X(\la))\right)$, where $\ms{B}(X(\la))$ is connected. 

Let $(\ms{L}(\U(\mf{l}_{m|0})^-),\ms{B}(\U(\mf{l}_{m|0})^-))$ and $(\ms{L}(V_{\mf{l}_{0|n}}(\la)),\ms{B}(V_{\mf{l}_{0|n}}(\la)))$ be the crystal bases of  $\U(\mf{l}_{m|0})^-$ and $V_{\mf{l}_{0|n}}(\la)$, respectively.
We have
\begin{equation*}
\begin{split}
\ms{L}(X(\la)) &\cong 
\left(\ms{L}(\K)\cdot \ms{L}(\U(\mf{l}_{m|0})^-)\right)\ot \ms{L}(V_{\mf{l}_{0|n}}(\la)),\\
\ms{B}(X(\la)) &\cong 
\left(\ms{B}(\K) \cdot \ms{B}(\U(\mf{l}_{m|0})^-)\right) \ot \ms{B}(V_{\mf{l}_{0|n}}(\la)),
\end{split}
\end{equation*}
under \eqref{eq:decomp of X(la)}.

For $\la\in {\tt P}^+$, let  
\begin{equation*}
 \xymatrixcolsep{3pc}\xymatrixrowsep{0pc}\xymatrix{
  X(\la) \ \ar@{->}[r]^{\pi_{K(\la)} } & K(\la)} 
\end{equation*}
be the $\U(\mf{l})$-module homomorphism such that $\pi_{K(\la)}(1\ot v_{\la})=1\ot v_{\la}$.
Then $\pi_{K(\la)}$ preserves the crystal bases in the sense of Theorem \ref{thm: compatibility 1}. 

Let $u\in \U(\mf{l})^-$ be given with $u=u_1 u_2 u_3$ for $u_1\in \mc{K}$, $u_2\in \U(\mf{l}_{m|0})^-$ and $u_3\in \U(\mf{l}_{0|n})^-$ by \eqref{eq:Levi decomp 2}. For $i\in I$ with $i>m$, we define
\begin{equation}\label{eq:Kashiwara operators for U(l)^-:>m}
\te_i u =  u_1u_2(\te_iu_3),\quad \tf_i u =  u_1u_2(\tf_iu_3),
\end{equation} 
where $\te_i$ and $\tf_i$ are those on $\U(\mf{l}_{0|n})^-$.
Then 
we define a crystal base of $\U(\mf{l})^-$ as in $X(\la)$, where $\te_i, \tf_i$ for $i>m$ are replaced by \eqref{eq:Kashiwara operators for U(l)^-:>m}. 

\begin{thm}{\em (\cite[Theorem 6.1]{JKU24}) }\label{thm:crystal base of Levi}
Let 
\begin{equation*}\label{eq: B l infty}
\begin{split}
\ms{L}(\U(\mf{l})^-) & = \ms{L}(\K)\cdot \ms{L}(\U(\mf{l}_{m|0})^-)\cdot \ms{L}(\U(\mf{l}_{0|n})^-), \\
\ms{B}(\U(\mf{l})^-)&=\ms{B}(\K) \cdot \ms{B}(\U(\mf{l}_{m|0})^-) \cdot \ms{B}(\U(\mf{l}_{0|n})^-),
\end{split}
\end{equation*}
where $\cdot$ denotes the multiplication in $\U^-$ and the induced one in $\ms{L}(\U(\mf{l})^-)/q\ms{L}(\U(\mf{l})^-)$ respectively.
Then $(\ms{L}(\U(\mf{l})^-),\ms{B}(\U(\mf{l})^-))$ is a crystal base of $\U(\mf{l})^-$, where $\ms{B}(\U(\mf{l})^-)$ is isomorphic to $\ms{B}(K(\infty))$. 
\end{thm}

Let $\la\in {\tt P}$ with $\la_-\in {\tt P}^+$.
Since $\left(\ms{L}({X(\la)}),\ms{B}(X(\la))\right)$ is a upper crystal base of $X(\la)$ as a $\U(\mf{l}_{0|n})$-module,
it is not compatible with that of $\U(\mf l)^-$ under the canonical projection from $\U(\mf{l})^-$ to $X(\la)$.

Let $(\ms{L}(\U(\mf{l}_{0|n})^-),\ms{B}(\U(\mf{l}_{0|n})^-))$ be the crystal base of $\U(\mf{l}_{0|n})^-$ (see \cite[(3.8)]{JKU23} for more details).
For $\la\in {\tt P}$ with $\la_-\in {\tt P}^+$, let $\pi_{\la}: \U(\mf{l}_{0|n})^- \longrightarrow V_{\mf{l}_{0|n}}(\la)$ be the canonical projection, and let 
\begin{equation*}
 \xymatrixcolsep{3pc}\xymatrixrowsep{0pc}\xymatrix{
  V_{\mf{l}_{0|n}}(\la) \ \ar@{->}[r]^{\pi^\vee_{\la}} & \U(\mf{l}_{0|n})^-}
\end{equation*}
be the dual of $\pi_{\la}$ induced from non-degenerate bilinear forms on $\U(\mf{l}_{0|n})^-$ and $V_{\mf{l}_{0|n}}(\la)$, which is injective. Then we have
\begin{equation}\label{eq:dual embedding}
\pi_\la^\vee(\ms{L}(V_{\mf{l}_{0|n}}(\la)))\subset \ms{L}(\U(\mf{l}_{0|n})^-),
\quad 
\ov{\pi}_\la^\vee(\ms{B}(V_{\mf{l}_{0|n}}(\la)))\subset \ms{B}(\U(\mf{l}_{0|n})^-),
\end{equation}
where $\ov{\pi}_\la^\vee$ is the induced map at $q=0$ and $\ov{\pi}_\la^\vee (\te_i b)= \te_i \ov{\pi}_\la^\vee (b)$ for $i>m$ and $b\in \ms{B}(V_{\mf{l}_{0|n}}(\la))$ (cf.~\cite[Section 3.2.2]{JKU23} and references therein).

For $\la\in {\tt P}$ with $\la_-\in {\tt P}^+$, let $\pi_{X(\la)}: \U(\mf{l})^- \longrightarrow X(\la)$
be the $\U(\mf{l})^-$-linear map such that $\pi_{X(\la)}(1)=1\ot v_{\la}$. By \eqref{eq:Levi decomp 2}, \eqref{eq:decomp of X(la)}, and \eqref{eq:dual embedding}, we immediately have 
\begin{equation}\label{eq: dual of proj to X}
 \xymatrixcolsep{3pc}\xymatrixrowsep{0pc}\xymatrix{
  X(\la)   \ar@{->}^{\pi_{X(\la)}^\vee}[r] & \U(\mf{l})^-}
\end{equation}
given by $\pi_{X(\la)}^\vee=1_{\mc{K}\ot \U(\mf{l}_{m|0})^-}\ot \pi^\vee_\la$. Then we have
\begin{equation}\label{eq: dual of proj to X crystal base}
 \pi_{X(\la)}^\vee(\ms{L}(X(\la)))\subset \ms{L}(\U(\mf{l})^-),\quad
 \ov{\pi}_{X(\la)}^\vee(\ms{B}(X(\la)))\subset \ms{B}(\U(\mf{l})^-),
\end{equation}
where $ \ov{\pi}_{X(\la)}^\vee$ is the induced map at $q=0$.
Moreover, $\ov{\pi}^\vee_{X(\la)}$ restricts to a bijection 
\begin{equation}\label{eq: dual of proj to X crystal}
 \ov{\pi}^\vee_{X(\la)} : \ms{B}(X(\la))  \longrightarrow \{\,b=b_1b_2b_3\in \ms{B}(\U(\mf{l})^-)\,|\ b_3 \ot t_{\la} \in \ms{B}(V_{\mf{l}_{0|n}}(\la)) \,\},
\end{equation}
where $b=b_1b_2b_3\in \ms{B}(\U(\mf{l})^-)$ is as in Theorem \ref{thm:crystal base of Levi} and $\{t_\la\}$ is the (abstract) crystal with $\te_i t_\la=\tf_i t_\la=0$, $\varepsilon_i(t_\la)=\varphi_i(t_\la)=-\infty$ for all $i\in I$, and ${\rm wt}(t_\la)=\La_{0,\la}$.
Here $\ov{\pi}^\vee_{X(\la)}$ commutes with $\te_i$, $\tf_i$ for $i<m$, while it commutes with $\te_i$, $\tf_i$ for $i\ge m$ and $b'=\ov{\pi}^\vee_{X(\la)}(b)$ with $\la$ sufficiently large (see \cite[Remark 6.6]{JKU23}).

Let $\la\in \cP_{m|n}$ be given. 
Summarizing, the crystal bases mentioned in this subsection are compatible with the following maps:
\begin{equation}\label{eq:projections gl}
 \xymatrixcolsep{3pc}\xymatrixrowsep{0pc}
 \xymatrix{
  \U(\mf{l})^- \ar@{<-}^{\ \ \ \pi^\vee_{X(\la)}}[r] &  X(\la) \ar@{->}^{\pi_{K(\la)}}[r] & K(\la)  \ar@{->}^{\pi_{V_{\mf{l}}(\la)}}[r] & V_{\mf{l}}(\la)
  }.
\end{equation}
The maps in \eqref{eq:projections gl} induce injective maps at $q=0$:
\begin{equation}\label{eq:embeddings gl crystal}
 \xymatrixcolsep{3pc}\xymatrixrowsep{0pc}
 \xymatrix{
  \ms{\U(\mf{l})^-} \ar@{<-}^{\ \ \ \iota_{X(\la)}}[r] &  \ms{B}(X(\la)) \ar@{<-}^{\iota_{K(\la)}}[r] & \ms{B}(K(\la))  \ar@{<-}^{\iota_{V_{\mf{l}}(\la)}}[r] & \ms{B}(V_{\mf{l}}(\la))
  }.
\end{equation}
A combinatorial description of $\iota_{V_{\mf{l}}(\la)}$ can be seen in \cite{K14}, and $\iota_{K(\la)}$ is essentially the one from $\ms{B}(V_{\mf{l}_{\ov 0}}(\la^+))$ into $\ms{B}(\U(\mf{l}_{m|0})^-)$ for $\la=\la^+\in \texttt{P}^+$ given by counting the number of entries in each row
of a semistandard Young tableau to give multiplicities of root vectors (and hence Lusztig data for a PBW type crystal), and the description of $\iota_{X(\la)}$ is similar.

We remark that the crystal $\ms{B}(\U(\mf{l})^-)$ is not connected for $n\ge 2$ unlike $\ms{B}(X(\la))$, $\ms{B}(K(\la))$, and $\ms{B}(V_{\mf l}(\la))$. 
More precisely, it is not connected if and only if $n\ge 2$, and \begin{equation*}
 \ms{B}(\U(\mf{l})^-) \cong \left(\ms{B}(\U(\mf{l}_{m|1})^-)\times \ms{B}(\U(\mf{l}_{0|n})^-)\right)^{\oplus 2^{m(n-1)}}
\end{equation*}
(see \cite[Theorem 6.8]{JKU23} for more details).

\section{Crystal base of $\U(\g)^-$}\label{sec:Main1}

For $\la\in \cP_{m|n}$, we have the following diagram:
\begin{equation*}\label{eq:main diagram}
\xymatrixcolsep{2.2pc}\xymatrixrowsep{3pc}\xymatrix{
  \U(\g)^- \ \ar@{->}[r]^{\pi_{R(\la)}}\ar@{->}[d]^{\cong} &\ R(\la) \ \ar@{->}[r]^{\pi_{Q(\la)}}\ar@{->}[d]^{\cong} &\ Q(\la) \ \ar@{->}[r]^{\pi_{P(\la)}}\ar@{->}[d]^{\cong} &\ P(\la) \ \ar@{->}[r]^{\pi_{V(\la,\ell)}}\ar@{->}[d]^{\cong} &\ V(\la,\ell)  \\
   \mc{N}\ot \U(\mf{l})^- \ar@{->}^{\ {}_{{\rm id}_{\mc{N}}\ot \pi_{X(\la)}}}[r] &  \mc{N}\ot X(\la) \ar@{->}^{{}_{{\rm id}_{\mc{N}}\ot\pi_{K(\la)}}}[r] & \mc{N}\ot K(\la)  \ar@{->}^{{}_{{\rm id}_{\mc{N}}\ot\pi_{V_{\mf{l}}(\la)}}}[r] & \mc{N}\ot V_{\mf{l}}(\la)
},
\end{equation*}
where the vertical maps are $\Bbbk$-linear isomorphisms. We construct crystal bases of $Q(\la)$, $R(\la)$, and $\U(\g)^-$ using Theorem \ref{thm:crystal base of N} and the results in Section \ref{sec:review on type A}.

\subsection{Crystal base of $Q(\la)$}

Let $B_q(\mf{l}_{1|1})$ be the $q$-boson algebra associated with $\mf{l}_{1|1}$, that is, an associative $\Bbbk$-algebra generated by $e'$ and $f$ subject to the relations
\begin{equation*}
	e'^2=f^2=0, \quad e'f+fe'=1.
\end{equation*}
Note that $\langle e'_m, f_m \rangle \cong B_q(\mf{l}_{1|1})$, where $e'_m$ acts as in \eqref{eq: q-derivation}, and $f_m$ acts by left multiplication on $\U(\g)^-$.

It is easy to see that $\U(\mf{l}_{1|1})^-=\Bbbk 1 \oplus \Bbbk f$ is a unique irreducible $B_q(\mf{l}_{1|1})$-module up to isomorphism, where $\U(\mf{l}_{1|1})=\langle e,f,k^{\pm 1} \rangle$.
For a $B_q(\mf{l}_{1|1})$-module, we define crystal operators $\td{e}$ and $\td{f}$ as \eqref{eq: odd crystal operator}. 
Then $\U(\mf{l}_{1|1})^-$ has a unique crystal base $(L,B)$ up to scalar multiplication, where $L=A_0 1 \oplus A_0 f$ and $B=\{1, f \ (\mathrm{mod} \ qL) \}$ with $\tf 1= f (\mathrm{mod} \ qL)$.

For $\la=a \delta_1 + b \delta_2$ with $a, b\in \Z_{\ge 0}$, let $V_{\mf{l}_{1|1}}(\la)$ be the irreducible highest weight $\U(\mf{l}_{1|1})$-module with highest weight $\la$. 
Then $V_{\mf{l}_{1|1}}(\la)\ot  \U(\mf{l}_{1|1})^-$ is a $B_q(\mf{l}_{1|1})$-module under the $\Bbbk$-algebra homomorphism $\Delta_{B_q(\mf{l}_{1|1})} : B_q(\mf{l}_{1|1}) \longrightarrow \U(\mf{l}_{1|1})\ot B_q(\mf{l}_{1|1})$, which is defined by 
\begin{equation*} \label{eq: comul on B}
	\begin{split}
		\Delta_{B_q(\mf{l}_{1|1})}(e') &= (q^{-1}-q) k e \ot 1 + k \ot e', \\
		\Delta_{B_q(\mf{l}_{1|1})}(f) &= f \ot 1 + k \ot f.
	\end{split}
\end{equation*}
	
\begin{lem}\label{eq: tensor product rule for odd isotropic}
Let $V_1=V_{\mf{l}_{1|1}}(\la)$ and $V_2=\U(\mf{l}_{1|1})^-$ with crystal bases $(L_1,B_1)$ and $({L}_2,{B}_2)$ respectively.
Then $(L_1 \ot L_2,B_1 \ot B_2)$ is a crystal base of $V_1\ot V_2$, where $\te$ and $\tf$ act on $B_1 \ot B_2$ by
{\allowdisplaybreaks
	\begin{equation*}
	\begin{split}
	\te(b_1\otimes b_2)=&
	\begin{cases}
	\te  b_1\otimes b_2, & \text{if } a+b >0, \\ 
	 b_1\otimes  \te b_2, & \text{if } a=b=0,
	\end{cases}
	\\
	\tf(b_1\otimes b_2)=&
	\begin{cases}
	 \tf b_1\otimes b_2, & \text{if } a+b >0, \\ 
	b_1\otimes  \tf b_2, & \text{if } a=b=0.
	\end{cases}
	\end{split}
\end{equation*}}
	
\end{lem}
\pf The proof for the case when $a=b=0$ is clear.
Suppose that $a+b>0$. Then we have $V_1=\Bbbk v_1 \oplus \Bbbk f v_1$ and it has a unique crystal base up to scalar multiplication.
So we may assume that $L_1= A_0 v_1 \oplus A_0 f v_1$ and $L_2=A_0 v_2 \oplus A_0 fv_2$ with $v_2=1\in \U(\mf{l}_{1|1})^-$.
We have 
\begin{gather*}
e' \cdot (v_1 \ot v_2)=0,\quad 
f \cdot (v_1 \ot v_2) = f v_1 \ot v_2 + q^{a+b} v_1 \ot f v_2 \equiv f v_1 \ot v_2 \ (\mathrm{mod} \ qL_1\ot L_2),\\
e' \cdot X=0,\quad f \cdot X = fv_1 \ot fv_2,
\end{gather*}
where $X=(1-q^{2(a+b)}) v_1 \ot f v_2 - q^{a+b} f v_1 \ot v_2 \equiv v_1 \ot f v_2  \ (\mathrm{mod} \ qL_1\ot L_2)$.
\qed
\medskip

For $\la\in {\tt P}^+$, consider $Q(\la)$ in \eqref{eq:Q}. 
As a $\Bbbk$-space, we have
\begin{equation}\label{eq: Q decomp}
Q(\la) 
\cong \mc{N}\ot K(\la)
\cong \mc{N}\ot \mc{K}\ot V_{\mf{l}_{\ov 0}}(\la),
\end{equation}
where the second isomorphism is given by sending $u_1\ot u_2 \ot v$ to $u_1 \ot u_2 v$ for $u_1\in \mc{N}$, $u_2\in \mc{K}$ and $v\in V_{\mf{l}_{\ov 0}}(\la)$ from the right to the left.

\begin{lem}\label{lem: Q as l_even}
The $\Bbbk$-linear map in \eqref{eq: Q decomp} is an isomorphism of $\U({\mf l}_{\ov 0})$-modules. 
\end{lem}
\pf By the same arguments as in \cite[Proposition 4.20]{JKU24}, the first $\Bbbk$-linear map in \eqref{eq: Q decomp} is $\U({\mf l}_{\ov 0})$-linear. Also by \cite[Proposition 4.1]{K14}, the $\Bbbk$-linear isomorphism \eqref{eq:Kac as l-module} is $\U({\mf l}_{\ov 0})$-linear, and hence the second $\Bbbk$-linear map in \eqref{eq: Q decomp} is $\U({\mf l}_{\ov 0})$-linear. 
\qed
\medskip

We define a crystal base of $Q(\la)$ as in Section \ref{subsubsec: crystal base of parabolic P} for $P(\la)$, where $\te_m, \tf_m$ are replaced by the ones in \eqref{eq: odd crystal operator}.

\begin{thm}\label{thm: crystal base of Q}
For $\la\in {\tt P}^+$, let
\begin{equation*}
\begin{split}
\ms{L}(Q(\la))&= \left(\ms{L}(\mc{N})\cdot \ms{L}(\mc{K})\right)\ot \ms{L}(V_{\mf{l}_{\ov 0}}(\la))\subset Q(\la),\\ 
\ms{B}(Q(\la))&=\left(\ms{B}(\mc{N})\cdot  \ms{B}(\mc{K}) \right)\ot \ms{B}(V_{\mf{l}_{\ov 0}}(\la))\subset \ms{L}(Q(\la))/q\ms{L}(Q(\la)).
\end{split}
\end{equation*}
where $\cdot$ denotes the multiplication in $\U(\g)^-$ and the induced one at $q=0$, respectively.
Then $(\ms{L}(Q(\la)), \ms{B}(Q(\la)))$ is a crystal base of 
$Q(\la)$, and the crystal $\ms{B}(Q(\la))$ is connected. 
\end{thm}
\pf Put $\ms{L}=\ms{L}(Q(\la))$ and $\ms{B}=\ms{B}(Q(\la))$. 
First we show that $(\ms{L},\ms{B})$ is a crystal base. It suffices to show that $\td{x}_i \ms{L}\subset \ms{L}$ and 
$\td{x}_i \ms{B}\subset \ms{B}\cup \{0\}$ for $i\in I$ and $x=e,f$,
which follows from similar arguments as in  \cite[Theorem 4.22]{JKU24} 
when $i=0$, and Lemma \ref{lem: Q as l_even} when $i\in I\setminus\{m\}$.
If $i=m$, it follows from Lemma \ref{eq: tensor product rule for odd isotropic}, since $\mc{N}$ is a semisimple $\U(\mf{l}_{1|1})$-module and $\mc{K}$ is a semisimple $B_q(\mf{l}_{1|1})$-module.

Next, we show that $\ms{B}$ is connected. 
As a set, we have
\begin{equation}\label{eq: crystal of Q}
 \ms{B}\cong
 \ms{B}(\mc{N})\times \ms{B}(\mc{K}) \times \ms{B}(V_{\mf{l}_{\ov 0}}(\la))=\ms{B}(\mc{N})\times  \ms{B}(K(\la)).
\end{equation}
Note that the crystal graph $\ms{B}(K(\la))$ is equal to that of $\ms{B}(K(\la+\ell\de_+))$ (up to shifting weights by $\ell\de_+$), where $\de_+=\de_1+\cdots+\de_m$ for all $\ell\ge 0$.
Hence we may identify $\ms{B}$ with $\ms{B}(\mc{N})\times   \ms{B}(K(\la+\ell\de_+))$ (up to shifting weights by $\ell\de_+$).

On the other hand, if $\ell$ is large enough, then $K(\la+\ell\de_+)=V_{\mf l}(\la+\ell\de_+)$, and $\ms{B}(K(\la+\ell\de_+))=\ms{B}(V_{\mf l}(\la+\ell\de_+))$.  
Since $\ms{B}(P(\la+\ell\de_+))=\ms{B}(\mc{N})\times \ms{B}(V_{\mf{l}}(\la+\ell\de_+))$ is connected by Theorem \ref{thm: crystal base of P}, $\ms{B}$ is connected
(see \cite[Section 6.4]{K14} and references therein).

\qed

\begin{cor}
For $\la\in \cP_{m|n}$, the map $\pi_{P(\la)}: Q(\la) \longrightarrow P(\la)$ preserves the crystal bases in the sense of Theorem \ref{thm: compatibility 1}.
\end{cor}
\pf It follows from the fact that the map \eqref{eq: Kac to poly} preserves the crystal bases and $\varepsilon_i, \varphi_i$ for $ i \in I \setminus \{ 0, m \}$ since $Q(\la)$ and $P(\la)$ are semisimple $\U(\mf{l}_{\ov 0})$-modules.
\qed
\medskip

The crystal structure on $\ms{B}(Q(\la))$ can be described explicitly from \eqref{eq: crystal of Q}. As an $\mf{l}_{\ov 0}$-crystal, we have
\begin{equation*}
\ms{B}(Q(\la))\cong \ms{B}(\mc{N})\otimes \ms{B}(K(\la)), 
\end{equation*}
while $\te_0, \tf_0$ act only on $\ms{B}(\mc{N})$, and the actions of $\te_m, \tf_m$ are given by Lemma \ref{eq: tensor product rule for odd isotropic}. A combinatorial description of $\ms{B}(\mc{N})$ and $\ms{B}(K(\la))$ can be found in \cite[Section 4.2]{JKU24} and \cite[Section 5.1]{K14}, respectively.

\subsection{Crystal base of $R(\la)$}

Let $B_q=B_q({\mf l}_{m|0})$ be the $q$-boson algebra associated with $\mf{l}_{m|0}$ generated by $e'_i$ and $f_i$ ($i=1,\dots,m-1$) \cite[Section 3.3]{Kas91}. Then $\U(\mf{l}_{m|0})^-$ is the unique irreducible $B_q$-module, where $e'_i$ acts as \eqref{eq: q-derivation} and $f_i$ acts as left multiplication. 

For $\la\in {\tt P}^+$ with $\la=\la_+$, let $V_{\mf{l}_{m|0}}(\la)$ be the irreducible highest weight $\U(\mf{l}_{m|0})$-module with highest weight $\la$. Then $V_{\mf{l}_{m|0}}(\la)\ot \U(\mf{l}_{m|0})^-$ is a $B_q$-module under the $\Bbbk$-algebra homomorphism $\cmB : B_q \longrightarrow \U(\mf{l}_{m|0})\ot B_q$ given by
\begin{equation*} \label{eq: comul on B}
	\cmB(e_i') = (q_i^{-1}-q_i) k_ie_i \ot 1 + k_i \ot e_i', \quad 
	\cmB(f_i) = f_i \ot 1 + k_i \ot f_i.
\end{equation*}
Let $V_1=V_{\mf{l}_{m|0}}(\la)$ and $V_2=\U(\mf{l}_{0|n})^-$ with crystal bases $(L_1,B_1)$ and $({L}_2,{B}_2)$, respectively.
Then $(L_1 \ot L_2,B_1 \ot B_2)$ is a crystal base of $V_1\ot V_2$, where $\te_i$ and $\tf_i$ $(i \in I)$ act on $B_1 \ot B_2$ by the same formula as in the case of a tensor product of crystal bases of $\U(\mf{l}_{m|0})$-modules \cite[Remark 3.5.1]{Kas91}.

For $\la\in {\tt P}$ with $\la_-\in {\tt P}^+$, consider $R(\la)$ in \eqref{eq:R}.
As a $\Bbbk$-space, we have
\begin{equation}\label{eq: R decomp}
R(\la) \cong \mc{N}\ot \K\ot \U(\mf{l}_{m|0})^-\ot V_{\mf{l}_{0|n}}(\la),
\end{equation}
where the isomorphism is given in a standard way (cf.~\eqref{eq: Q decomp}).

\begin{lem}\label{lem: R as l_m}
The $\Bbbk$-linear map in \eqref{eq: R decomp}  is an isomorphism of $\U({\mf l}_{0|n})$-modules, where we regard $\U(\mf{l}_{m|0})^-$ as a trivial $\U({\mf l}_{0|n})$-module, and also it is an isomorphism of $B_q$-modules, where we regard  $V_{\mf{l}_{0|n}}(\la)$ as a trivial $B_q$-module.
\end{lem}
\pf The proof is almost identical to that of Lemma \ref{lem: Q as l_even}. We leave it to the reader.
\qed
\medskip

We define a crystal base of $R(\la)$ as in $Q(\la)$, where $\te_i, \tf_i$ for $i=1,\dots,m-1$ is replaced by those defined as in \eqref{eq: crystal operator for 0} with respect to $e'_i$.
 
\begin{thm}\label{thm: crystal base of R}
For $\la\in {\tt P}$ with $\la_-\in {\tt P}^+$, let
\begin{equation*}
\begin{split}
\ms{L}(R(\la))&= \left(\ms{L}(\mc{N})\cdot \ms{L}(\K)\cdot \ms{L}(\U(\mf{l}_{m|0})^-)\right)\ot \ms{L}(V_{\mf{l}_{0|n}}(\la))\subset R(\la),\\ 
\ms{B}(R(\la))&=\left(\ms{B}(\mc{N})\cdot \ms{B}(\K) \cdot \ms{B}(\U(\mf{l}_{m|0})^-)\right) \ot \ms{B}(V_{\mf{l}_{0|n}}(\la))\subset \ms{L}(R(\la))/q\ms{L}(R(\la)).
\end{split}
\end{equation*}
Then $(\ms{L}(R(\la)), \ms{B}(R(\la)))$
is a crystal base of $R(\la)$, where the crystal $\ms{B}(R(\la))$ is connected. 
\end{thm}
\pf By the same arguments as in Theorem \ref{thm: crystal base of Q} using Lemma \ref{lem: R as l_m} in this case, 
we have that $(\ms{L}(R(\la)), \ms{B}(R(\la)))$ is a crystal base of $R(\la)$.
The proof of connectedness of $\ms{B}(R(\la))$ is almost identical to the proof of Theorem 6.2 in \cite{JKU23},
since $\ms{B}(R(\la))$ is isomorphic to the limit of $\ms{B}(P(\la))$ when $\la_+ \rightarrow \infty $.

As a set, we have 
\begin{equation}\label{eq: crystal of R}	
 \ms{B}(R(\la))\cong
 \ms{B}(\mc{N})\times \ms{B}(\K) \times \ms{B}(\U(\mf{l}_{m|0})^-) \times \ms{B}(V_{\mf{l}_{0|n}}(\la)).
\end{equation}
Recall that $\te_0, \tf_0$ act only on $\ms{B}(\mc{N})$.
The actions of $\te_i, \tf_i$ for $i\in I\setminus\{0\}$ can be described by applying them to
\begin{equation*}
\begin{cases}
 \ms{B}(\mc{N})\otimes \ms{B}(\K) \otimes \ms{B}(V_{\mf{l}_{0|n}}(\la)) & \text{using \cite[Theorem 1]{Kas91} for $i>m$},\\
 \ms{B}(\mc{N})\otimes \ms{B}(\K) \otimes \ms{B}(\U(\mf{l}_{m|0})^-) & \text{using \cite[Remark 3.5.1]{Kas91} for $i<m$},\\
 \ms{B}(\mc{N})\otimes \ms{B}(\K) & \text{using Lemma \ref{eq: tensor product rule for odd isotropic} for $i=m$}.
\end{cases}
\end{equation*}

\begin{cor}
For $\la\in {\tt P}^+$, the map $\pi_{Q(\la)}: R(\la) \longrightarrow Q(\la)$ preserves the crystal bases in the sense of Theorem \ref{thm: compatibility 1}.
\end{cor}
\pf Note that $\pi_{Q(\la)}={\rm id}_{\mc{N}\cdot\mc{K}}\ot \pi_{\la}\ot {\rm id}_{V_{\mf{l}_{0|n}}(\la)}$, where $\pi_{\la}: \U(\mf{l}_{m|0})^- \longrightarrow V_{\mf{l}_{m|0}}(\la)$ is the canonical projection. Since $\pi_\la$ preserves the crystal bases and $\varepsilon_i$, $\pi_{Q(\la)}$ also preserves crystal bases by Theorem \ref{thm: crystal base of R}.
\qed

\subsection{Crystal base of $\U(\g)^-$}

We define a crystal base of $\U(\g)^-$ as in $R(\la)$ except for  $\te_i$ and $ \tf_i$ for $i>m$, which are defined by
\begin{equation*}\label{eq:Kashiwara operators for U^-:>m}
\te_i u =  u_1(\te_iu_2),\quad \tf_i u =  u_1(\tf_iu_2),
\end{equation*}
for $u=u_1u_2\in \U(\g)^-$ with $u_1\in \mc{N}$ and $u_2\in \U(\mf l)^-$, where $\te_iu_2$ and $\tf_iu_2$ are given in \eqref{eq:Kashiwara operators for U(l)^-:>m}.

Let 
\begin{equation*}
\ms{B}(Q(\infty))=\ms{B}(\mc{N})\times \ms{B}(\U(\mf{l})^-)
\end{equation*}
as a set. We define a $\U(\g)$-crystal $\ms{B}(Q(\infty))$ as follows: 
Let $b=(S,b_0)$ be given. 
We let
${\rm wt}(b) = {\rm wt}(S) + {\rm wt}(b_0)$,
and let $\varepsilon_i(b)$ and $\varphi_i(b)$ be as in \eqref{eq:ep phi}. 
For $i\in I$ and $x=e,f$, we define
\begin{equation*}
\begin{split}
 \tilde{x}_i b &= 
 \begin{cases}
 (\tilde{x}_0S,b_0) & \text{if $i=0$},\\
 (S',b'_0) & \text{if $1\le i\le m$ and $\tilde{x}_i (S\ot b_0)= S'\ot b'_0$},\\
 (S,\tilde{x}_ib_0) & \text{if $i>m$},\\
 \end{cases}
\end{split}
\end{equation*}
where $\tilde{x}_i (S\ot b_0)$ is defined by using \cite[Remark 3.5.1]{Kas91} for $1\le i<m$ and using Lemma \ref{eq: tensor product rule for odd isotropic} for $i=m$.
Here we assume that $\tilde{x}_i b={\bf 0}$ if any of its component on the right-hand side is ${\bf 0}$.

\begin{thm}\label{thm:crystal base of U^-}
Let 
\begin{equation*}
\begin{split}
\ms{L}(\infty) & = \ms{L}(\mc{N})\cdot \ms{L}(\U(\mf l)^-), \ \
\ms{B}(\infty) = \ms{B}(\mc{N})\cdot \ms{B}(\U(\mf l)^-).
\end{split}
\end{equation*}
Then $(\ms{L}(\infty),\ms{B}(\infty))$ is a crystal base of $\U(\g)^-$, and $\ms{B}(\infty)\cong \ms{B}(Q(\infty))$ as $\U(\g)$-crystals.
\end{thm}

For $\la\in {\tt P}$ with $\la_-\in {\tt P}^+$, let
\begin{equation*} 
 \xymatrixcolsep{3pc}\xymatrixrowsep{0pc}\xymatrix{
  R(\la)   \ar@{->}^{\pi_{R(\la)}^\vee}[r] & \U(\g)^-}
\end{equation*}
given by ${\rm id}_{\mc{N}}\ot \pi_{X(\la)}^\vee$ (see \eqref{eq: dual of proj to X} for $\pi_{X(\la)}^\vee$).  

\begin{cor} We have the following. 
\begin{itemize}
 \item[(1)] $\pi^\vee_{R(\la)}(\ms{L}(R(\la)))\subset \ms{L}(\infty)$.
 
 \item[(2)] $\ov{\pi}^\vee_{R(\la)}(\ms{B}(R(\la)))\subset \ms{B}(\infty)$, where $\ov{\pi}^\vee_{R(\la)}$ is the $\Q$-linear map induced from $\pi_{R(\la)}^\vee$.
 
 \item[(3)] $\ov{\pi}^\vee_{R(\la)}$ restricts to a bijection 
\begin{equation*}
 \ov{\pi}^\vee_{R(\la)} : \ms{B}(R(\la))  \longrightarrow \{\,b=b_1b_2b_3b_4\in \ms{B}(\infty)\,|\ b_4 \ot t_{\la_-} \in \ms{B}(V_{\mf{l}_{0|n}}(\la)) \,\},
\end{equation*} 
which commutes with $\te_i$, $\tf_i$ for $i<m$, where $b=b_1b_2b_3b_4 \in \ms{B}(\infty)$ is as in \eqref{eq: crystal of R}.

\end{itemize}
\end{cor}
\pf It follows from \eqref{eq: dual of proj to X crystal base} and \eqref{eq: dual of proj to X crystal}.
\qed

\begin{rem} \label{rem: compatibility}
{\rm 
The map $\ov{\pi}^\vee_{R(\la)}$ commutes with $\te_i, \tf_i$ for $i\ge m$ only when $\la_-$ is sufficiently large. More precisely,
let $b_0 \in \ms{B}(\infty)$ be given with $b_0=\ov{\pi}^\vee_{R(\la)}(b)$ for some $\la\in {\tt P}$ with $\la_-\in \tt{P}^+$ and $b \in \ms{B}(R(\la))$.
Then 
$\ov{\pi}^\vee_{R(\la)}(\tilde{x}_i b)= \tilde{x}_i \ov{\pi}^\vee_{R(\la)}(b)$, 
for $x=e, f$ if $i>m$ and $\la_- \gg 0$, or if $i=m$ and $( {\rm wt}(b),\alpha_{m+1})\ll 0$ (cf.~\cite[Remark 6.6]{JKU23}).
}
\end{rem}

Summarizing, the crystal bases constructed in Section \ref{sec:Main1} are compatible with the following maps (partly compatible in case of $\pi^\vee_{R(\la)}$):
\begin{equation*}\label{eq: seq of maps for g}
 \xymatrixcolsep{2.5pc}\xymatrixrowsep{3pc}\xymatrix{
  \U(\g)^- \ \ar@{<-}[r]^{\ \  \pi^\vee_{R(\la)}}  &\ R(\la) \ \ar@{->}[r]^{\! \pi_{ Q(\la)}}  &\ Q(\la) \ \ar@{->}[r]^{\! \pi_{P(\la)}}  &\ P(\la) \ \ar@{->}[r]^{\! \! \pi_{V(\la,\ell)}}  &\ V(\la,\ell)
}.
\end{equation*}
This induces the following commuting diagram at $q=0$:
\begin{equation*}\label{eq:main diagram}
\xymatrixcolsep{3.5pc}\xymatrixrowsep{3pc}\xymatrix{
& \ms{B}(P(\la)) \ \ar@{->}[ddr] \ar@{<-}[r]^{\iota_{V(\la,\ell)}} \ar@{->}[d]^{\cong} 
&\ \ms{B}(V(\la,\ell)) 
 \\
& \ms{B}(\mc{N})\ot \ms{B}(V_{\mf{l}}(\la)) \ \ar@{->}[ddr]
&  \\
  \ms{B}(\U(\g)^-) \ \ar@{<-}[r] \ar@{->}[d]^{\cong} 
& \ms{B}(R(\la)) \ \ar@{<-}[r] \ar@{->}[d]^{\cong} 
& \ms{B}(Q(\la)) \ \ar@{->}[d]^{\cong} 
\\
  \ms{B}(\mc{N})\ot \ms{B}(\U(\mf{l})^-) \ar@{<-}[r] 
& \ms{B}(\mc{N})\ot \ms{B}(X(\la)) \ar@{<-}[r] 
& \ms{B}(\mc{N})\ot \ms{B}(K(\la))    
}
\end{equation*}
Since the maps in the bottom line are obtained from \eqref{eq:embeddings gl crystal} by tensoring with ${\rm id}_{\ms{B}(\mc{N})}$, it is enough to describe 
\begin{equation}\label{eq:main embedding}
\iota_{V(\la,\ell)}: \ms{B}(V(\la,\ell)) \longrightarrow \ms{B}(P(\la)). 
\end{equation}

\section{Combinatorial description of crystal embeddings}\label{sec: Main 2}
In this section, we give a combinatorial description of \eqref{eq:main embedding}, which is a super-analogue of the results in \cite{JK21,K18}. So we may assume that $n\ge 1$. 

The cases of $\mf{g}=\mf{b}_{m|n}$ and $\mf{c}_{m|n}$ can be naturally generalized from the results for $\mf{b}_{m+n}$ and $\mf{c}_{m+n}$ \cite{K18}, while the case of $\mf{d}_{m|n}$ is not. 
In this case, we need a crystal theoretic super-analogue of Burge correspondence of type $D$ (Theorem \ref{thm: Burge}). It is essentially the main result in this section and its proof will be given in Section \ref{sec:proof of thm:Burge}.

\subsection{Crystal $\ms{B}(V_{\mf l}(\la))$}\label{subsec: type A crystal}
Let us recall the results on the crystal $\ms{B}(V_{\mf l}(\la))$ for $\la\in \ms{P}_{m|n}$ \cite{BKK}.

Let $\mathcal{A}$ be a $\Z_2$-graded linearly ordered set.
For a skew Young diagram $\lambda/\mu$, we denote by $SST_{\mc{A}}(\lambda/\mu)$ the set of all $\mathcal{A}$-semistandard tableaux of shape $\lambda/\mu$ with entries in $\mathcal{A}$.
We write $SST_\mathcal{A}(\lambda/\mu)$ by $SST_{m|n}(\lambda/\mu)$ when $\mathcal{A} = \I$.
For $T \in SST_\mathcal{A}(\lambda/\mu)$, let ${\rm sh}(T)=\la/\mu$ denote the shape of $T$, and let $w(T)$ be the word of $T$, which is the word with letters in $\mathcal{A}$ obtained from $T$ by reading its entries column by column from the right to the left, and from top to bottom in each column.

We may identify $\I =SST_{m|n}((1))$, a crystal of $V_{\mf l}((1))$. Then $SST_{m|n}(\la/\mu)$ has a structure of an (abstract) $\U(\mf{l})$-crystal (or simply $\mf l$-crystal) under the identification of $T\in SST_{m|n}(\la/\mu)$ with $w(T)$, where we regard $w(T)$ as an element in $SST_{m|n}((1))^{\otimes |\la/\mu|}$ ($|\la/\mu|$ is the number of the boxes in $\la/\mu$) and apply the tensor product rule. Then we have $SST_{m|n}(\la)\cong \ms{B}(V_{\mf l}(\la))$ for $\la\in \ms{P}_{m|n}$.

For a skew Young diagram $\la/\mu$, let $(\la/\mu)^\pi$ be the skew Young diagram obtained by rotating $\la/\mu$ by $180^\circ$. 
We also write $\la^\pi/\mu^\pi=(\la/\mu)^\pi$.
For $T\in SST_{m|n}(\de^\pi)$ with $\de \in \cP_{m|n}$ and $i\in \I$, we denote by $(T\leftarrow i)$ the tableau obtained by applying the Schensted's column insertion of $i$ into $T$ in a reverse way starting from the rightmost column of $T$ so that $(T\leftarrow i)\in SST_{m|n}(\gamma^\pi)$ for some $\gamma\supset \de$, which is given by adding a box in a corner of $\de$.

\subsection{Crystal $\ms{B}(V(\la,\ell))$}\label{subsec:spinor model}
Let us recall the results on the realization of the crystal $\ms{B}(V(\la,\ell))$ in \cite{K15,K16}.
Here we follow the convention in \cite[Section 3]{K18-2}, which are slightly different from those in \cite{K15,K16}.

\begin{rem}\label{rem:spinor convention}
{\rm
The definitions of relevant notions in \cite{K18-2} are given only in terms of $\mathbb{N}$-semistandard tableaux with usual linear order and $\Z_2$-grading $\N=\N_{\ov 0}$. But we may still obtain the definition for $\ms{B}(V(\la,\ell))$ simply by replacing $\N$-semistandard tableaux with $\I$-semistandard ones. So we only recall necessary notations briefly and refer the reader to \cite{K18-2} for full details.
} 
\end{rem}

First we assume the following notations:
\begin{itemize}
\item[$\bullet$] ${\bf T}^\g(a)$ for $a\in \Z_{\ge 0}$,

\item[$\bullet$] ${\bf T}^{\texttt{sp}}$ for $\g = \mf{b}_{m|n}$ and $\mf{d}_{m|n}$,

\item[$\bullet$] $\ov{\bf T}(0)$, ${\bf T}^{\texttt{sp}+}$, ${\bf T}^{\texttt{sp}-}$ for $\g=\mf{d}_{m|n}$.
\end{itemize}
For example, ${\bf T}^{\mf{c}_{m|n}}(a)=\bigsqcup_{c\ge 0}SST_{m|n}\left(\lambda(a,0,c)\right)$, where $\lambda(a, b, c) = (2^{b+c}, 1^a) / (1^b)$, while ${\bf T}^{\mf{b}_{m|n}}(a)$ and ${\bf T}^{\mf{d}_{m|n}}(a)$ consist of $T\in SST_{m|n}\left(\lambda(a,b,c)\right)$ ($a,b\in \Z_{\ge 0}$) satisfying certain mild constraints.

Let ${\bf B}$ be one of the sets above, and we regard it as an (abstract) $\g$-crystal such that the $\mf{l}$-crystal structure is given as in Section \ref{subsec: type A crystal}, and $\te_0$, $\tf_0$ are defined as in \cite[Section 3.1]{K18-2}, where for $T \in {\bf B}$, 
\begin{equation*}
\begin{split}
{\rm wt}(T) &=
\begin{cases}
	2\Lambda_0 + {\rm wt}_{\mf{l}}(T) & \text{if (${\bf B} = {\bf T}^{\g}(a)$ with $\g = \mf{b}_{m|n}$)} \\
	 & \text{or (${\bf B} = {\bf T}^{\g}(a),\, \ov{\bf T}(0)$ with $\g = \mf{d}_{m|n}$),} \\
	\Lambda_0 + {\rm wt}_{\mf{l}}(T) & \text{otherwise,}
\end{cases} \\
	\varepsilon_i(T) &= \max\left\{ k \, | \, \te_i^k T  \neq {\bf 0} \right\}, \qquad
	\varphi_i(T) = \max\left\{ k \, | \, \tf_i^k T  \neq {\bf 0} \right\},
\end{split}
\end{equation*}
and ${\rm wt}_{\mf{l}}(T)$ denotes the $\mf{l}$-weight of $T \in {\bf B}$.
Indeed, they are crystals of $V(\la,\ell)$ as follows:
\begin{itemize}
	\item[$\bullet$] when $\g = \mf{c}_{m|n}$, 
	$
		{\bf T}^{\mf{c}_{m|n}}(a) \cong \ms{B}\left( V((1^a),1) \right),
	$
	
	\item[$\bullet$] when $\g = \mf{b}_{m|n}$,  
	$
	{\bf T}^{\mf{b}_{m|n}}(a) \cong \ms{B}\left( V((1^a),2) \right)
	$ and
	$	 
	{\bf T}^{\texttt{sp}} \cong \ms{B}\left( V((0), 1) \right),
	$
	
	\item[$\bullet$] when $\g = \mf{d}_{m|n}$,
	\begin{equation*}
		\begin{cases}
		{\bf T}^{\mf{d}_{m|n}}(a) \cong \ms{B}\left( V((1^a),2) \right), \\
		\ov{\bf T}(0) \cong \ms{B}\left( V((2), 2) \right),
	\end{cases}
	\quad
	\begin{cases}
		{\bf T}^{\texttt{sp}+} \cong \ms{B}\left( V((0),1) \right), \\
		{\bf T}^{\texttt{sp}-} \cong \ms{B}\left( V((1),1) \right).
	\end{cases}
	\end{equation*}
\end{itemize}
Here $(0)$ denotes the empty partition (that is, $(0,0,\dots)$) and $(1^0)$ is understood as $(0)$.

For partition $\la=(\la_i)_{i \ge 1}$, let $\la'=(\la'_i)_{i \ge 1}$ be its conjugate.
For $\g=\mf{d}_{m|n}$ and given $(\lambda, \ell) \in \ms{P}(\mf{d}_{m|n})$,
let $\ov{\la} = (\ov{\la}_i)_{i \ge 1}$ be given by $\ov{\la}_1 = \ell - \la_1$ and $\ov{\la}_i = \la_i$ for $i \ge 2$. Let
$\mu = \la'$, $\ov{\mu} = \left( \ov{\la} \right)'$, $M_+ = \la_1$, $M_- = \ov{\la}_1$, 
and let $q_{\pm}$ and $r_{\pm}\in \left\{ 0, 1 \right\}$ be such that
\begin{equation*}
	\begin{cases}
		\ell - 2\la_1 = 2q_+ + r_+ & \text{if $\ell - 2\la_1 \ge 0$,} \\
		2\la_1 - \ell = 2q_- + r_- & \text{if $\ell - 2\la_1 < 0$},
	\end{cases}
\end{equation*}
We define $\widehat{\bf T}^\g(\lambda, \ell)$ as follows:
\begin{itemize}
 \item[$\bullet$] when $\g = \mf{b}_{m|n}$ or $\mf{c}_{m|n}$, 
\begin{equation*}
	\widehat{\bf T}^\g(\lambda, \ell) =
	\begin{cases}
		{\bf T}^\g(\la'_1) \times \dots \times {\bf T}^\g(\la'_\ell) & \text{if $\g = \mf{c}_{m|n}$,} \\
		{\bf T}^\g(\la'_1) \times \dots \times {\bf T}^\g(\la'_{\ell/2}) & \text{if $\g = \mf{b}_{m|n}$ and $\ell - 2\la_1$ is even,} \\
		{\bf T}^\g(\la'_1) \times \dots \times {\bf T}^\g(\la'_{(\ell-1)/2}) \times {\bf T}^{\texttt{sp}} & \text{if $\g = \mf{b}_{m|n}$ and $\ell - 2\la_1$ is odd,} \\
	\end{cases}
\end{equation*}

 \item[$\bullet$] when $\g = \mf{d}_{m|n}$,
 \begin{equation*}
	\widehat{\bf T}^\g(\lambda, \ell) =
	\begin{cases}
		{\bf T}^\g(\mu_1) \times \dots \times {\bf T}^\g(\mu_{M_+}) \times \left( {\bf T}(0) \right)^{q_+} \times \left( {\bf T}^{\texttt{sp}+} \right)^{r_+} & \text{if $\ell - 2\la_1 \ge 0$,} \\
		{\bf T}^\g(\ov{\mu}_1) \times \dots \times {\bf T}^\g(\ov{\mu}_{M_-}) \times \left( \ov{\bf T}(0) \right)^{q_-} \times \left( {\bf T}^{\texttt{sp}-} \right)^{r_-} & \text{if $\ell - 2\la_1 < 0$,}
	\end{cases}
\end{equation*}
\end{itemize}
where one may define a $\g$-crystal structure on $\widehat{\bf T}^\g(\lambda, \ell)$ by identifying $(\dots, T_2, T_1)$ with $T_1 \otimes T_2 \otimes \dots $.
Let 
\begin{equation*}
	{\bf T}^\g(\lambda, \ell) = 
	\left\{ \left. (\dots, T_2, T_1) \in \widehat{\bf T}^\g(\lambda, \ell) \, \right| \, T_{i+1} \prec T_i \text{ for $i \ge 1$ } \right\},
\end{equation*}
(see \cite[Definition 3.2]{K18-2} for the definition of $T_{i+1} \prec T_i$). 
Then ${\bf T}^\g(\lambda, \ell)$ is a $\g$-subcrystal of $\widehat{\bf T}^\g(\lambda, \ell)$, and it is shown in \cite{K15, K16} that 
\begin{equation*}
	{\bf T}^\g(\lambda, \ell) \cong \ms{B}(V(\la,\ell)).
\end{equation*}

\begin{ex}\label{ex: running ex 1}
{\rm
(1) Let $\g=\mf{c}_{2|5}$ and take $\la = (3, 2, 1, 1) \in \cP_{2|5}$ and $\ell = 3$.
Since $\ell - \la_1= 0$, we know that $(\la, \ell) \in \cP(\mf{c}_{2|5})$.
Let ${\bf T}=(T_3,T_2,T_1)$ given by
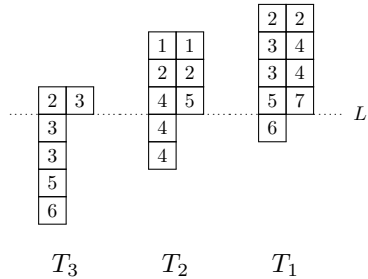
\begin{equation*}
    \ytableausetup {mathmode, boxsize=1.0em} 
    \begin{ytableau}
    \none & \none & \none & \none & \none & \none & \none & \none & \none & \tl{$2$} & \tl{$2$} & \none \\
    \none & \none & \none & \none & \none & \tl{$1$} & \tl{$1$} & \none & \none & \tl{$3$} & \tl{$4$} & \none \\
    \none & \none & \none & \none & \none & \tl{$2$} & \tl{$2$} & \none & \none & \tl{$3$} & \tl{$4$} & \none \\
	\none[\mathrel{\raisebox{-0.7ex}{$\scalebox{0.45}{\dots\dots}$}}] & \tl{$2$} & \tl{$3$} & 
	\none[\mathrel{\raisebox{-0.7ex}{$\scalebox{0.45}{\dots\dots}$}}] &  
	\none[\mathrel{\raisebox{-0.7ex}{$\scalebox{0.45}{\dots\dots}$}}] & \tl{$4$} & \tl{$5$} & 
	\none[\mathrel{\raisebox{-0.7ex}{$\scalebox{0.45}{\dots\dots}$}}] & 
	\none[\mathrel{\raisebox{-0.7ex}{$\scalebox{0.45}{\dots\dots}$}}] & \tl{$5$} & \tl{$7$} & 
	\none[\ \ \ \mathrel{\raisebox{-0.7ex}{$\scalebox{0.45}{\dots\dots}$\ ${}_{\scalebox{0.75}{$L$}}$}}] \\
    \none & \tl{$3$} & \none & \none & \none & \tl{$4$} & \none & \none & \none & \tl{$6$} & \none & \none \\
    \none & \tl{$3$} & \none & \none & \none & \tl{$4$} & \none & \none & \none & \none & \none & \none & \none \\
    \none & \tl{$5$} & \none & \none & \none & \none & \none & \none & \none & \none & \none & \none \\
    \none & \tl{$6$} & \none & \none & \none & \none & \none & \none & \none & \none & \none & \none \\
    \none & \none & \none & \none & \none & \none & \none & \none & \none & \none & \none & \none \\
    \none & \none[\quad T_3] & \none & \none & \none & \none[\quad T_2] & \none & \none & \none & \none[\quad T_1] & \none & \none \\
    \end{ytableau}
\end{equation*}
\vskip 2mm \noindent
    where the tableaux $T_1,\cdots,T_3$ are placed on the plane with a horizontal line $L$ drawn as the dotted line. 
 Let us check that $T_3 \prec T_2$. We have
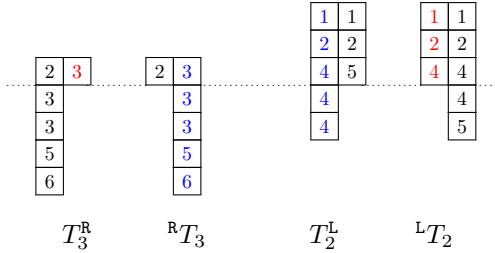
\begin{equation*}
    \begin{ytableau}
    \none & \none & \none & \none & \none & \none & \none & \none & \none & \none & \none & \tl{$\blue{1}$} & \tl{$1$} & \none & \none & \tl{$\red{1}$} & \tl{$1$} & \none \\
    \none & \none & \none & \none & \none & \none & \none & \none & \none & \none & \none & \tl{$\blue{2}$} & \tl{$2$} & \none & \none & \tl{$\red{2}$} & \tl{$2$} & \none \\
    \none[\mathrel{\raisebox{-0.7ex}{$\scalebox{0.45}{\dots\dots}$}}] & \tl{$2$} & \tl{$\red{3}$} & 
	\none[\ \ \ \mathrel{\raisebox{-0.7ex}{$\scalebox{0.45}{\dots\dots\dots\dots}$}}] & \none & \tl{$2$} & \tl{$\blue{3}$} & 
	\none[\ \ \ \ \ \ \ \ \ \ \mathrel{\raisebox{-0.7ex}{$\scalebox{0.45}{\dots\dots\dots\dots\dots\dots\dots}$}}] & \none & \none & \none & \tl{$\blue{4}$} & \tl{$5$} & 
	\none[\ \ \ \mathrel{\raisebox{-0.7ex}{$\scalebox{0.45}{\dots\dots\dots\dots}$}}] & \none & \tl{$\red{4}$} & \tl{$4$} & 
	\none[\mathrel{\raisebox{-0.7ex}{$\scalebox{0.45}{\dots\dots}$}}] \\
    \none & \tl{$3$} & \none & \none & \none & \none & \tl{$\blue{3}$} & \none & \none & \none & \none & \tl{$\blue{4}$} & \none & \none & \none & \none & \tl{$4$} & \none \\
    \none & \tl{$3$} & \none & \none & \none & \none & \tl{$\blue{3}$} & \none & \none & \none & \none & \tl{$\blue{4}$} & \none & \none & \none & \none & \tl{$5$} & \none \\
    \none & \tl{$5$} & \none & \none & \none & \none & \tl{$\blue{5}$} & \none & \none & \none & \none & \none & \none & \none & \none & \none & \none & \none \\
    \none & \tl{$6$} & \none & \none & \none & \none & \tl{$\blue{6}$} & \none & \none & \none & \none & \none & \none & \none & \none & \none & \none & \none \\
    \none & \none & \none & \none & \none & \none & \none & \none & \none & \none & \none & \none & \none & \none & \none & \none & \none & \none \\    
    \none & \none & \none[T_3^{\tt R}] & \none & \none & \none & \none[^{\tt R}T_3] & \none & \none & \none & \none & \none[T_2^{\tt L}] & \none & \none & \none & \none[^{\tt L}T_2] & \none & \none \\    
\end{ytableau}
\end{equation*}
\noindent (see \cite[Section 3.2]{K18-2} for the definition of $T^{\tt R}, T^{\tt L}, ^{\tt R}T$ and $^{\tt L}T$). 
Then $({}^{\tt R}T_3,T^{\tt L}_2)$ (in blue) and $(T^{\tt R}_3,{}^{\tt L}T_2)$ (in red) form semistandard tableaux
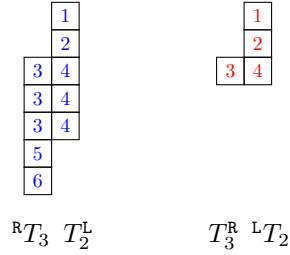
\begin{equation*}
    \begin{ytableau}
    \none & \tl{$\blue{1}$} & \none & \none & \none & \none & \none & \none & \tl{$\red{1}$}\\
    \none & \tl{$\blue{2}$} & \none & \none & \none & \none & \none & \none & \tl{$\red{2}$}\\
    \tl{$\blue{3}$} & \tl{$\blue{4}$} & \none & \none & \none & \none & \none & \tl{$\red{3}$} & \tl{$\red{4}$}\\
    \tl{$\blue{3}$} & \tl{$\blue{4}$} & \none & \none & \none & \none & \none & \none & \none\\
    \tl{$\blue{3}$} & \tl{$\blue{4}$} & \none & \none & \none & \none & \none & \none & \none\\
    \tl{$\blue{5}$} & \none & \none & \none & \none & \none & \none & \none & \none\\
    \tl{$\blue{6}$} & \none & \none & \none & \none & \none & \none & \none & \none\\
    \none & \none & \none & \none & \none & \none & \none & \none & \none\\    
    \none[\! \! \! ^{\tt R}T_3] & \none[\ \ \ T_2^{\tt L}] & \none & \none & \none & \none & \none & \none[\! \! \! T_3^{\tt R}] & \none[\ \ \ ^{\tt L}T_2]\\    
\end{ytableau}
\end{equation*}
\noindent In similar way, we can also check that $T_2 \prec T_1$, and hence ${\bf T} \in {\bf T}^{\g}(\lambda,\ell)$. \\
(2) Let $\g = \mf{d}_{4|4}$.
    Take $\la = (4, 4, 2) \in \cP_{4|4}$ and $\ell = 8$ so that $\ell - 2\la_1 = 0$ and $\mu = \la' = (3, 3, 2, 2)$. 
	Since $\ell - \la_1 - \la_2 = 0$ and $\ell - 2\la_1 = 0$, we know that $(\la, \ell) \in \cP(\mf{d}_{4|4})$ and $q_+ = r_+ = 0$.     
    Let ${\bf T} = (T_4, T_3, T_2, T_1)$ given by
    \begin{equation*}
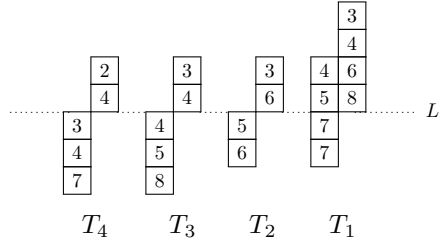
    \begin{split}
    & \begin{ytableau}
    \none &\none & \none & \none & \none & \none & \none & \none & \none & \none & \none & \tl{$3$} & \none  \\
    \none &\none & \none & \none & \none & \none & \none & \none & \none & \none & \none & \tl{$4$} & \none  \\
    \none &\none & \tl{$2$} & \none & \none & \tl{$3$} & \none & \none & \tl{$3$} & \none & \tl{$4$} & \tl{$6$} & \none \\
    \none[\!\!\!\!\mathrel{\raisebox{-0.7ex}{$\scalebox{0.45}{\dots\dots\dots\dots}$}}] &\none & \tl{$4$} & \none[\mathrel{\raisebox{-0.7ex}{$\scalebox{0.45}{\dots\dots}$}}] & \none & \tl{$4$} & \none[\mathrel{\raisebox{-0.7ex}{$\scalebox{0.45}{\dots\dots}$}}] & \none & \tl{$6$} & \none[\mathrel{\raisebox{-0.7ex}{$\scalebox{0.45}{\dots\dots}$}}]& \tl{$5$} & \tl{$8$} & \none[\ \mathrel{\raisebox{-0.7ex}{\quad $\scalebox{0.45}{\dots\dots\dots\dots}$\ ${}_{\scalebox{0.75}{$L$}}$}}] \\
    \none &\tl{$3$} & \none & \none & \tl{$4$} & \none & \none & \tl{$5$} & \none & \none & \tl{$7$} & \none & \none  \\
    \none &\tl{$4$} & \none & \none & \tl{$5$} & \none & \none & \tl{$6$} & \none & \none & \tl{$7$} & \none & \none  \\
    \none &\tl{$7$} & \none & \none & \tl{$8$} & \none & \none & \none & \none & \none & \none & \none & \none  \\
    \end{ytableau} \quad\quad \\
    &\ \hskip 5mm T_4 \hskip 8mm T_3 \hskip 7mm T_2 \hskip 7mm T_1 \hskip 13mm 
    \end{split}
    \end{equation*}
	One can check that $T_4 \prec T_3 \prec T_2 \prec T_1$, and hence ${\bf T} \in {\bf T}^{\g}(\lambda,\ell)$. Note that $T_4 \ntriangleleft T_3 \triangleleft T_2 \triangleleft T_1$, where the definition of $\triangleleft$ is given in \cite[Definition 3.5]{JK19}.
} 
\end{ex}

\subsection{Crystal $\ms{B}(P(\la))$}
Let us give some comments on the combinatorial realization of $\ms{B}(P(\la))$.
By \eqref{eq: crystal base of Pla}, it suffices to describe $\ms{B}(\mc{N})$. We may regard $\ms{B}(\mc{N})= M(\mf{u})$ by Theorem \ref{thm:crystal base of N}, and also identify it with the set ${\bf M}^\g$ of upper-triangular matrices ${\bf c}=\left(c_{ij}\right)_{i, j \in \I, i \le j}$ such that
\begin{equation*}
	c_{ij} \in
	\begin{cases}
		\Z_{\ge 0} & \text{if $\epsilon_i = \epsilon_j$,} \\
		\left\{ 0, 1 \right\} & \text{if $\epsilon_i \neq \epsilon_j$,}
	\end{cases}
	\quad \text{ and } \quad
	c_{ii} =
	\begin{cases}
		0 & \text{if $\g = \mf{c}_{m|n}$ and $i \in \I_{\ov{1}}$,} \\
		0 & \text{if $\g = \mf{d}_{m|n}$ and $i \in \I_{\ov{0}}$,}
	\end{cases}
\end{equation*}
where $c_{ij}$ is equal to the multiplicity $c_\beta=c_{(i,j)}$ of $\bfF_\beta$ with $\beta$ corresponding to $(i,j)$ following \eqref{eq:nilradical roots}.
For example, when $\g = \mf{b}_{2|3}$,
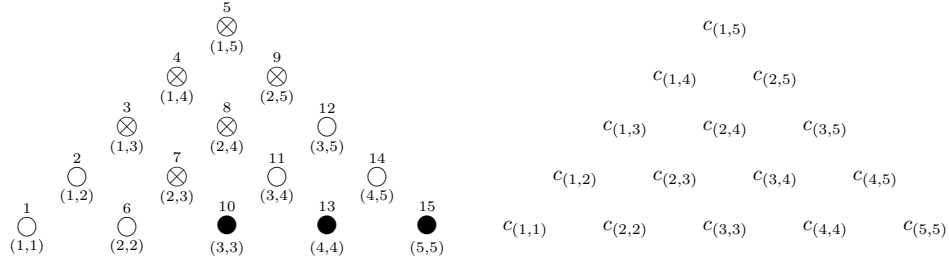
\begin{equation*}
\begin{tikzpicture}[baseline=(current  bounding  box.center), every node/.style={scale=0.8}, scale=0.95]
		\node (glw_00f) at (0,0) {$\underset{(1,1)}{\overset{\tiny 1}{\bigcirc}}$};
		\node (glw_01m) at (0.7,0.7) {$\underset{(1,2)}{\overset{\tiny 2}{\bigcirc}}$};
		\node (glw_02m) at (1.4,1.4) {$\underset{(1,3)}{\overset{\tiny 3}{\scalebox{1.4}{$\ot$}}}$};
		\node (glw_03m) at (2.1,2.1) {$\underset{(1,4)}{\overset{\tiny 4}{\scalebox{1.4}{$\ot$}}}$};
		\node (glw_04m) at (2.8,2.8) {$\underset{(1,5)}{\overset{\tiny 5}{\scalebox{1.4}{$\ot$}}}$};
		\node (glw_01f) at (1.4,0) {$\underset{(2,2)}{\overset{\tiny 6}{\bigcirc}}$};
		\node (glw_12m) at (2.1,0.7) {$\underset{(2,3)}{\overset{\tiny 7}{\scalebox{1.4}{$\ot$}}}$};
		\node (glw_13m) at (2.8,1.4) {$\underset{(2,4)}{\overset{\tiny 8}{\scalebox{1.4}{$\ot$}}}$};
		\node (glw_14m) at (3.5,2.1) {$\underset{(2,5)}{\overset{\tiny 9}{\scalebox{1.4}{$\ot$}}}$};
		\node (glw_02f) at (2.8,0) {$\underset{(3,3)}{\overset{\tiny 10}{{\scalebox{2.3}{$\bullet$}}}}$};
		\node (glw_23m) at (3.5,0.7) {$\underset{(3,4)}{\overset{\tiny 11}{\bigcirc}}$};		
		\node (glw_24m) at (4.2,1.4) {$\underset{(3,5)}{\overset{\tiny 12}{\bigcirc}}$};
		\node (glw_03f) at (4.2,0) {$\underset{(4,4)}{\overset{\tiny 13}{{\scalebox{2.3}{$\bullet$}}}}$};
		\node (glw_34m) at (4.9,0.7) {$\underset{(4,5)}{\overset{\tiny 14}{\bigcirc}}$};
		\node (glw_04m) at (5.6,0) {$\underset{(5,5)}{\overset{\tiny 15}{{\scalebox{2.3}{$\bullet$}}}}$};
		\end{tikzpicture}
		\,\,\quad\,\,
		\begin{tikzpicture}[baseline=(current  bounding  box.center), every node/.style={scale=0.8}, scale=0.95]
		\node (glw_00f) at (0,0) {$c_{(1,1)}$};
		\node (glw_01m) at (0.7,0.7) {$c_{(1,2)}$};
		\node (glw_02m) at (1.4,1.4) {$c_{(1,3)}$};
		\node (glw_03m) at (2.1,2.1) {$c_{(1,4)}$};
		\node (glw_04m) at (2.8,2.8) {$c_{(1,5)}$};
		\node (glw_01f) at (1.4,0) {$c_{(2,2)}$};
		\node (glw_12m) at (2.1,0.7) {$c_{(2,3)}$};
		\node (glw_13m) at (2.8,1.4) {$c_{(2,4)}$};
		\node (glw_14m) at (3.5,2.1) {$c_{(2,5)}$};
		\node (glw_02f) at (2.8,0) {$c_{(3,3)}$};
		\node (glw_23m) at (3.5,0.7) {$c_{(3,4)}$};		
		\node (glw_24m) at (4.2,1.4) {$c_{(3,5)}$};
		\node (glw_03f) at (4.2,0) {$c_{(4,4)}$};
		\node (glw_34m) at (4.9,0.7) {$c_{(4,5)}$};
		\node (glw_04m) at (5.6,0) {$c_{(5,5)}$};
	\end{tikzpicture}
\end{equation*}
and when $\g = \mf{d}_{2|4}$,
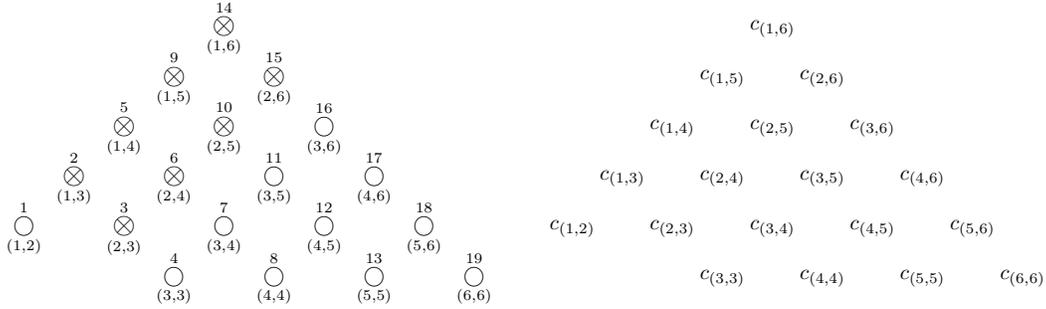
\begin{equation*}
\begin{tikzpicture}[baseline=(current  bounding  box.center), every node/.style={scale=0.8}, scale=0.95]
		\node (glw_01m) at (0.7,0.7) {$\underset{(1,2)}{\overset{\tiny 1}{\bigcirc}}$};
		\node (glw_02m) at (1.4,1.4) {$\underset{(1,3)}{\overset{\tiny 2}{\scalebox{1.4}{$\ot$}}}$};
		\node (glw_03m) at (2.1,2.1) {$\underset{(1,4)}{\overset{\tiny 5}{\scalebox{1.4}{$\ot$}}}$};
		\node (glw_04m) at (2.8,2.8) {$\underset{(1,5)}{\overset{\tiny 9}{\scalebox{1.4}{$\ot$}}}$};
		\node (glw_04m) at (3.5,3.5) {$\underset{(1,6)}{\overset{\tiny 14}{\scalebox{1.4}{$\ot$}}}$};
		\node (glw_12m) at (2.1,0.7) {$\underset{(2,3)}{\overset{\tiny 3}{\scalebox{1.4}{$\ot$}}}$};
		\node (glw_13m) at (2.8,1.4) {$\underset{(2,4)}{\overset{\tiny 6}{\scalebox{1.4}{$\ot$}}}$};
		\node (glw_14m) at (3.5,2.1) {$\underset{(2,5)}{\overset{\tiny 10}{\scalebox{1.4}{$\ot$}}}$};
		\node (glw_14m) at (4.2,2.8) {$\underset{(2,6)}{\overset{\tiny 15}{\scalebox{1.4}{$\ot$}}}$};
		\node (glw_02f) at (2.8,0) {$\underset{(3,3)}{\overset{\tiny 4}{\bigcirc}}$};
		\node (glw_23m) at (3.5,0.7) {$\underset{(3,4)}{\overset{\tiny 7}{\bigcirc}}$};		
		\node (glw_24m) at (4.2,1.4) {$\underset{(3,5)}{\overset{\tiny 11}{\bigcirc}}$};
		\node (glw_24m) at (4.9,2.1) {$\underset{(3,6)}{\overset{\tiny 16}{\bigcirc}}$};
		\node (glw_03f) at (4.2,0) {$\underset{(4,4)}{\overset{\tiny 8}{\bigcirc}}$};
		\node (glw_34m) at (4.9,0.7) {$\underset{(4,5)}{\overset{\tiny 12}{\bigcirc}}$};
		\node (glw_34m) at (5.6,1.4) {$\underset{(4,6)}{\overset{\tiny 17}{\bigcirc}}$};
		\node (glw_04m) at (5.6,0) {$\underset{(5,5)}{\overset{\tiny 13}{\bigcirc}}$};
		\node (glw_04m) at (6.3,0.7) {$\underset{(5,6)}{\overset{\tiny 18}{\bigcirc}}$};
		\node (glw_04m) at (7.0,0) {$\underset{(6,6)}{\overset{\tiny 19}{\bigcirc}}$};
		\end{tikzpicture}
		\,\,\quad\,\,
		\begin{tikzpicture}[baseline=(current  bounding  box.center), every node/.style={scale=0.8}, scale=0.95]
		\node (glw_01m) at (0.7,0.7) {$c_{(1,2)}$};
		\node (glw_02m) at (1.4,1.4) {$c_{(1,3)}$};
		\node (glw_03m) at (2.1,2.1) {$c_{(1,4)}$};
		\node (glw_04m) at (2.8,2.8) {$c_{(1,5)}$};
		\node (glw_04m) at (3.5,3.5) {$c_{(1,6)}$};
		\node (glw_12m) at (2.1,0.7) {$c_{(2,3)}$};
		\node (glw_13m) at (2.8,1.4) {$c_{(2,4)}$};
		\node (glw_14m) at (3.5,2.1) {$c_{(2,5)}$};
		\node (glw_14m) at (4.2,2.8) {$c_{(2,6)}$};
		\node (glw_02f) at (2.8,0) {$c_{(3,3)}$};
		\node (glw_23m) at (3.5,0.7) {$c_{(3,4)}$};		
		\node (glw_24m) at (4.2,1.4) {$c_{(3,5)}$};
		\node (glw_24m) at (4.9,2.1) {$c_{(3,6)}$};
		\node (glw_03f) at (4.2,0) {$c_{(4,4)}$};
		\node (glw_34m) at (4.9,0.7) {$c_{(4,5)}$};
		\node (glw_34m) at (5.6,1.4) {$c_{(4,6)}$};
		\node (glw_04m) at (5.6,0) {$c_{(5,5)}$};
		\node (glw_04m) at (6.3,0.7) {$c_{(5,6)}$};
		\node (glw_04m) at (7.0,0) {$c_{(6,6)}$};
	\end{tikzpicture}
\end{equation*}
Here the notation $\underset{\beta}{\overset{k}{\scalebox{1.3}{$\odot$}}}$ means the $k$-th root $\beta$ with respect to $\prec$ and $\scalebox{1.3}{$\odot$}$ denotes the type of $\beta$, where $\scalebox{1.3}{$\odot$} = \bigcirc$ (even), $\ot$ (isotropic) or $\raisebox{-0.1cm}{{\scalebox{2.3}{$\bullet$}}}$ (non-isotropic odd).
An explicit description of $\te_i$ and $\tf_i$ on ${\bf M}^\g$ can be found in \cite[Section 4.2]{JKU24}. 

For $\la\in \cP_{m|n}$, we set
\begin{equation} \label{eq: def of Mla}
	{\bf M}^\g_\la = {\bf M}^\g \times SST_{m|n}(\la),
\end{equation}
which we may identify with $\ms{B}(P(\la))$.

\subsection{RSK of BCD types}
Let us give an alternative description of ${\bf M}^\g$ in terms of $\I$-semistandard tableaux. Let  
\begin{equation*} \label{eq:def of Vg and Vgla}
	{\bf V}^\g  = \bigsqcup_{\delta \in \ms{P}_\g} SST_{m|n}\left( \delta^\pi \right), 
\end{equation*}
where
\begin{equation*}
	\ms{P}_\g = 
	\begin{cases}
		\ms{P}_{m|n} & \text{if $\g = \mf{b}_{m|n}$,} \\
		\left\{ \la=(\la_i)_{i\ge 1} \in \ms{P}_{m|n} \, | \, \la_i \in 2 \Z_{\ge 0} \right\}  & \text{if $\g = \mf{c}_{m|n}$,} \\
		\left\{ \la=(\la_i)_{i\ge 1} \in \ms{P}_{m|n} \, | \, \la'_i \in 2 \Z_{\ge 0} \right\}  & \text{if $\g = \mf{d}_{m|n}$.} \\
	\end{cases}
\end{equation*}
Then ${\bf V}^\g$ is a $\g$-crystal, where the $\mf{l}$-crystal structure is given as in Section \ref{subsec: type A crystal}, and the actions of $\te_0$ and $\tf_0$ are defined as in \cite[Remark 4.1]{K18} for $\g=\mf{b}_{m|n}$, $\mf{c}_{m|n}$, and as in \cite[Section 5.2]{K13} for $\g=\mf{d}_{m|n}$.
For example, $\tf_0$ acts on $T\in {\bf V}^\g$ by adding
\begin{equation*}
\begin{cases}
\ \scalebox{0.8}{$\boxed{1}$} & \text{if $\g = \mf{b}_{m|n}$,} \\
\ \scalebox{1.1}{$\hdomino$} & \text{if $\g = \mf{c}_{m|n}$,} \\
\ \scalebox{1.1}{$\vdomino$} & \text{if $\g = \mf{d}_{m|n}$,}
\end{cases}
\end{equation*}
on $T$ at certain corner according to a signature rule.
In this subsection, we give a combinatorial description of an isomorphism of $\g$-crystals
\begin{equation}\label{eq: iso from M to V}
 {\bf M}^\g\cong {\bf V}^\g.
\end{equation}

\subsubsection{The case when $\g=\mf{b}_{m|n}$, $\mf{c}_{m|n}$}
In this case, the isomorphism \eqref{eq: iso from M to V} is obtained from the RSK correspondence on symmetric matrices as an $\mf{l}$-crystal isomorphism \cite{K07} and then extending it to an isomorphism of $\g$-crystals. The proof is essentially the same as in \cite[Section 4.3]{K18}. So we briefly explain how to construct an isomorphism without proof.

Let 
\begin{equation*}\label{eq: def of Mmn}
{\bf M}= 
\left\{\, {\bf m}=\left( m_{ij} \right)_{i,j \in \I} \Big|\,\, m_{ij} \in \Z_{\ge 0} \  (\epsilon_i = \epsilon_j),\  m_{ij} \in \left\{ 0, 1 \right\}\  (\epsilon_i \neq \epsilon_j)\, \right\}.
\end{equation*}
Then it has an $(\mf{l},\mf{l})$-bicrystal structure, where we denote by $\tte_i$ and $\ttf_i$ (resp. $\tte^*_i$ and $\ttf^*_i$) the crystal operators for the first (resp. second) $\mf{l}$ \cite{K07}.

Let ${\bf m}\in {\bf M}$ be given. Let ${\bf m}^j$ denote the $j$-th column of ${\bf m}$ and $m_j=\sum_{i}m_{ij}$. We identify ${\bf m}^j$ with an element in $SST_{m|n}((m_j))$ (resp. $SST_{m|n}((1^{m_j}))$) when $\e_j=0$ (resp. $\e_j=1$), and regard ${\bf m}={\bf m}^1\ot \dots \ot {\bf m}^{m+n}$ as an element of an $\mf{l}$-crystal with respect to $\tte_i$ and $\ttf_i$ ($i\in I$).  
Then $\tte^*_i$ and $\ttf^*_i$ are given by $\tte^*_i{\bf m} =\left(\tte_i {\bf m}^t\right)^t$ and $\ttf^*_i{\bf m} =\left(\ttf_i {\bf m}^t\right)^t$, where ${\bf m}^t$ denotes the transpose of ${\bf m}$.

Let $w_1\dots w_{\ell}=w^1 \cdots w^{m+n}$ be the word with letters $w_k\in \I$ ($1\le k\le \ell:=\sum_{i,j}m_{ij})$, where $w^j=w({\bf m}^j)$ is the word of ${\bf m}^j$ regarding ${\bf m}^j$ as an element in $SST_{m|n}((m_j))$ or $SST_{m|n}((1^{m_j}))$.
Let 
\begin{equation*}
 P({\bf m})=(((w_{\ell}\leftarrow w_{\ell-1})\cdots)\leftarrow w_1),\quad Q({\bf m})=P({\bf m}^t).
\end{equation*}
Then the following RSK correspondence
\begin{equation} \label{eq:RSK}
\begin{split}
	\kappa :
\xymatrixcolsep{3.5pc}\xymatrixrowsep{0pc}\xymatrix{
	{\bf M} \ar@{->}[r] & \displaystyle \bigsqcup_{\la \in \ms{P}_{m|n}} SST_{m|n}(\la^\pi) \times SST_{m|n}(\la^\pi) \\
	{\bf m} \ar@{|->}[r] & (P({\bf m}),Q({\bf m})) 
	}
\end{split}
\end{equation}
is an isomorphism of $(\mf{l},\mf{l})$-bicrystals \cite[Theorem 3.11]{K07}, where $\tte_i$ and $\ttf_i$ (resp. $\tte^*_i$ and $\ttf^*_i$) act on $P({\bf m})$  (resp. $Q({\bf m})$) on the right-hand side.

Now, we regard as a set
\begin{equation}\label{eq: M^g into sym}
{\bf M}^\g\subset {\bf M}^{\rm sym}:=\{\,{\bf m}\in {\bf M}\,|\,{\bf m}={\bf m}^t\,\} \quad (\g=\mf{b}_{m|n}, \mf{c}_{m|n}) 
\end{equation}
by identifying ${\bf c}=(c_{ij})$ with ${\bf m}=(m_{ij})$ where $m_{ii}=2c_{ii}/r_\g$ and $m_{ij}=c_{ij}$ ($i\neq j$). Under \eqref{eq: M^g into sym}, we can check by the $\mf{l}$-crystal structure on ${\bf M}^{\g}=\ms{B}(\mc{N})$ \cite[Section 4.2]{JKU24} that
$\te_i=\tte_i\tte_i^*$ and $\tf_i=\ttf_i\ttf_i^*$ $(i\in I\setminus\{0\})$.

Then we have the following which proves \eqref{eq: iso from M to V}.

\begin{thm}\label{thm: iso from M to V for BC}
The restriction of \eqref{eq:RSK} on ${\bf M}^\g$ is an isomorphism of $\g$-crystals
\begin{equation*}\label{eq:folding of RSK}
\begin{split}
\kappa^\g : 
\xymatrix@R=0em{
{\bf M}^\g \ar@{->}[r] & \displaystyle {\bf V}^\g \\ 
{\bf c} \ar@{|->}[r] & P({\bf m}),
}
\end{split}
\end{equation*}
where ${\bf m}\in {\bf M}^{\rm sym}$ corresponds to ${\bf c}$ under \eqref{eq: M^g into sym}
\end{thm}
\pf It follows from almost the same argument as in \cite[Theorem 4.9]{K18}.
\qed

\begin{ex}\label{ex: running ex 3}
{\rm
Let $\mf{g}=\mf{c}_{2|5}$. The matrix ${\bf m}\in {\bf M}^{\rm sym}$ and $\tf_0({\bf m})$ are given by
\begin{equation*}
{\bf m}=	
   \begin{bmatrix} 
   0 & 2 & 0 & 0 & 0 & 0 & 0 \\
   2 & 0 & 1 & 1 & 0 & 0 & 1 \\
   0 & 1 & 0 & 1 & 0 & 1 & 0 \\
   0 & 1 & 1 & 0 & 1 & 0 & 0 \\
   0 & 0 & 0 & 1 & 0 & 0 & 0 \\
   0 & 0 & 1 & 0 & 0 & 0 & 0 \\
   0 & 1 & 0 & 0 & 0 & 0 & 0
   \end{bmatrix} \quad , \quad
\tf_0 {\bf m}=	
   \begin{bmatrix} 
   \red{\bf 2} & 2 & 0 & 0 & 0 & 0 & 0 \\
   2 & 0 & 1 & 1 & 0 & 0 & 1 \\
   0 & 1 & 0 & 1 & 0 & 1 & 0 \\
   0 & 1 & 1 & 0 & 1 & 0 & 0 \\
   0 & 0 & 0 & 1 & 0 & 0 & 0 \\
   0 & 0 & 1 & 0 & 0 & 0 & 0 \\
   0 & 1 & 0 & 0 & 0 & 0 & 0
   \end{bmatrix} \quad .
\end{equation*}
Then $P({\bf m})$ and $P(\tf_0{\bf m})$ are given by
\begin{equation*}
    \begin{ytableau}
		\none & \none & \none & \none & \tl{$2$} & \tl{$2$} \\
		\none & \none & \tl{$1$} & \tl{$1$} & \tl{$3$} & \tl{$4$} \\
		\none & \none & \tl{$2$} & \tl{$2$} & \tl{$3$} & \tl{$4$} \\
		\tl{$2$} & \tl{$3$} & \tl{$4$} & \tl{$5$} & \tl{$6$} & \tl{$7$} \\
		\none & \none & \none & \none & \none & \none \\
		\none & \none & \none & \none[P({\bf m})] & \none & \none \\
    \end{ytableau}
	\quad \qquad \quad
    \begin{ytableau}
		\none & \none & \none & \none & \tl{$2$} & \tl{$2$} \\
		\none & \none & \tl{$1$} & \tl{$1$} & \tl{$3$} & \tl{$4$} \\
		\tl{$\red{\bf 1}$} & \tl{$\red{\bf 1}$} & \tl{$2$} & \tl{$2$} & \tl{$3$} & \tl{$4$} \\
		\tl{$2$} & \tl{$3$} & \tl{$4$} & \tl{$5$} & \tl{$6$} & \tl{$7$} \\
		\none & \none & \none & \none & \none & \none \\
		\none & \none & \none & \none[P(\tf_0{\bf m})] & \none & \none \\
    \end{ytableau}
\end{equation*}
\noindent where we can also see that $\tf_0(P({\bf m}))=P(\tf_0{\bf m})$.
} 
\end{ex}

\subsubsection{The case when $\g=\mf{d}_{m|n}$}
In this case, we give a crystal theoretic super-analogue of the bijection of type $D$ \cite{B}, which proves \eqref{eq: iso from M to V}.

We first remark that the set ${\bf M}^{\mf{d}_{m|n}}$ is in bijection with the following set of biwords with with letters in $\I$;
\begin{equation} \label{eq: df of biws}
	\biws =
	\left\{ \,
				\left( {\bf i}, {\bf j} \right) \, \left| \,\, 
				\parbox{20em}{(1)\ ${\bf i} = i_1\cdots i_r$, ${\bf j} = j_1\cdots j_r$ for some $r \ge 0$  \\
				(2)\, $(i_1, j_1) \le' \dots \le' (i_r, j_r)$  \\
				(3)\, $\epsilon_{i_k} \neq \epsilon_{j_k}$ implies $(i_k, j_k) \neq (i_{k\pm1}, j_{k\pm1})$ }
			\, \right.
	\right\},
\end{equation}
where $(i,j) <' (k,l)$ if either one of the following conditions is satisfied: (i) $i < k$, (ii) $i = k \in \I_{\ov{0}}$ and $j > l$, (iii) $i = k \in \I_{\ov{1}}$ and $j < l$. Here we identify $\left( {\bf i}, {\bf j} \right)\in \biws$ with ${\bf m}({\bf i}, {\bf j}) = (m_{ij}) \in {\bf M}^{\mf{d}_{m|n}}$, the unique element of ${\bf M}^{\mf{d}_{m|n}}$ such that $m_{ij} = \left| \left\{ k \, | \, (i_k, j_k) = (i, j) \right\} \right|$.

For a biletter $\binom{\,i\,}{\,j\,}$ and $T \in SST_{m|n}(\la^\pi)$ for $i,j \in \I$ and $\la \in \ms{P}_{\mf{d}_{m|n}}$, let
\begin{equation*}
 T \overset{\,\,{\rm B}}{\longleftarrow} \binom{\,i\,}{\,j\,}
\end{equation*}
be the tableau obtained by adding $\scalebox{0.8}{$\boxed{i}$}$ at the corner of $T \leftarrow j$ located above the box in the skew Young diagram ${\rm sh}\left(T \leftarrow j \right) / {\rm sh}\left( T \right)$. It is an analogue of the insertion in \cite{B} with respect to $\I$-semistandard tableaux.

For $({\bf i}, {\bf j}) \in \biws$ with ${\bf i} = i_1 \dots i_r$ and ${\bf j} = j_1 \dots j_r$, we define
\begin{equation}\label{df: type D insertion}
 P({\bf i}, {\bf j})=
 \left(\left(\left(\vd{i_{r}}{j_r}\overset{\,\,{\rm B}}{\longleftarrow} \binom{\,i_{r-1}\,}{\,j_{r-1}\,}\right)\cdots\right)\overset{\,\,{\rm B}}{\longleftarrow} \binom{\,i_{1}\,}{\,j_{1}\,}\right).
\end{equation}

\begin{ex} \label{ex:RSK for d23}
{\em 
Suppose that $m=4$ and $n=4$ and let $({\bf i}, {\bf j}) \in \Omega^{\mf{d}_{4|4}}$ be given by
\begin{equation*}
\begin{split}
	\left(
		\begin{array}{c}
			{\bf i} \\
			{\bf j} 
		\end{array}
	\right)
	&=
	\left(
		\begin{array}{ccccc|c}
			2 & 3 & 3 & 3 & 4 & 5 \\
			8 & 4 & 4 & 4 & 6 & 7
		\end{array}
	\right).
\end{split}
\end{equation*}
Let $P_k=P(i_k\dots i_6,j_k\dots j_6)$ for $k=1,\dots,6$. Then we have
\begin{equation*}
P_6= \vd{5}{7},\quad 
P_5=\left(P_6 \overset{\,\,{\rm B}}{\longleftarrow} \binom{\,4\,}{\,6\,}\right) = \scalebox{0.8}{\raisebox{0.36cm}{\!\!\xymatrix@R=-0.82em @C=-0.82em{ \boxed{4} & \boxed{6} \\ \boxed{5} & \boxed{7}}}},
\end{equation*}
since $P_6 \leftarrow 6 = \scalebox{0.8}{\raisebox{0.36cm}{\!\!\xymatrix@R=-0.82em @C=-0.82em{ & \boxed{6} \\ \boxed{5} & \boxed{7}}}}$, and by iterating the preceding procedure, we have that $P_1=P({\bf i}, {\bf j})$ is equal to
\begin{equation*}
\begin{ytableau}
    \none & \none & \none & \none & \tl{$3$}\\
    \none & \none & \none & \none & \tl{$4$}\\
    \tl{$2$} & \tl{$3$} & \tl{$3$} & \tl{$4$} & \tl{$6$}\\
    \tl{$4$} & \tl{$4$} & \tl{$5$} & \tl{$7$} & \tl{$8$}\\
\end{ytableau}
\end{equation*}
\\
and $P({\bf i}, {\bf j}) \in SST_{4|4}(\delta^\pi)$ where $\delta = (5,5,1,1) \in \ms{P}_{\mf{d}_{4|4}}$.
}
\end{ex}

Now we have the following, which proves \eqref{eq: iso from M to V}.
\begin{thm} \label{thm: Burge}
The map
\begin{equation*} \label{eq:Burge correspondence}
\begin{split}
\kappa^{\mf{d}_{m|n}} :
\xymatrix@R=0em{
	{\bf M}^{\mf{d}_{m|n}} \ar@{->}[r] & {\bf V}^{\mf{d}_{m|n}} \\
	{\bf m} \ar@{|->}[r] & P({\bf i}, {\bf j})
},
\end{split}
\end{equation*}
where ${\bf m} = {\bf m}({\bf i}, {\bf j})$ for $({\bf i}, {\bf j}) \in \biws$, is an isomorphism of $\mf{d}_{m|n}$-crystals.
\end{thm}

When $n=0$, it is proved in \cite[Theorem 4.6]{JK19}.
But unlike the case of $\mf{b}_{m|n}$ and $\mf{c}_{m|n}$, the result  when $n=0$ does not admit an immediate generalization to the case of $n\ge 1$, whose proof requires additional non-trivial (combinatorial) arguments. The proof is given in Section \ref{sec:proof of thm:Burge}.

\begin{ex} \label{ex:type d f0}
{\em 
Suppose that $m=4$ and $n=4$ and $({\bf i}, {\bf j})$ is as in Example \ref{ex:RSK for d23}.
We have $\tf_0 ({\bf m}({\bf i}, {\bf j}))= {\bf m}({\bf i}', {\bf j}')$, where
\begin{equation*}
\begin{split}
	\left(
		\begin{array}{c}
			{\bf i}' \\
			{\bf j}' 
		\end{array}
	\right)
	&=
	\left(
		\begin{array}{cccccc|c}
			\red{\bf 1} & 2 & 3 & 3 & 3 & 4 & 5 \\
			\red{\bf 2} & 8 & 4 & 4 & 4 & 6 & 7
		\end{array}
	\right),
\end{split}
\end{equation*}
and $P({\bf i}',{\bf j}')$ is given by
\begin{equation*}
\begin{ytableau}
    \none & \none & \none & \none & \tl{$\red{\bf 1}$}\\
    \none & \none & \none & \none & \tl{$\red{\bf 2}$}\\
    \none & \none & \none & \none & \tl{$3$}\\
    \none & \none & \none & \none & \tl{$4$}\\
    \tl{$2$} & \tl{$3$} & \tl{$3$} & \tl{$4$} & \tl{$6$}\\
    \tl{$4$} & \tl{$4$} & \tl{$5$} & \tl{$7$} & \tl{$8$}\\
\end{ytableau}
\end{equation*}
which is equal to $\tf_0(P({\bf i},{\bf j}))$.
}
\end{ex}

\subsection{Separation and the embedding of $\ms{B}(V(\la,\ell))$ into $\ms{B}(P(\la))$} \label{subsec: separation}

For $\la\in \cP_{m|n}$, set
\begin{equation*} \label{eq: def of Vla}
	{\bf V}^\g_\la = {\bf V}^\g \times SST_{m|n}(\la).
\end{equation*}
By \eqref{eq: def of Mla} and \eqref{eq: iso from M to V}, we have 
\begin{equation*}
 {\bf V}^\g_\la \cong {\bf M}^\g_\la.
\end{equation*}
So it remains to describe the embedding of ${\bf T}^\g(\la,\ell)$ into ${\bf V}^\g_\la$ in order to have \eqref{eq:main embedding}.

Let $(\la,\ell)\in \mc{P}(\g)$ be given. For ${\bf T} \in {\bf T}^\g(\lambda, \ell)$, we apply the combinatorial algorithms (called separation) in \cite[Section 5.1]{K18} for $\g=\mf{b}_{m|n}$, $\mf{c}_{m|n}$, and in \cite[Sections 4.2 and 4.3]{JK21} for $\g=\mf{d}_{m|n}$, which indeed do not depend on $n\ge 0$ (cf.~Remark \ref{rem:spinor convention}). 

Then we obtain an $\I$-semistandard tableau $\ov{\bf T}$ of skew shape satisfying the following properties;

\begin{itemize}
	\item[(S1)] $\ov{\bf T}$ is equivalent to ${\bf T}$ as elements of $\mf{l}$-crystals, 
	
	\item[(S2)] $\ov{\bf T}^{\texttt{tail}} \in SST_{m|n}(\la)$ and $\ov{\bf T}^{\texttt{body}} \in SST_{m|n}(\delta^{\pi})$ for some $\delta \in \ms{P}_\g$,
	
	\item[(S3)] ${\bf T}$ is equivalent to $\ov{\bf T}^{\texttt{body}} \otimes \ov{\bf T}^{\texttt{tail}}$ as elements of $\mf{l}$-crystals,
\end{itemize}
where $\ov{\bf T}^{\texttt{body}}$ and $\ov{\bf T}^{\texttt{tail}}$ are subtableaux  obtained by dividing $\ov{\bf T}$ with respect to a certain horizontal line $L$ (see Example \ref{ex: running ex 2}).

\begin{rem}{\rm
Let ${\bf T} = (T_l, \dots, T_1, T_0)\in {\bf T}^\g(\la, \ell)$ be given.
Recall that each $T_i$ is an $\I$-semistandard tableau of skew shape with at most two columns. 
Roughly speaking, $\ov{\bf T}$ is obtained from ${\bf T}$ by applying a sequence of jeu-de-taquin (for $\I$-semistandard tableaux) as far as possible so that no subtableau below the horizontal line $L$ is movable to the left.
} 
\end{rem}

\begin{ex}\label{ex: running ex 2}
{\rm
(1) Let $\g=\mf{c}_{2|5}$ and ${\bf T}=(T_3,T_2,T_1)$ as in Example \ref{ex: running ex 1} (1).
We first regard ${\bf T}$ as the following sequence of single columned tableaux
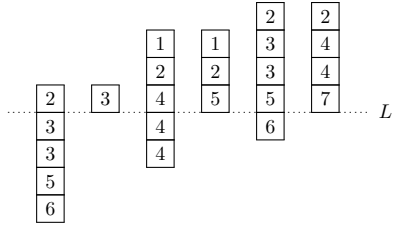
\begin{equation*}
    \begin{ytableau}
    \none & \none & \none & \none & \none & \none & \none & \none & \none & \tl{$2$} & \none & \tl{$2$} & \none \\
    \none & \none & \none & \none & \none & \tl{$1$} & \none & \tl{$1$} & \none & \tl{$3$} & \none & \tl{$4$} & \none \\
    \none & \none & \none & \none & \none & \tl{$2$} & \none & \tl{$2$} & \none & \tl{$3$} & \none & \tl{$4$} & \none \\
    \none[\mathrel{\raisebox{-0.7ex}{$\scalebox{0.45}{\dots\dots}$}}] & \tl{$2$} & 
	\none[\mathrel{\raisebox{-0.7ex}{$\scalebox{0.45}{\dots\dots}$}}] & \tl{$3$} & 
	\none[\mathrel{\raisebox{-0.7ex}{$\scalebox{0.45}{\dots\dots}$}}] & \tl{$4$} & 
	\none[\mathrel{\raisebox{-0.7ex}{$\scalebox{0.45}{\dots\dots}$}}] & \tl{$5$} & 
	\none[\mathrel{\raisebox{-0.7ex}{$\scalebox{0.45}{\dots\dots}$}}] & \tl{$5$} & 
	\none[\mathrel{\raisebox{-0.7ex}{$\scalebox{0.45}{\dots\dots}$}}] & \tl{$7$} & 
	\none[\ \ \ \mathrel{\raisebox{-0.7ex}{$\scalebox{0.45}{\dots\dots}$\ ${}_{\scalebox{0.75}{$L$}}$}}] \\
    \none & \tl{$3$} & \none & \none & \none & \tl{$4$} & \none & \none & \none & \tl{$6$} & \none & \none & \none \\
    \none & \tl{$3$} & \none & \none & \none & \tl{$4$} & \none & \none & \none & \none & \none & \none & \none \\
    \none & \tl{$5$} & \none & \none & \none & \none & \none & \none & \none & \none & \none & \none & \none \\
    \none & \tl{$6$} & \none & \none & \none & \none & \none & \none & \none & \none & \none & \none & \none[\quad \quad \quad \quad .] \\
    \end{ytableau}
\end{equation*}
\noindent We then apply jeu-de-taquin as far as possible so that no box below $L$ is movable to the left to get 
\vskip 2mm
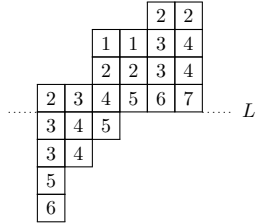
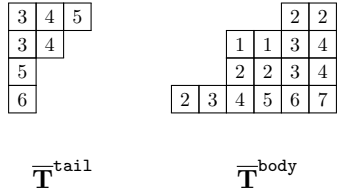
\begin{equation*}
    \begin{ytableau}
    \none & \none & \none & \none & \none & \tl{$2$} & \tl{$2$} & \none \\
    \none & \none & \none & \tl{$1$} & \tl{$1$} & \tl{$3$} & \tl{$4$} & \none \\
    \none & \none & \none & \tl{$2$} & \tl{$2$} & \tl{$3$} & \tl{$4$} & \none \\
    \none[\mathrel{\raisebox{-0.7ex}{$\scalebox{0.45}{\dots\dots}$}}] & \tl{$2$} & \tl{$3$} & \tl{$4$} & \tl{$5$} & \tl{$6$} & \tl{$7$} & \none[\ \ \ \mathrel{\raisebox{-0.7ex}{$\scalebox{0.45}{\dots\dots}$\ ${}_{\scalebox{0.75}{$L$}}$}}] \\
    \none & \tl{$3$} & \tl{$4$} & \tl{$5$} & \none & \none & \none & \none \\
	\none & \tl{$3$} & \tl{$4$} & \none & \none & \none & \none & \none\\
	\none & \tl{$5$} & \none & \none & \none & \none & \none & \none\\
	\none & \tl{$6$} & \none & \none & \none & \none & \none & \none\\
	\none & \none & \none & \none & \none & \none & \none & \none\\
	\none & \none & \none & \none & \none[\ov{\bf T}] & \none & \none & \none\\
    \end{ytableau}
\end{equation*}
 \noindent and hence \vskip 2mm
\begin{equation*}
    \qquad\quad
    \begin{ytableau}
		\tl{$3$} & \tl{$4$} & \tl{$5$} \\
		\tl{$3$} & \tl{$4$} & \none \\
		\tl{$5$} & \none & \none \\
		\tl{$6$} &  \none & \none \\
		\none & \none & \none \\
		\none & \none & \none \\
		\none & \none[\quad \ov{\bf T}^{\tt tail}] & \none \\
    \end{ytableau}
    \qquad\quad 
    \begin{ytableau}
		\none & \none & \none & \none & \tl{$2$} & \tl{$2$} \\
		\none & \none & \tl{$1$} & \tl{$1$} & \tl{$3$} & \tl{$4$} \\
		\none & \none & \tl{$2$} & \tl{$2$} & \tl{$3$} & \tl{$4$} \\
		\tl{$2$} & \tl{$3$} & \tl{$4$} & \tl{$5$} & \tl{$6$} & \tl{$7$} \\
		\none & \none & \none & \none & \none & \none \\
		\none & \none & \none & \none & \none & \none \\
		\none & \none & \none & \none[\ov{\bf T}^{\tt body}] & \none & \none[\quad \quad \quad.] \\
    \end{ytableau}
\end{equation*}
Thus we obtain two semistandard tableaux $\ov{\bf T}^{\rm body}$ and\, $\ov{\bf T}^{\rm tail}$ of shape $\delta^{\pi}$ and $\mu$, respectively, where $\delta = (6,4,4,2)$ and $\la = (3,2,1,1)$. \\
\noindent 
(2) Let $\g=\mf{d}_{4|4}$ and ${\bf T}=(T_4,T_3,T_2,T_1)$ as in Example \ref{ex: running ex 1} (2).
By applying the separation in \cite[Section 4.2]{JK19} to ${\bf T}$, we have
    \begin{equation*}
    \qquad\quad
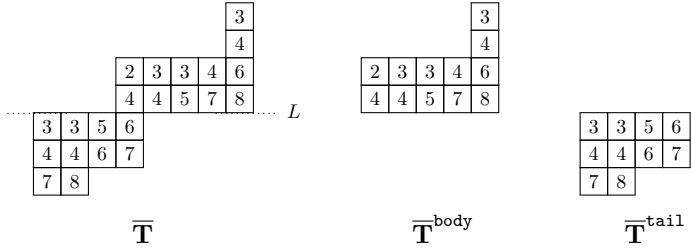
    \begin{ytableau}
    \none & \none & \none & \none & \none & \none & \none & \tl{$3$} & \none  \\
    \none & \none & \none & \none & \none & \none & \none & \tl{$4$} & \none \\
    \none & \none &\none & \tl{$2$} & \tl{$3$} & \tl{$3$} & \tl{$4$} & \tl{$6$}  & \none\\
    \none[\!\!\!\!\mathrel{\raisebox{-0.7ex}{$\scalebox{0.45}{\dots\dots\dots\dots}$}}] &\none & \none &\tl{$4$} & \tl{$4$} & \tl{$5$} & \tl{$7$} & \tl{$8$}   & \none[\!\!\!\!\mathrel{\raisebox{-0.7ex}{$\scalebox{0.45}{\dots\dots\dots\dots}$\ ${}_{\scalebox{0.75}{$L$}}$}}] \\
    \tl{$3$} & \tl{$3$} & \tl{$5$} & \tl{$6$} &\none & \none & \none & \none   & \none\\
    \tl{$4$} & \tl{$4$} & \tl{$6$} & \tl{$7$} &\none & \none & \none & \none  & \none \\
    \tl{$7$} & \tl{$8$} & \none & \none & \none & \none & \none & \none  & \none\\ 
    \none & \none & \none & \none & \none & \none & \none & \none  & \none\\ 
    \none & \none & \none & \none[\quad \ov{\bf T}] & \none & \none & \none & \none  & \none\\ 
    \end{ytableau}
    \qquad\quad 
    \begin{ytableau}
    \none & \none & \none & \none & \tl{$3$} & \none  \\
    \none & \none & \none & \none & \tl{$4$} & \none \\
    \tl{$2$} & \tl{$3$} & \tl{$3$} & \tl{$4$} & \tl{$6$}  & \none\\
    \tl{$4$} & \tl{$4$} & \tl{$5$} & \tl{$7$} & \tl{$8$}   & \none \\
    \none &\none & \none & \none & \none   & \none\\
    \none &\none & \none & \none & \none  & \none \\
    \none & \none & \none & \none & \none  & \none\\ 
    \none &\none & \none & \none & \none  & \none \\
    \none & \none & \none[\quad \ov{\bf T}^{\texttt{body}}] & \none & \none  & \none\\ 
    \end{ytableau}
    \qquad
    \begin{ytableau}
    \none & \none & \none & \none & \none & \none & \none & \none & \none  \\
    \none & \none & \none & \none & \none & \none & \none & \none & \none \\
    \none & \none &\none & \none & \none & \none & \none & \none  & \none\\
    \none &\none & \none & \none & \none & \none & \none & \none   & \none \\
    \tl{$3$} & \tl{$3$} & \tl{$5$} & \tl{$6$} &\none & \none & \none & \none   & \none\\
    \tl{$4$} & \tl{$4$} & \tl{$6$} & \tl{$7$} &\none & \none & \none & \none  & \none \\
    \tl{$7$} & \tl{$8$} & \none & \none & \none & \none & \none & \none  & \none\\ 
    \none &\none & \none & \none &\none & \none & \none & \none  & \none \\
    \none &\none & \none[\,\,\, \ov{\bf T}^{\texttt{tail}}] & \none & \none & \none & \none & \none  & \none[.]\\ 
    \end{ytableau}
    \end{equation*}
    Hence we obtain two semistandard tableaux $\ov{\bf T}^{\rm body}$ and\, $\ov{\bf T}^{\rm tail}$ of shape $\delta^{\pi}$ and $\mu$, respectively, where $\delta = (5,5,1,1)$ and $\la = (4,4,2)$.
}
\end{ex}

For $\Lambda\in P$, let $T_{\Lambda}=\{\,t_{\Lambda}\,\}$ be a $\g$-crystal with ${\rm wt}(t_\Lambda)=\Lambda$, and $\varepsilon_i(t_\Lambda)=\varphi_i(t_\Lambda)=-\infty$. Then we have the following.

\begin{thm} \label{thm:T to V}
For $(\la, \ell) \in \mc{P}(\g)$, the map
\begin{equation*}
\xymatrix@C=3em @R=0.5em{
{\bf T}^\g(\la, \ell) \otimes T_{-\ell\Lambda_0} \ar@{->}[r] & {\bf V}^\g_\la \\
{\bf T} \ot t_{-\ell\Lambda_0} \ar@{|->}[r] 
&  \left(\ov{\bf T}^{\texttt{\em body}},\ov{\bf T}^{\texttt{\em tail}} \right)
}
\end{equation*}
is an injective morphism of $\g$-crystals.
\end{thm}
\pf It follows from the same arguments in \cite[Theorem 5.7]{K18} for $\g = \mf{b}_{m|n},  \mf{c}_{m|n}$, and \cite[Theorem 4.12]{JK21} for $\g = \mf{d}_{m|n}$ simply by replacing $\I_{\ov 0}$-semistandard tableaux in \cite{JK21,K18} with $\I$-semistandard ones.
\qed
\smallskip

Summarizing, we have the following realization of \eqref{eq:main embedding}:
\begin{equation*}
\xymatrix@C=2.5em @R=0.5em{
\ms{B}(V(\la,\ell))\otimes T_{-\ell\Lambda_0}\ \cong\ {\bf T}^\g(\la, \ell) \otimes T_{-\ell\Lambda_0} \ar@{->}[r] & {\bf V}^\g_\la \ar@{->}^{\hskip -10mm\cong}[r] & {\bf M}^\g_\la \ \cong\ \ms{B}(P(\la)) \\
\quad\quad\quad\quad\quad\quad\quad\quad {\bf T} \ot t_{-\ell\Lambda_0} \ar@{|->}[r] 
&  \left(\ov{\bf T}^{\texttt{body}},\ov{\bf T}^{\texttt{tail}} \right) \ar@{|->}[r] & \left({\bf m},\ov{\bf T}^{\texttt{tail}} \right)\quad\quad\quad\quad
},
\end{equation*}
where $\kappa^\g({\bf m})=\ov{\bf T}^{\texttt{body}}$.

\begin{ex}\label{ex: running ex 4}
{\rm
\noindent (1) Let $\g={\mf c}_{2|5}$ and ${\bf T}$ is as in Example \ref{ex: running ex 1}(1).
Then corresponding ${\mf m} \in \bM^{\g}$ is given in Example \ref{ex: running ex 3},
and $\ov{\bf T}^{\texttt{tail}}$ is given in Example \ref{ex: running ex 2}(1).
\\
\noindent (2) Let $\g={\mf d}_{4|4}$ and ${\bf T}$ is as in Example \ref{ex: running ex 1}(2). 
Then corresponding ${\mf m} \in \bM^{\g}$ is associated with $(\bi, \bj)$ in Example \ref{ex:RSK for d23},
and $\ov{\bf T}^{\texttt{tail}}$ is given in Example \ref{ex: running ex 2}(2).
} 
\end{ex}

\begin{rem}\label{rem: proof of ch of P(la)}
{\rm Let us give a proof of Lemma \ref{lem: character of P(la)}. Let $\xi=\sum_{i\in I}c_i\alpha_i \in Q^-$ be given. 
Since the action of $\tf_0$ is given by adding a vertical domino, we see that the number of boxes in ${\bf T}$ above the horizontal line $L$ is bounded by $-2c_0$ for any $(\lambda,\ell)\in \cP(\g)$ and ${\bf T}\in {\bf T}^\g(\la, \ell)_{\La(\la,\ell)+\xi}$. 

Let $(\lambda,\ell)\in \cP(\g)$ be given.
From the above observation, we may choose a sufficiently large $t$ such that for any ${\bf T}\in {\bf T}^\g(\la, \ell)_{\La(\la,\ell+t)+\xi}$, all boxes in ${\bf T}$ above $L$ belong to the columns which are strictly to the right of those the boxes below $L$ belong to. This implies that $\ov{\bf T}^{\texttt{body}}$ (resp. $\ov{\bf T}^{\texttt{tail}}$) is simply obtained by the sliding the boxes above $L$ (resp. below $L$) to the right (resp. left) separately. This yields a bijection
\begin{equation*}
\xymatrix@C=2.5em @R=0.5em{
{\bf T}^\g(\la, \ell)_{\La(\la,\ell+t)+\xi} \otimes T_{-(\ell+t)\Lambda_0} \ar@{->}[r] & ({\bf V}^\g_\la)_{\La(\la,0)+\xi},
}
\end{equation*}
and hence a bijection from $\ms{B}(V(\la,\ell))_{\La(\la,\ell+t)+\xi}\otimes T_{-(\ell+t)\Lambda_0}$ to $\ms{B}(P(\la))_{\La(\la,0)+\xi}$.
This proves the lemma.
}
\end{rem}


\section{Proof of Theorem \ref{thm: Burge}}\label{sec:proof of thm:Burge}
Throughout this section, we assume that $\g = \mf{d}_{m|n}$. 
For simplicity, we omit $\g$ as super- or subscripts in certain notations if there is no confusion. 
For example, we simply write $\bM = \bM^{\mf{d}_{m|n}}$ and $\bV = \bV^{\mf{d}_{m|n}}$.
But we denote $\kappa^{\mf{d}_{m|n}}$ by $\burge$ following \cite{JK19} to avoid notational confusion. 
This section is organized as follows.  
In Section \ref{subsec:super Burge}, we check that $\burge$ is a well-defined weight-preserving bijection between $\bM$ and $\bV$. 
In Section \ref{subsec: Burge and RSK}, we investigate a relation between Burge and RSK correspondences, which plays a crucial role in the proof of Theorem \ref{thm: Burge},
while we present the proofs of some delicate combinatorial arguments in Section \ref{subsec: combinatorics of Q bar}.
Then in Section \ref{sec: proof of Burge}, we give the proof of Theorem \ref{thm: Burge}. 
\smallskip

\noindent
{\bf Convention}.
Thorughout Section \ref{sec:proof of thm:Burge}, we adopt the following conventions.
\begin{enumerate}[(1)]
	\item For $T \in SST_{m|n}(\la^\pi)$, we denote by $T_i$ (resp.~$T^i$) the $i$-th column (resp.~row) of $T$ from the right (resp.~bottom) so that the length of $T_i$ (resp.~$T^i$) is greater than or equal to that of $T_{i+1}$ (resp.~$T^{i+1}$). 
	\smallskip

	\item For two tableaux $S \in SST_{m|n}(\la^\pi)$ and $T \in SST_{m|n}(\mu^\pi)$, we say that $T$ is a \emph{subtableau} of $S$, denoted by $T \subset S$, if $\mu \subset \la$ and each letter of $T$ is exactly same with that of $S$ on $\mu^\pi$. When $T$ is a subtableau of $S$, we define $S / T$ to be the skew tableau of shape $(\la / \mu)^\pi$ obtained by removing $T$ inside $S$. 
	We often write $S = (S / T) \sqcup T$.
\end{enumerate}

\subsection{Well-definedness of Burge correspondence} \label{subsec:super Burge}

For $T \in SST_{m|n}(\la^\pi)$ and $a \in \I$, let $i_k$ be the letter bumped out from $T_k$ along the insertion of $a$ into $T$ to have $\left( T \leftarrow a \right)$, where we set $i_0 = a$. 
The sequence of boxes including $i_k$'s in $T$ is called the \emph{bumping route} of $a$ in $T$, and is denoted by $R_T(a)$.
We denote by $B_T(a)$ the box corresponding to $\mathrm{sh}\left( T \leftarrow a \right)^\pi / \mathrm{sh}\left( T \right)^\pi$. 
We say that a route $R$ is \emph{strictly below} (resp.~\emph{weakly below}) a route $R'$ if in each column which contains a box of $R'$, $R$ has a box which is below (resp.~below or equal to) the box in $R'$. 
The notions of strictly above (resp. weakly above) are defined analogously.
The following lemma is a super analogue of \cite[Appendix A.2]{F}, and its proof is straightforward.

\begin{lem} [Column Bumping Lemma] \label{lem: bumping routes}
Let $T \in SST_{m|n}(\la^\pi)$ be given. 
For $a, b \in \I$, put $U = \left( T \leftarrow a \right)$ and $V = \left( U \leftarrow b \right)$.
\begin{enumerate}[\em (1)]
    \item If $a, b \in \I_{\ov{0}}$ with $a \le b$, then $R_U(b)$ is weakly below $R_T(a)$, and $B_U(b)$ is southwest of $B_T(a)$. \label{it:bump2}
    \item If $a, b \in \I_{\ov{0}}$ with $a > b$, then $R_U(b)$ is strictly above $R_T(a)$, and $B_U(b)$ is northeast of $B_T(a)$. \label{it:bump1}
    \item If $a \in \I_{\ov{0}}$ and $b \in \I_{\ov{1}}$, then $R_U(b)$ is weakly below $R_T(a)$ and $B_U(b)$ is southwest of $B_T(a)$.\!\label{it:bump3}
    \item If $a \in \I_{\ov{1}}$ and $b \in \I_{\ov{0}}$, then $R_U(b)$ is strictly above $R_T(a)$ and $B_U(b)$ is northeast of $B_T(a)$.\!\label{it:bump4}
    \item If $a, b \in \I_{\ov{1}}$ with $a < b$, then $R_U(b)$ is weakly below $R_T(a)$ and $B_U(b)$ is southwest of $B_T(a)$. \label{it:bump6}
    \item If $a, b \in \I_{\ov{1}}$ with $a \ge b$, then $R_U(b)$ is strictly above $R_T(a)$ and $B_U(b)$ is northeast of $B_T(a)$. \label{it:bump5}
\end{enumerate}
\end{lem}

\begin{prop} [cf.~\cite{B, JK19, Su}] \label{prop:well-definedness}
	Let $\bc \in \bM$ be given. 
	\begin{enumerate}[\em (1)]
		\item $\burge \left( \bc \right) \in \bV$.
		\item $\burge$ is a weight-preserving bijection between $\bM$ and $\bV$.
	\end{enumerate}  
\end{prop}
\pf
(1) can be proved directly by induction on the sum of $c_{(i,j)}$ in $\bc$, and applying Lemma \ref{lem: bumping routes} at each step.
To prove (2), we describe the reverse map of $\burge$. 
First, associate $\bc$ with $({\bf i}, {\bf j}) \in \biws$ for some ${\bf i}=i_1 \dots i_r$, ${\bf j}=j_1 \dots j_r$ and $r \in \Z_{\ge 0}$, 
and consider the positions of $i_k$'s in $\burge(\bc)=P({\bf i},{\bf j})$.
Suppose $i_k = i_{k+1} = i \in \I_{\ov 0}$. 
Then, by \eqref{eq: df of biws}, we have $j_k\geq j_{k+1}$ ($j_k > j_{k+1}$ if $j_k \in \I_{\ov 1}$). 
Applying Lemma \ref{lem: bumping routes} \eqref{it:bump2}, \eqref{it:bump3} and \eqref{it:bump6}, it follows that $i_k$ is located southwest of $i_{k+1}$ in $\burge({\bc})$.
Hence, the leftmost $i$ of $\burge({\bc})$ is the most recent one.
Similarly, 
if we assume $i\in \I_{\ov 1}$, then the topmost $i$ is the most recent one by Lemma \ref{lem: bumping routes} \eqref{it:bump6}.

Now let $T \in \bV$ be given and let $i$ denote the smallest entry in $T$. 
If $i \in \I_{\ov 0}$ (resp. $i \in \I_{\ov 1}$), remove the leftmost (resp. topmost) $i$ and then apply the reverse column insertion with respect to the entry immediately below it.
In this way, the removed entry is identified with $i_1$ and the element bumped out is $j_1$.
Repeating this procedure until the tableau becomes empty yields $({\bf i}, {\bf j}) \in \biws$ with ${\bf i}=i_1 \dots i_r$ and ${\bf j}=j_1 \dots j_r$.
Therefore, $\burge$ is a bijection map.
\qed

\subsection{Relation between Burge and RSK correspondences} \label{subsec: Burge and RSK}
We investigate the relation between Burge and RSK correspondences generalizing the work of the case when $n = 0$ \cite{JK19}.
\smallskip

Let us follow the convention in \cite[Section 3.2]{JKU24}.
For $\bc \in \bM$ and $i \in \I$, we use the following notations:
\begin{enumerate}[\qquad$\bullet$]
	\item $\bti$ : the triangluar-shaped arrangement of $(k,l) \in \Phi(\mf{u})$ for $i \le k \le l \le m+n$,
	\smallskip
	
	\item $\btic = \bt_1 \setminus \bti$ : the complement of $\bti$ in $\bt_1$, 
	\smallskip
		
	\item $\ti$ : the subset of $\bti$ given by
		\begin{equation*}
			\begin{cases}
				\{ (i, i+1) \} & \text{if $i < m$,} \\
				\{ (i, i+1), (i+1, i+1) \} & \text{if $i = m$,} \\
				\{ (i, i), (i, i+1), (i+1, i+1) \} & \text{if $i > m$,}
			\end{cases}
		\end{equation*}
	
	\item $\lozenge_i=\bti \setminus ( \ti \cup \btiII)$ : the complement of $\ti \cup \btiII$ in $\bti$, 
	\smallskip
	
	\item $\bc_X$ : the element of $\bM$ whose components on $X$ coincide with those of $\bc$, and $0$ elsewhere,
	
	\item $\bM_X$ : the subset of $\bM$ consisting of $\bc \in \bM$ such that every component of $\bc_{X^c}$ is $0$,
\end{enumerate}
where $X$ is a subset of $\bt_1$ and $X^c$ is the complement of $X$ in $\bt_1$.
\smallskip

\begin{ex}
	{\em
	Let $\g=\mf{d}_{2|4}$. 
	Then $\bt_2$ is decomposed as $\bt_2= {\bt_4} \sqcup {\triangle_2} \sqcup{{\lozenge}_2}$ as illusted in the following diagram:
	
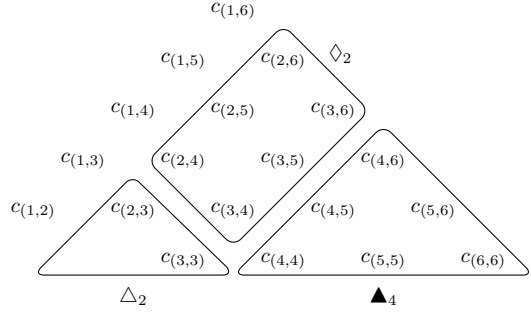
\begin{equation*}
		\begin{tikzpicture}[baseline=(current  bounding  box.center), every node/.style={scale=0.8}, scale=0.95]
			\node (glw_01m) at (0.7,0.7) {$c_{(1,2)}$};
			\node (glw_02m) at (1.4,1.4) {$c_{(1,3)}$};
			\node (glw_03m) at (2.1,2.1) {$c_{(1,4)}$};
			\node (glw_04m) at (2.8,2.8) {$c_{(1,5)}$};
			\node (glw_04m) at (3.5,3.5) {$c_{(1,6)}$};
			\node (glw_12m) at (2.1,0.7) {$c_{(2,3)}$};
			\node (glw_13m) at (2.8,1.4) {$c_{(2,4)}$};
			\node (glw_14m) at (3.5,2.1) {$c_{(2,5)}$};
			\node (glw_14m) at (4.2,2.8) {$c_{(2,6)}$};
			\node (glw_02f) at (2.8,0) {$c_{(3,3)}$};
			\node (glw_23m) at (3.5,0.7) {$c_{(3,4)}$};		
			\node (glw_24m) at (4.2,1.4) {$c_{(3,5)}$};
			\node (glw_24m) at (4.9,2.1) {$c_{(3,6)}$};
			\node (glw_03f) at (4.2,0) {$c_{(4,4)}$};
			\node (glw_34m) at (4.9,0.7) {$c_{(4,5)}$};
			\node (glw_34m) at (5.6,1.4) {$c_{(4,6)}$};
			\node (glw_04m) at (5.6,0) {$c_{(5,5)}$};
			\node (glw_04m) at (6.3,0.7) {$c_{(5,6)}$};
			\node (glw_04m) at (7.0,0) {$c_{(6,6)}$};
			\draw[rounded corners] (4.2-0.7,0-0.2) -- (7.0+0.7,0-0.2) -- (5.6,1.4 + 0.5) -- cycle;
			\draw (5.6,-0.5) node{$\bt_4$};
	        \begin{scope}[shift={(-2.8,0)}]
				\draw[rounded corners] (4.2-0.7,0-0.2) -- (5.6+0.7,0-0.2) -- (4.9,0.7 + 0.5) -- cycle;
			\end{scope}
			\draw (2.1,-0.5) node{$\triangle_2$};
			\draw[rounded corners] (2.8-0.5,1.4) -- (4.2,2.8+0.5) -- (4.9+0.5,2.1) -- (3.5,0.7-0.5) -- cycle;
			\draw (5.0,2.9) node{${\lozenge}_2$};
		\end{tikzpicture}
	\end{equation*}
	}
\end{ex}
When we consider $\te_i$ or $\tf_i$ on $\bc$ for $i \in I \setminus \{ 0 \}$, we often decompose $\bc $ into $\bc_{\btic} \ot \bc_{\bti}$.
Thus we first focus on $\bc_{\bti} \in \bM_\bti$ in this subsection.

Let $T \in SST_{m|n}(\mu^\pi)$ and $\bc \in \bM$ with $\bc = \bc(\ba ,\bb)$, where $\mu \in \cP_{m|n}$ and $(\ba,\bb) = (a_1 \dots a_s\,,\, b_1 \dots b_s) \in \biws$. Then we define 
\begin{enumerate}[\qquad$\bullet$]
	\item $\tP(T, \bc) = \left( \left( T \leftarrow b_s \right) \leftarrow \cdots \leftarrow b_1 \right)$,
	
	\item $\tQ(T, \bc) \in SST\left( \left( \la / \mu \right)^\pi \right)$, where the subtableau of $\tQ(T, \bc)$ corresponding to ${\rm sh}(\tP_k) / {\rm sh}(\tP_{k+1})$ is filled with $a_k$ for $1 \le k \le s$,
\end{enumerate}
where $\la = {\rm sh}\left( \tP(T, \bc) \right)^\pi$ and $\tP_k = \left( \left( T \leftarrow b_s \right) \leftarrow \cdots \leftarrow b_k \right)$ with $\tP_{s+1} = T$. 

Fix $i \in I \setminus \{ 0, m \}$.
Let $\bc \in \bM_{\bti}$ be given.
We define 
\begin{equation*}
\tP_i(\bc)=\tP(T,\bc_{\lozenge_i}) \quad \text{and} \quad \tQ_i(\bc)=\tQ(T,\bc_{\lozenge_i}),
\end{equation*}
where $T = \burge(\bc_\btiII)$. 
For a tableau $U$, let $U_{\geq k}$ be the subtableau of $S$ consisting of entries greater than or equal to $k$.

\begin{lem} \label{lem: P-tab}
We have $\tP_i(\bc)=\burge(\bc)_{\geq i+2}$.
\end{lem}
\pf
We proceed by induction on $s$, the length of biword corresponding to $\bc_{\lozenge_i}$.
When $s = 0$, our assertion clearly holds.
Let ${\bf 1}_{(a_1, b_1)} \in \bM$ be defined by $c_{(a_1, b_1)} = 1$ and $c_\gamma = 0$ for all $\gamma \in \bt_1 \setminus \{ (a_1, b_1) \}$.
By induction hypothesis, we have $\tP_i(\bc-{\bf 1}_{(a_1, b_1)}) = \burge\left( \bc - {\bf 1}_{(a_1, b_1)} \right)_{\geq i+2}$, 
where $\bc - {\bf 1}_{(a_1, b_1)}$ is understood as the component-wise subtraction. 
Then we can show that the bumping route of $(\tP_i\left( \bc - {\bf 1}_{(a_1, b_1)} \right) \leftarrow b_1)$
is contained in the bumping route of $(\burge\left( \bc - {\bf 1}_{(a_1, b_1)} \right) \leftarrow b_1)$, which yields the desired statement.
\qed

\smallskip

The following lemma is one of key ingredients in the proof of Theorem \ref{thm: Burge}.

\begin{lem} \label{lem: tfi on P-tab and Q-tab}
Let us identify $\bc$ with $\bc_\ti \ot \bc_{\lozenge_i} \ot \bc_\btiII$. 
If
\begin{equation} \label{eq:technical assumption}
	\tf_i \left( \bc_\ti \ot \bc_{\lozenge_i} \ot \bc_\btiII \right) = \bc_\ti \ot \left( \tf_ i \bc_{\lozenge_i} \right) \ot \bc_\btiII,
\end{equation}
and $\tf_i \bc_{\lozenge_i} \neq {\bf 0}$, then we have
\begin{enumerate}[{\em (1)}]
	\item $\tP_i( \tf_i \bc ) = \tP_i(  \bc )$,
	
	\item $\tQ_i( \tf_i \bc ) = \tf_i \tQ_i( \bc )$.
\end{enumerate}
\end{lem}
\pf
Recall that $i \in I \setminus \{ 0, m \}$.
Note that for $1 \leq i<m$, \eqref{eq:technical assumption} always holds, as both $\bc_\ti$ and $\bc_\btiII$ have a trivial crystal structure with respect to $\te_i$ and $\tf_i$ \cite[Section 4.2]{JKU24}.

Let $T=\burge(\bc_{\btiII})$ and $\nu = {\rm sh}\left( T \right)$.
Each element of $\bM_{\lozenge_i}$ can be regarded as an element of $\mathbf{M}$ in \cite[Section 3.2]{K07},
and the action of $\tf_i$ on $\bM_{\lozenge_i}$ coincides with that on $\mathbf{M}$ by \cite[Section 4.2]{JKU24}.
Under this identification, statements (1) and (2) follow by applying the same argument as in \cite[Proposition 4.6]{K09}, 
where \cite[Proposition 3.3]{K09} is replaced by \cite[Theorem 3.11]{K07} (see also the paragraph of \cite[Theorem 3.11]{K07} above).
This is because the map $(\bc_{\lozenge_i},\emptyset,T) \mapsto \left( \tP_i(\bc), \tQ_i(\bc) \right)$ for $\bc_{\lozenge_i} \in \bM_{\lozenge_i}$ is a crystal morphism,
which may be regarded as a super-analogue of the inverse of $\Psi_{\emptyset, \nu; s}$ in \cite[Proposition 4.6]{K09} for $s \ge \ell(\nu)$.
\qed
\smallskip

We now present another construction of $\burge(\bc)$ in terms of $\tP_i(\bc)$, $\tQ_i(\bc)$ and $\bc_\ti$ in the following algorithms.
Recall that $\burge(\bc)$ and $\tP_i(\bc)$ only differ by their subtableaux with entries $i$ and $i+1$.
We show that $\burge(\bc)$ can be given by gluing $\tP_i(\bc)$ and a tableau which is determined by $\tQ_i(\bc)$ and $\bc_\ti$.

Algorithm \ref{alg: Q bar i<m} for the case of $1\leq i<m$ is essentially identical to that in \cite[Section 5.4]{JK19},
whereas Algorithm \ref{alg: Q bar i>m} provides its super-analogue.
Let $\tP(\bc)=\tP_i(\bc)$ and $\tQ(\bc)=\tQ_i(\bc)$ for simplicity.
\begin{alg} \label{alg: Q bar i<m}
{\em
Suppose that $1 \leq i<m$. We construct a tableau $\tT(\bc)$ as follows:
\begin{itemize}
	\item[(1)] Let $\la=\mathrm{sh}(\tP(\bc))^\pi$ and $L=\lfloor \ell(\la)/2 \rfloor$.
		For $1 \le l \le L$, let $\tQ(\bc)^l$ be the subtableau consisting of the $(2l-1)$-th and $2l$-th rows of $\tQ(\bc)$ from the bottom row.
		Apply jeu de taquin to $\tQ(\bc)^l$, to obtain a normal-shaped tableau,
		which we denote by $\ov{\tQ}(\bc)^l$.
	\item[(2)] 	Let $\tP(\bc)^l$ be the subtableau of $\tP(\bc)$ corresponding to the same positions as $\tQ(\bc)^l$. 
		We also denote by $\tP(\bc)^{\ge l}$ the subtableau of $\tP(\bc)$ lying strictly above the $(2l-2)$-th row from the bottom, with the convention $\tP(\bc)^{\ge 1} = \tP(\bc)$.
		Define $\tT_L(\bc)$ to be the two-rowed tableau obtained by gluing $\ov{\tQ}(\bc)^{L}$ to $\tP(\bc)^{L}$ and $c_{(i,i+1)}$ copies of \scalebox{0.8}{$\vd{\phantom{+}i\phantom{\,+}}{i+1}$} to the left of the resulting tableau.
	\item[(3)] Next, we construct a tableau from $\tP(\bc)^{\ge L-1}$ by replacing the subtableau $\tP(\bc)^{\ge L}$ with $\tT_L(\bc)$,
		and then gluing $\ov{\tQ}(\bc)^{L-1}$ to the corresponding subtableau $\tP(\bc)^{L-1}$. 
		If the number of columns in $\tT_L(\bc)$ exceeds $\mu_{2L-3}-\mu_{2L-1}$, where $\mu=\mathrm{sh}(T)^\pi$,
		then move dominos \scalebox{0.8}{$\vd{\phantom{+}i\phantom{\,+}}{i+1}$} down to the next two rows as many as it exceeds, 
		and denote the resulting tableau by $\tT_{L-1}(\bc)$.
	\item[(4)] Repeat the step (3) to obtain $\tT_{L-2}(\bc), \dots, \tT_1(\bc)$, and finally set $\tT(\bc)=\tT_1(\bc)$.
\end{itemize}
}
\end{alg}

\begin{ex} \label{ex: gluing i<m}
{\em 
Suppose $\g = \mf{d}_{2|4}$ and set $i=1$.
Let us consider ${\bf c} ={\bf c}({\bf i}, {\bf j}) \in \bM_{\bt_1}$ where
\begin{equation*}
\begin{split}
	\left(
		\begin{array}{c}
			{\bf i} \\
			{\bf j} 
		\end{array}
	\right)
	&=
	\left(
		\begin{array}{ccccccccccccccc}
			1 & 1 & 1 & 1 & 1 & 2 & 2 & 2 & 3 & 3 & 3 & 3 & 4 & 4 & 5 \\
			6 & 5 & 4 & 2 & 2 & 6 & 4 & 3 & 4 & 5 & 6 & 6 & 4 & 6 & 5
		\end{array}
	\right).
\end{split}
\end{equation*}
Let us compute $\tT(\bc)$ following Algorithm \ref{alg: Q bar i<m}.
By \eqref{df: type D insertion}, $T=\burge(\bc_{\bt_{i+2}})$ is given by
\begin{equation}\label{ex: D24-1}
\vcenter{%
\xymatrix@R=-0.82em @C=-0.82em{
	&  & \boxed{3} \\
	&  & \boxed{4} \\
	\boxed{3} & \boxed{4} & \boxed{5} \\
	\boxed{3} & \boxed{4} & \boxed{6} \\
	\boxed{3} & \boxed{5} & \boxed{6} \\
	\boxed{4} & \boxed{5} & \boxed{6}
	}
}%
\in SST_{2|4}(\mu^\pi),
\end{equation}
where $\mu=(3,3,3,3,1,1)$, and
\begin{equation*}
	\tP(\bc_{\bt_3},\bc_{\lozenge_1}) = 
	\vcenter{\xymatrix@R=-0.82em @C=-0.82em{
	&  &  & \boxed{\blue{4}} \\
	&  & \boxed{\blue{3}} & \boxed{4} \\
	&  \boxed{\blue{3}} & \boxed{\blue{4}} & \boxed{6} \\
	&  \boxed{3} & \boxed{5} & \boxed{6} \\
	&  \boxed{4} & \boxed{5} & \boxed{6} \\
	\boxed{\blue{3}} &  \boxed{4} & \boxed{5} & \boxed{6} \\
	\boxed{\blue{3}} &  \boxed{4} & \boxed{5} & \boxed{6}
	}}, \qquad
	\tQ(\bc_{\bt_3},\bc_{\lozenge_1}) = 
	\vcenter{\xymatrix@R=-0.82em @C=-0.82em{
	&  &  & \boxed{\blue{2}} \\
	&  & \boxed{\blue{1}} & \\
	&  \boxed{\blue{1}} & \boxed{\blue{2}} & \\
	&  &  & \phantom{\boxed{9}}\\
	&  &  & \phantom{\boxed{9}}\\
	\boxed{\blue{1}} &  &  & \\
	\boxed{\blue{2}} &  &  &
	}}.
\end{equation*}
Then $\tQ(\bc)^l$'s and $\ov{\tQ}(\bc)^l$'s are given by
\begin{equation*}
	\tQ(\bc)^1 = 
	\vcenter{\xymatrix@R=-0.82em @C=-0.82em{
	\boxed{\blue{1}} \\
	\boxed{\blue{2}}
	}}, \qquad
	\tQ(\bc)^2 = 
	\emptyset, \qquad
	\tQ(\bc)^3 = 
	\vcenter{\xymatrix@R=-0.82em @C=-0.82em{
	& \boxed{\blue{1}} \\
	\boxed{\blue{1}} & \boxed{\blue{2}}
	}}, \qquad
	\tQ(\bc)^4 = 
	\vcenter{\xymatrix@R=-0.82em @C=-0.82em{
	\boxed{\blue{2}}
	}},
\end{equation*}
\begin{equation*}
	\ov{\tQ}(\bc)^1 = 
	\vcenter{\xymatrix@R=-0.82em @C=-0.82em{
	\boxed{\blue{1}}\\
	\boxed{\blue{2}}
	}}, \qquad
	\ov{\tQ}(\bc)^2 = 
	\emptyset, \qquad
	\ov{\tQ}(\bc)^3 = 
	\vcenter{\xymatrix@R=-0.82em @C=-0.82em{
	\boxed{\blue{1}} & \boxed{\blue{1}} \\
	\boxed{\blue{2}} &
	}}, \qquad
	\ov{\tQ}(\bc)^4 = 
	\vcenter{\xymatrix@R=-0.82em @C=-0.82em{
	\boxed{\blue{2}}
	}}.
\end{equation*}
Thus $\tT_4(\bc)$, $\tT_3(\bc)$, $\tT_2(\bc)$ and $\tT_1(\bc)$ are as follows:
\begin{equation*}
	\tT_4(\bc) = 
	\vcenter{\xymatrix@R=-0.82em @C=-0.82em{
	\boxed{\red{1}} & \boxed{\red{1}} & \boxed{\blue{2}} \\
	\boxed{\red{2}} & \boxed{\red{2}} & \boxed{\blue{4}}
	}}, \qquad
	\tT_3(\bc) = 
	\vcenter{\xymatrix@R=-0.82em @C=-0.82em{
	&  &  &  &  & \boxed{\blue{2}} \\
	&  &  &  &  & \boxed{\blue{4}} \\
	\boxed{\red{1}} & \boxed{\red{1}} &  \boxed{\blue{1}} &  \boxed{\blue{1}} & \boxed{\blue{3}} & \boxed{4} \\
	\boxed{\red{2}} & \boxed{\red{2}} &  \boxed{\blue{2}} &  \boxed{\blue{3}} & \boxed{\blue{4}} & \boxed{6}
	}}, \qquad
	\tT_2(\bc) = 
	\vcenter{\xymatrix@R=-0.82em @C=-0.82em{
	&  &  &  &  & \boxed{\blue{2}} \\
	&  &  &  &  & \boxed{\blue{4}} \\
	&  &  &  \boxed{\blue{1}} & \boxed{\blue{3}} & \boxed{4} \\
	&  &  &  \boxed{\blue{3}} & \boxed{\blue{4}} & \boxed{6} \\
	\boxed{\red{1}} &  \boxed{\red{1}} &  \boxed{\red{1}}  &  \boxed{3} & \boxed{5} & \boxed{6} \\
	\boxed{\red{2}} &  \boxed{\red{2}} &  \boxed{\red{2}}  &  \boxed{4} & \boxed{5} & \boxed{6}
	}},
\end{equation*}
\begin{equation*}
\tT(\bc)=\tT_1(\bc)=
\vcenter{%
\xymatrix@R=-0.82em @C=-0.82em{
	&  &  &  &  &  &  & \boxed{\blue{2}} \\
	&  &  &  &  &  &  & \boxed{\blue{4}} \\
	&  &  &  &  &  \boxed{\blue{1}} & \boxed{\blue{3}} & \boxed{4} \\
	&  &  &  &  &  \boxed{\blue{3}} & \boxed{\blue{4}} & \boxed{6} \\
	&  &  &  &  &  \boxed{3} & \boxed{5} & \boxed{6} \\
	&  &  &  &  &  \boxed{4} & \boxed{5} & \boxed{6} \\
	\boxed{\red{1}} & \boxed{\red{1}} & \boxed{\red{1}} & \boxed{\blue{1}} & \boxed{\blue{3}} &  \boxed{4} & \boxed{5} & \boxed{6} \\
	\boxed{\red{2}} & \boxed{\red{2}} & \boxed{\red{2}} & \boxed{\blue{2}} & \boxed{\blue{3}} &  \boxed{4} & \boxed{5} & \boxed{6}
	}
}%
,
\end{equation*}
where the red \scalebox{0.8}{$\vd{1}{2}$}'s in $\tT_{4}(\bc)$ are the $c_{(1,2)}=2$ copy of dominos described in Algorithm \ref{alg: Q bar i<m} (2),
while the red \scalebox{0.8}{$\vd{1}{2}$}'s in $\tT_{l}(\bc)$ for $1\leq l<4$ are those moved down from $\tT_{l+1}(\bc)$ as described in Algorithm \ref{alg: Q bar i<m} (3).
One can check that $\burge(\bc)$ coincides with $\tT(\bc)$ by \eqref{df: type D insertion}.
}
\end{ex}

For ${\sf z} \in \Z_+$, let $\qo{\sf z} \in \Z_{\ge 0}$ and $\re{\sf z} \in \{ 0, 1 \} $ be such that ${\sf z} = 2\qo{\sf z} + \re{\sf z}$.

\begin{alg} \label{alg: Q bar i>m}
{\em
Suppose that $i>m$. We construct a tableau $\tT(\bc)$ as follows:
\smallskip

\begin{itemize}
	\item[(1)] Let $\bsfT(\bc)=({\sf x}_k; {\sf y}_k)_{k \ge -1}$ be a two-rowed array as follows:
	\begin{equation} \label{eq: df of two array form} 
		\bsfT(\bc) = \quad
	\begin{tabular}{c|ccccccc}
			& $\cdots$ & $k$ & $\cdots$ & $2$ & $1$ & $0$ & $ -1 $  \\ \hline
		${\sf x}_k$	& $\cdots$ & ${\sf a}_k$ & $\cdots$ & ${\sf a}_2$ & ${\sf a}_1$ & $c_{(i , i+1)} + 2 c_{(i,i)}$ & $ 0 $ \\ 
		${\sf y}_k$    & $\cdots$ & ${\sf b}_k$ & $\cdots$ & ${\sf b}_2$ & ${\sf b}_1$ & $2 c_{(i+1, i+1)} + \qo{c_{(i,i+1)}}$ & $2 \re{c_{(i,i+1)}}$
		\end{tabular},
	\end{equation}
	where ${\sf a}_k$ (resp. ${\sf b}_k$) is the number of $i$'s (resp.~$i+1$'s) in $k$-th column of $\tQ(\bc)$ from the right.
	\vskip 2mm
	\item[(2)] Define a two-rowed array $\mc{F}(\bsfT(\bc))= (\widetilde{\sf{X}}_k ; \widetilde{\sf{Y}}_k)_{k \ge -1}$ by the following steps: for each $k\ge -1$
	\begin{enumerate}
		\item[(I)] \label{it:shift II} Replace ${\sf y}_k$ in $\bsfT(\bc)$ by $\re{{\sf y}_k} + 2\qo{{\sf y}_{k-1}}$.
		\item[(II)] \label{it:shift iI} Set  ${\sf z}_k = \min \left\{\, {\sf x}_k,\, \re{{\sf y}_k} + 2\qo{{\sf y}_{k-1}} \,\right\}$. 
			Then replace ${\sf x}_k$ (resp.~$\re{{\sf y}_k} + 2\qo{{\sf y}_{k-1}}$) by ${\sf x}_k - {\sf z}_k + {\sf z}_{k-1}$ (resp.~$\re{{\sf y}_k} + 2\qo{{\sf y}_{k-1}} - {\sf z}_k + {\sf z}_{k-1}$), say ${\sf X}_k$ (resp.~$\widetilde{\sf Y}_k$).
		\item[(III)] \label{it:shift ii} Replace ${\sf X}_k$ by $\re{{\sf X}_k} + 2\qo{{\sf X}_{k-1}}$, say $\widetilde{\sf X}_k$.
	\end{enumerate}
	Here we assume that ${\sf x}_{-2}={\sf y}_{-2}={\sf X}_{-2}=0$ as a convention.
	It is easy to see that $\widetilde{\sf X}_k=\widetilde{\sf Y}_k =0 $ unless $ 1 \leq k \leq n$, since ${\sf x}_k={\sf y}_k=0$ for $k \geq n-2$.
	\vskip 2mm
	\item[(3)] 
	Consider $\mc{F}(\bsfT(\bc))= (\widetilde{\sf{X}}_k ; \widetilde{\sf{Y}}_k)_{k \ge -1}$.
	For each $k\ge 1$, attach to the top of $k$-th column of $\tP(\bc)$ a tableau $C_k \in {SST}_{\{i,i+1\}}((1^{c_k}))$, where $c_k=\widetilde{\sf{X}}_k+\widetilde{\sf{Y}}_k$ and the number of $i$'s (resp. $i+1$'s) in $C_k$ is $\widetilde{\sf{X}}_k$ (resp. $\widetilde{\sf{Y}}_k$).
	We denote the resulting tableau by $\tT(\bc)$.
\end{itemize}
}
\end{alg}

\begin{ex} \label{ex: gluing i>m}
{\em 
Suppose $\g = \mf{d}_{2|4}$ and set $i=3$.
Let us consider ${\bf c} ={\bf c}({\bf i}, {\bf j}) \in \bM_{\bt_3}$ where
\begin{equation*}
\begin{split}
	\left(
		\begin{array}{c}
			{\bf i} \\
			{\bf j} 
		\end{array}
	\right)
	&=
	\left(
		\begin{array}{ccccccc}
			3 & 3 & 3 & 3 & 4 & 4 & 5 \\
			4 & 5 & 6 & 6 & 4 & 6 & 5
		\end{array}
	\right).
\end{split}
\end{equation*}
By \eqref{df: type D insertion}, $\burge(\bc)$ coincides with \eqref{ex: D24-1}.

Let us compute $\tT(\bc)$ following Algorithm \ref{alg: Q bar i>m}.
We have $\burge(\bc_{\bt_5})=\vd{5}{5}$ and
\begin{equation*}
	\tP(\bc_{\bt_5},\bc_{{\lozenge}_3}) = 
	\vcenter{\xymatrix@R=-0.82em @C=-0.82em{
	 & \boxed{5} \\
	 & \boxed{6} \\
	\boxed{5} & \boxed{6} \\
	\boxed{5} & \boxed{6} \\
	}}, \qquad
	\tQ(\bc_{\bt_5},\bc_{{\lozenge}_3}) = \vcenter{\xymatrix@R=-0.82em @C=-0.82em{
	& \boxed{3} \\
	& \boxed{3} \\
	\boxed{3} & \\
	\boxed{4} &
	}}.
\end{equation*}
Since $\bc_{\ti}=(c_{(3, 3)},c_{(3, 4)}, c_{(4, 4)})=(0,1,1)$, we have $\bsfT(\bc)$ and $\mc{F}(\bsfT(\bc))$ as shown below.
\begin{equation*}  
	\bsfT(\bc) = 
    \begin{tabular}{cccccccc}
		 $\cdots$ & $3$ & $2$ & $1$ & $0$ & $-1$  \\ \hline        
         $\cdots$ & $0$ & $1$ & $2$ & $1$ & $0$   \\ 
         $\cdots$ & $0$ & $1$ & $0$ & $3$ & $0$  
    \end{tabular}
    \quad
	,
	\quad
	\mc{F}(\bsfT(\bc))=
	\begin{tabular}{ccccccccc}
		 $\cdots$ & $4$ & $3$ & $2$ & $1$ & $0$  & $-1$ \\ \hline        
         $\cdots$ & $0$ & $3$ & $0$ & $1$ & $0$   & $0$ \\ 
         $\cdots$ & $0$ & $1$ & $2$ & $1$ & $0$   & $0$
    \end{tabular}
\end{equation*}
Then $\tT(\bc)$ is obtained by gluing 
\begin{equation*}
\xymatrix@R=-0.82em @C=-0.82em{
	&  & \boxed{3} \\
	&  & \boxed{4} \\
	\boxed{3} & \boxed{4} &  \\
	\boxed{3} & \boxed{4} &  \\
	\boxed{3} &  &  \\
	\boxed{4} &  &  
}
\end{equation*}
to $\tP(\bc_{{\lozenge}_3},\bc_{\bt_5})$, which coincides with \eqref{ex: D24-1}.
}
\end{ex}

Now we have the following.

\begin{prop} \label{prop: gluing algorithm}
For $i \in I \setminus \{ 0 , m \}$ and $\bc \in \bM_\bti$,
we have $\burge(\bc) = \tT(\bc)$.
\end{prop}
\pf
The proof for $1 \leq i < m$ is almost identical to that of \cite[Lemma 5.14(1)]{JK19}, using Lemma \ref{lem: bumping routes}.
The proof for $i > m$ is more involved, so we will present it in Section \ref{subsubsec: proof of gluing alg}.
\qed

\begin{rem} \label{rem: why (3) holds in general}
{\em
When $i > m$, one may check that $\burge( \bc_\ti )=\tT(\bc_\ti)$ is one of the following forms:
\begin{equation} \label{eq: image of abc in type d for i > m}
	\begin{split}
	\qquad
	\scalebox{0.7}{%
    \begin{tabular}{r@{}ll}
    \raisebox{-24ex}{$\ciI\left\{\vphantom{\begin{array}{c}~\\[3ex] ~
    \end{array}}\right.$} & 
	\begin{ytableau}
		\none & \text{\scriptsize $i$} \\
		\none & \raisebox{-0.2em}{$\vdots$} \\
		\none & \text{\scriptsize $i$} \\
		\none & \text{\scriptsize $i$+1} \\
		\none & \raisebox{-0.2em}{$\vdots$} \\
		\none & \text{\scriptsize $i$+1} \\
		\text{\scriptsize $i$} & \text{\scriptsize $i$+1} \\
		\raisebox{-0.2em}{$\vdots$} & \raisebox{-0.2em}{$\vdots$} \\
		\text{\scriptsize $i$} & \text{\scriptsize $i$+1} \\
	\end{ytableau}
	\!\!\raisebox{-2.3ex}{$\left.\vphantom{\begin{array}{c}~\\[3ex] ~
    \end{array}}\right\}2\cii$}
    \shifttext{-50pt}{\raisebox{-13ex}{$\left.\vphantom{\begin{array}{c}~\\[3ex] ~
    \end{array}}\right\}\ciI$}}
    \shifttext{-88pt}{\raisebox{-24ex}{$\left.\vphantom{\begin{array}{c}~\\[3ex] ~
    \end{array}}\right\}2\cII$}}
    \end{tabular}
    }%
    \hspace{-0.5cm}
    \text{ or }
    \quad
    \scalebox{0.7}{%
    \begin{tabular}{r@{}ll}
    \raisebox{-24ex}{$\ciI-1\left\{\vphantom{\begin{array}{c}~\\[3ex] ~
    \end{array}}\right.$} & 
    \ytableausetup {mathmode, boxsize=1.5em} 
	\begin{ytableau}
		\none & \text{\scriptsize $i$} \\
		\none & \raisebox{-0.2em}{$\vdots$} \\
		\none & \text{\scriptsize $i$} \\
		\none & \text{\scriptsize $i$+1} \\
		\none & \raisebox{-0.2em}{$\vdots$} \\
		\none & \text{\scriptsize $i$+1} \\
		\text{\scriptsize $i$} & \text{\scriptsize $i$+1} \\
		\raisebox{-0.2em}{$\vdots$} & \raisebox{-0.2em}{$\vdots$} \\
		\text{\scriptsize $i$} & \text{\scriptsize $i$+1} \\
	\end{ytableau}
	\!\!\raisebox{-2.3ex}{$\left.\vphantom{\begin{array}{c}~\\[3ex] ~
    \end{array}}\right\}2\cii+1$}
    \shifttext{-84pt}{\raisebox{-13ex}{$\left.\vphantom{\begin{array}{c}~\\[3ex] ~
    \end{array}}\right\}\ciI$}}
    \shifttext{-121pt}{\raisebox{-24ex}{$\left.\vphantom{\begin{array}{c}~\\[3ex] ~
    \end{array}}\right\}2\cII$}}
    \end{tabular}
    }%
    \end{split}
    \end{equation}
where $\cii = c_{(i,i)}$, $\ciI = c_{(i,i+1)}$ and $\cII = c_{(i+1,i+1)}$. Recall from \cite[Remark 4.14]{JKU24} that
	\begin{equation} \label{eq: tfi on abc in type d for i > m}
		\tf_i (\cii, \ciI, \cII) = 
		\begin{cases}
			(\cii-1, \ciI+1, \cII) & \text{if $\cii \ge 1$ and $\ciI$ is even,} \\
			(\cii, \ciI-1, \cII+1) & \text{if $\ciI$ is odd with $\ciI \ge 1$,} \\
			{\bf 0} & \text{otherwise.}
		\end{cases}
	\end{equation}
Then it is straightforward to check that $\burge( \tf_i \bc_\ti ) = \tf_i \burge( \bc_\ti )$ by comparing \eqref{eq: image of abc in type d for i > m} and \eqref{eq: tfi on abc in type d for i > m}. 
}
\end{rem}

\begin{lem} \label{lem: tfi on T-tab}
	For $i \in I \setminus \{ 0, m \}$ and $\bc \in \bM_\bti$, we have $\tT( \tf_i \bc ) = \tf_i \tT( \bc )$. 
\end{lem}
\pf
The proof for the case of $1 \leq i < m$ is almost identical to the proof of \cite[Lemma 5.14 (2)]{JK19}\footnote{There is a minor error in the proof of \cite[Lemma 5.14(2)]{JK19}. The corrected proof can be found in arXiv:1810.02103.}, using Lemma \ref{lem: tfi on P-tab and Q-tab}. 
For the case of $i>m$, the proof is given in Section \ref{subsubsec: proof of tfi on Q and Q bar for i > m}.
\qed

\begin{cor} \label{cor: compatibility on bti}
For $i \in I \setminus \{ 0 \}$ and $\bc \in \bM_\bti$, we have $\burge( \tf_i \bc ) = \tf_i \burge( \bc )$.
\end{cor}
\pf
For $i \ne m$, it follows immediately from Proposition \ref{prop: gluing algorithm} and Lemma \ref{lem: tfi on T-tab}.
Now assume $i = m$. 
For $\bc = (c_{(i,j)}) \in \bM_{\bt_m}$, define
\begin{equation*}
	k_0 = \min\{ k \geq m+1 \,|\, c_{(m,k)} = 1 \} \,\,\text{ and }\,\,
	k_1 = \min\{ k \geq m+1 \,|\, c_{(m+1,k)} \neq 0 \},
\end{equation*}
and set $k_a=\infty$ for $a=0,1$ if no such $k$ exists.
If $k_0=\infty$, then $\tf_m \bc={\bf 0}$ and $\tf_m \burge(\bc)={\bf 0}$ since there is no $m$ in $\burge(\bc)$.
So we may assume $k_0<\infty$ and consider the following two cases.

 {\it Case 1}. Suppose $k_1 < k_0$.
By the tensor product rule for $\tf_m$ and \cite[Remark 4.14]{JKU24}, we have $\tf_m \bc = {\bf 0}$.
On the other hand, by Lemma \ref{lem: bumping routes}\eqref{it:bump6} and \eqref{it:bump5}, 
all entries in $\burge(\bc)$ equal to $m$ lies southwest of the entry $m+1$ obtained by the insertion of $\binom{m+1}{k_1}$.
Hence $\tf_m \burge( \bc ) = {\bf 0}$.

 {\it Case 2}. Suppose $k_1 \ge k_0$.
By \cite[Remark 4.14]{JKU24}, $\tf_m(\bc)$ is obtained from $\bc$ by decreasing $\bc_{(m,k)}$ by $1$ and increasing $\bc_{(m+1,k)}$ by $1$.
On the other hand, since $k_1 \ge k_0$, Lemma \ref{lem: bumping routes} \eqref{it:bump5} implies that 
the northeastern most $m$ in $\burge(\bc)$ lies northeast of all $m+1$'s.
Thus, $\tf_m \burge(\bc)$ is obtained from $\burge(\bc)$ by changing this $m$ to $m+1$.
If we apply the reverse algorithm for the inverse of $\burge$ to $\tf_m \burge(\bc)$ 
(see the proof of Proposition \ref{prop:well-definedness})
the entry $m+1$ obtained by applying $\tf_m$ still bumps out $k_0$. 
This implies $\tf_m \burge(\bc) = \burge(\tf_m \bc)$.
\qed

\begin{ex} \label{ex: tfi}
{\em 
Let us illustrate {\it Case 1} and {\it Case 2} in the proof of Corollary \ref{cor: compatibility on bti}.
Suppose that $\g = \mf{d}_{4|3}$ and $i=4$. 
\smallskip

{\it Case 1}. Let us consider $\bc \in \bM_{\bt_4}$ whose biword is given by 
\begin{equation*}
	\left(
		\begin{array}{cccccc}
			4 & \blue{\bf 4} & \red{\bf 5} & 5 & 5 & 6 \\
			7 & \blue{\bf 6} & \red{\bf 5} & 7 & 7 & 6
		\end{array}
	\right).
\end{equation*}
Then $k_0 = 6$, $k_1 = 5$, and we have
\begin{equation*}
\burge(\bc) = 
\vcenter{%
\xymatrix@R=-0.82em @C=-0.82em{
	\boxed{4} & \boxed{\blue{4}} & \boxed{\red{5}} \\
	\boxed{5} & \boxed{6} & \boxed{7} \\
	\boxed{5} & \boxed{6} & \boxed{7} \\
	\boxed{5} & \boxed{6} & \boxed{7} \\
	}
},%
\end{equation*}
where $\red{\bf 5}$ (resp.~$\blue{\bf 4}$) is the recording letter from the insertion of $\red{\bf 5}$ (resp.~$\blue{\bf 6}$). Thus $\tf_4 \burge( \bc ) = {\bf 0}$. \vskip 2mm

{\it Case 2}. Let us consider $\bc \in \bM_{\bt_4}$ and $\tf_4 \bc$ whose biword are given by 
\begin{equation*}
	\left(
		\begin{array}{cccccc}
			4 & \blue{4} & 5 & 5 & 5 & 6 \\
			7 & 6 & 6 & 7 & 7 & 6
		\end{array}
	\right)
	\quad ,\quad
	\left(
		\begin{array}{ccccccc}
			4 & \red{5} & 5 & 5 & 5 & 6 \\
			7 & 6 & 6 & 7 & 7 & 6
		\end{array}
	\right),
\end{equation*}
respectively, where $k_0 = 6$, $k_1 = 6$.
We have
\begin{equation*}
\burge(\bc) = 
\vcenter{%
\xymatrix@R=-0.82em @C=-0.82em{
	\boxed{4} & \boxed{\blue{4}} & \boxed{6} \\
	\boxed{5} & \boxed{6} & \boxed{7} \\
	\boxed{5} & \boxed{6} & \boxed{7} \\
	\boxed{5} & \boxed{6} & \boxed{7} \\
	}
}, \quad \burge(\tf_4 \bc) = 
\vcenter{%
\xymatrix@R=-0.82em @C=-0.82em{
	\boxed{4} & \boxed{\red{5}} & \boxed{6} \\
	\boxed{5} & \boxed{6} & \boxed{7} \\
	\boxed{5} & \boxed{6} & \boxed{7} \\
	\boxed{5} & \boxed{6} & \boxed{7} \\
	}
},
\end{equation*}
from which we see that $\tf_4 \burge(\bc) = \burge(\tf_4 \bc)$.

}
\end{ex}

\subsection{Proof of Theorem \ref{thm: Burge}} \label{sec: proof of Burge}
Let $\bc \in \bM$ and $i \in I$ be given. 
Let $\U(\g_i)$ denote the subalgebra of $\U(\g)$ generated by $e_i,f_i$ and $k_\mu$ for $\mu \in P$.
By Proposition \ref{prop:well-definedness}, it is enough to show that $\burge( \tf_i \bc ) = \tf_i \burge( \bc )$. 

We first suppose that $i \neq 0$. The map
\begin{equation} \label{eq:decomp of bc}
	\bc \mapsto \bc_\btic \ot \bc_\bti
\end{equation}
is an isomorphism of $\U(\g_i)$-crystals with respect to $\te_i$ and $\tf_i$ \cite[Section 4.2]{JKU24}.
By definition, one may write
	\begin{equation} \label{eq:sep of burge}
		\burge(\bc) = \left( \burge(\bc_\bti) \Bins \biw{\ba}{\bb} \right),
	\end{equation}
	where $(\ba, \bb) \in \biws$ such that $\bc_\btic = \bc(\ba, \bb)$. 
	Put $\bV_{\ge i}$ to be the subcrystal of $\bV$ (with respect to $\te_i$ and $\tf_i$) consisting of $T \in \bV$ all of whose entries are greater than or equal to $i$.
	The map $\burge$ can be factorized as follows:
	\begin{equation*} \label{eq:composition of maps to burge}
	\begin{gathered}
		\xymatrix@R=3em @C=3em{
		&  \bM \ar@{->}[d]_{\burge_1}^{\eqref{eq:decomp of bc}} \ar@{->}[r]^\burge & \bV & \\
		& \bM_\btic \ot \bM_\bti \ar@{->}[r]_{\burge_2} & \bM_\btic \ot \bV_{\ge i} \ar@{->}[u]^{\eqref{eq:sep of burge}}_{\burge_3} &
		}
	\end{gathered}
	\end{equation*}
	where $\burge_2$ is defined by $\bc_\btic \ot \bc_\bti \mapsto \bc_\btic \ot \burge(\bc_\bti)$. 
	Note that $\burge_2$ is a\, $\U(\g_i)$-crystal isomorphism by Corollary \ref{cor: compatibility on bti}.
	Since there are no letters $i$ and $i+1$ in $\ba$, if we ignore the letters smaller than $i$, then $\left( \burge(\bc_\bti) \Bins \biw{\ba}{\bb} \right)$ is essentially equal to the Schensted's column insertion.
	Then it follows from \cite[Section 4.4]{BKK} (especially, \cite[Theorem 4.1]{BKK}) that $\burge_3$ is an isomorphism of $\U(\g_i)$-crystals. 
	Hence we have $\burge(\tf_i \bc) = \tf_i \burge( \bc )$ for all $i \neq 0$.

Next suppose that $i = 0$. By the argument in the last paragraph of \cite[Section 5.4]{JK19}, we have the following commuting diagram:
\begin{equation*}
	\xymatrix@R=3em @C=3em{
		\bc \ar@{|->}[r] \ar@{|->}[d] & \bc_{\bt_0 \setminus \bt_2} \ot \bc_{\bt_2} \ar@{|->}[r] & \bc_{\bt_0 \setminus \bt_2} \ot \burge( \bc_{\bt_2} ) \ar@{|->}[d]^{(\dagger)}  \\
		\tT( \bc ) & & \left( \texttt{P}( \bc ), \texttt{Q}( \bc ) \right) \ar@{|->}[ll]_{(\ddagger)\qquad}
	}
\end{equation*}
where we regard $\bc_{\bt_0 \setminus \bt_2}$ by letting  $c_{(0,k)} = 0$ for all $k \ge 1$,
and define $(\tP(\bc),\tQ(\bc))=(\tP_0(\bc),\tQ_0(\bc))$ in the same way as in $i\ne 0$. 
First, we can show that the map $(\dagger)$ commutes with $\tf_0$ by the same argument as in the proof of \cite[Theorem 3.6]{K09}, 
where the set of $(\tP(\bc),\tQ(\bc))$ has a $\U(\g)$-crystal structure (in particular, with respect to $\te_0$ and $\tf_0$) in the sense of \cite[Section 3.2]{K09}. 
Second, we can easily check that the map $(\ddagger)$ commutes with $\tf_0$ by comparing the actions of $\tf_0$ on both sides. 
Finllay, we know $\burge(\bc) = \tT(\bc)$ by Proposition \ref{prop: gluing algorithm}. 
This completes the proof of Theorem \ref{thm: Burge}. 
\qed
	
\begin{rem}
{\em
It is known that $\mc{N}$ is decomposed into irreducible polynomial $\U(\mf{l})$-modules:
\begin{equation*} \label{eq: decomp of N}
\begin{split}
	\mc{N} \cong \bigoplus_{\mu \in \ms{P}_{\mf{d}_{m|n}}} V_{\mf{l}}(\mu).
\end{split}
\end{equation*}
Since the crystals of the above modules (as $\mf{l}$-crystals) are realized as $\bM$ and $\bV$, respectively, the map $\burge$ can be viewed as a crystal-theoretic interpretation of the above decomposition.
}
\end{rem}

\subsection{Proof of Proposition \ref{prop: gluing algorithm} and Lemma \ref{lem: tfi on T-tab}} \label{subsec: combinatorics of Q bar}
The purpose of this section is to complete the proofs of Proposition \ref{prop: gluing algorithm} and Lemma \ref{lem: tfi on T-tab} for the case of $i > m$.
To prove these, it suffices to consider only the numbers of $i$ and $i+1$ in each column of the tableau $\tQ(\bc)$ or $\burge(\bc)$,
which we represent as a two-rowed array with non-negative integers as $\bsfT(\bc)$ in \eqref{eq: df of two array form}.
We assume $i>m$ throughout Section \ref{subsec: combinatorics of Q bar}.

Here we introduce some conventions.
We assume that all two-rowed arrays considered here have only finitely many non-zero entries.
For a given two rowed array 
\begin{equation} \label{eq:two-rowed array} 
		\bsfX = \quad
	\begin{tabular}{ccccc}
		$\cdots$ & $k+1$ & $k$ & $k-1$ & $\cdots$ \\ \hline
		$\cdots$ & ${\sf x}_{k-1}$ & ${\sf x}_{k}$ & ${\sf x}_{k+1}$ & $\cdots$ \\ 
		$\cdots$ & ${\sf y}_{k-1}$ & ${\sf y}_{k}$ & ${\sf y}_{k+1}$ & $\cdots$
	\end{tabular},
\end{equation}
we often write $\bsfX = ({\sf x}_k\,; {\sf y}_k)_{k \geq -1}= ({\sf x}_k\,; {\sf y}_k)$. 
For $\bsfX = ({\sf q}_k\,; {\sf r}_k)$ and $\bsfY = ({\sf s}_k\,; {\sf t}_k)$, we define $\bsfX \pm \bsfY = ({\sf q}_k \pm {\sf s}_k\,; {\sf r}_k \pm {\sf t}_k)$. 
For $a \in \Z_{\geq 0}$ and $ k \geq -1 $, we denote $\bsfx_k(a)$ and $\bsfy_k(a)$ by
\begin{equation*} 
	\bsfx_k(a) = \quad
	\begin{tabular}{ccccc}
	$\cdots$ & $k+1$ & $k$ & $k-1$ & $\cdots$\\ \hline
	$\cdots$ & $0$ & $a$ & $0$ & $\cdots$\\ 
	$\cdots$ & $0$ & $0$ & $0$ & $\cdots$
	\end{tabular}\ \ , \quad \quad
		\bsfy_k(a) = \quad
	\begin{tabular}{ccccc}
	$\cdots$ & $k+1$ & $k$ & $k-1$ & $\cdots$\\ \hline
	$\cdots$ & $0$ & $0$ & $0$ & $\cdots$\\ 
	$\cdots$ & $0$ & $a$ & $0$ & $\cdots$
	\end{tabular} \ \ .
\end{equation*}

Recall that for ${\sf z} \in \Z_+$, we define $\qo{\sf z} \in \Z_+$ and $\re{\sf z} \in \{ 0, 1 \} $ by ${\sf z} = 2\qo{\sf z} + \re{\sf z}$.

\subsubsection{Proof of Proposition \ref{prop: gluing algorithm}} \label{subsubsec: proof of gluing alg}
For $\bc \in \bM_{\bti}$, we prove in this subsection that $\burge(\bc)=\tT(\bc)$.
Recall from Lemma \ref{lem: P-tab} that $\tP(\bc)$ is precisely the subtableau of $\burge(\bc)$ consisting of entries greater than or equal to $i+2$.
Hence, by Algorithm \ref{alg: Q bar i>m}, it remains to show that the number of $i$'s (resp. $i+1$'s) in each column of the tableau $\burge(\bc)$
coincides with $\widetilde{\sf{X}}_k$ (resp. $\widetilde{\sf{Y}}_k$) in $\mc{F}(\bsfT(\bc))= (\widetilde{\sf{X}}_k ; \widetilde{\sf{Y}}_k)_{k \in \Z}$ described in Algorithm \ref{alg: Q bar i>m} (2).

For convenience, we set $\btriangle_j = \emptyset $ for all $ j > m+n$.
The case of $i=m+n-1$ is trivial, and the case of $ m+n-2$ are already treated in \eqref{eq: image of abc in type d for i > m}, since $\triangle_{m+n-2}=\bt_{m+n-2}$.
Let $({\ba}, {\bb}) = (a_1 \dots a_s, b_1 \dots b_s) \in \biws$ be such that $\bc_{\bti \setminus \btiII} = \bc({\ba}, {\bb})$.
By \eqref{df: type D insertion}, we have
\begin{equation*} \label{eq: btiII <-B- ti and pai}
\burge(\bc) = \left( \left( \left( \burge(\bc_\btiII) \Bins \binom{\,a_s\,}{\,b_s\,} \right) \Bins \binom{\,a_{s-1}\,}{\,b_{s-1}\,} \right) \Bins \,\cdots\, \Bins \binom{\,a_1\,}{\,b_1\,} \right).
\end{equation*}
For $1 \leq t \leq s$, we define $(\ba_t , \bb_t)= (a_t \dots a_s , b_t \dots b_s)$ and set $\bc_t = \bc( \ba_t , \bb_t ) \ot \bc_\btiII$.
By convention, we let $\bc_{s+1}=\bc_{\btiII}$.
We proceed by induction on $t$ in decreasing order to prove that $\burge(\bc_t)=\tT(\bc_t)$.

Note that $\burge(\bc_t)= \left( \burge(\bc_{t+1}) \Bins \binom{\,a_t\,}{\,b_t\,} \right)$.
Suppose that $\burge(\bc_{t+1}) = \tT(\bc_{t+1})$.
For simplicity, we write $\tP_t = \tP(\bc_t)$ and $\tQ_t = \tQ(\bc_t)$, and let
\begin{equation*}
\bsfT_t = \bsfT(\bc_t) = ({\sf x}_{t,k} \, ; {\sf y}_{t,k})_{k \in \Z}, \quad
\mc{F}(\bsfT_t) = (\widetilde{{\sf X}}_{t,k} \, ; \widetilde{{\sf Y}}_{t,k})_{k \in \Z}. 
\end{equation*}
To show $\burge(\bc_{t}) = \tT(\bc_{t})$, we analyze $\tT(\bc_{t})$ by considering the difference $\mc{F}(\bsfT_t) - \mc{F}(\bsfT_{t+1})$, since this difference captures the change in the positions of $i$ and $i+1$ between $\tT(\bc_t)$ and $\tT(\bc_{t+1})$.

If $(a_t,b_t)$ lies on ${\lozenge_i}$, then let $k$ be the column index of the unique box in $\tQ_t \setminus \tQ_{t+1}$, which is filled with $a_t$. 
Then it follows from the definition of $\bsfT_t$ in \eqref{eq: df of two array form} that 
\begin{equation*}
\bsfT_t =
\begin{cases}
	\bsfT_{t+1} + \bsfx_k(1) & \text{if $a_t = i$,} \\
	\bsfT_{t+1} + \bsfy_k(1) & \text{if $a_t = i+1$.}
\end{cases}
\end{equation*}
Now we consider five cases according to the position of $(a_t, b_t)$ in $\bti$ as follows. 
\smallskip

{\it Case 1}. Suppose that $a_t=i+1$ and $b_t>i+1$.
Then all the boxes of $\tQ_{t+1}$ are filled with $i+1$, and, by Lemma \ref{lem: bumping routes} \eqref{it:bump5}, their column indices are greater than or equal to $k$.
Thus, we have 
\begin{equation*} 
	\bsfT_{t+1}= \quad
    \begin{tabular}{ccccccc}
         $\cdots$ & $n+1$ & $n$ & $\cdots$ & $k$ & $k-1$ & $\cdots$ \\ \hline
         $\cdots$ & $0$ & $0$ & $\cdots$ & $0$ & $0$ & $\cdots$ \\ 
         $\cdots$ & $0$ & ${\sf y}_{t+1,n}$ & $\cdots$ & ${\sf y}_{t+1,k}$ & $0$ & $\cdots$
    \end{tabular}
	,
\end{equation*}
and $\bsfT_t = \bsfT_{t+1} + \bsfy_k(1)$.
We have $\widetilde{\sf Y}_{t+1,k} = \re{{\sf y}_{t+1,k}}$ by Algorithm \ref{alg: Q bar i>m} (2), which denotes the number of $i+1$'s in the $k$-th column of $\burge(\bc_{t+1})=\tT(\bc_{t+1}) $. 
Consider the following two subcases, illustrating how the bumping route behaves differently.

	{\it Case 1-1}. If $\widetilde{\sf Y}_{t+1,k}=0$, then the bumping route of $b_t$ in $(\burge(\bc_{t+1}) \leftarrow b_t)$ coincides with that in $(\tP_{t+1} \leftarrow b_t)$.
Indeed, the only difference between $\burge(\bc_{t+1})$ and $\tP_{t+1}$ is therein subtableaux of entries $i+1$, 
all of which are located in the columns strictly to the left of the $k$-th column of $\burge(\bc_{t+1})$.
Hence $\burge(\bc_t)= \left( \burge(\bc_{t+1}) \Bins \binom{\,a_t\,}{\,b_t\,} \right)$ is obtained from $\burge(\bc_{t+1})$ 
by adding $b_t$ and then $i+1$ directly above it at the top of the $k$-th column.
On the other hand, we have $\mc{F}(\bsfT_t) - \mc{F}(\bsfT_{t+1}) = \bsfy_k(1) $ since $\widetilde{\sf Y}_{t+1,k}=0$.
This shows that a new letter $i+1$ appears in the $k$-th column of $\tT(\bc_t)$, and thus $\burge(\bc_{t}) = \tT(\bc_{t})$.
	\smallskip
	
	{\it Case 1-2}. If $\widetilde{\sf Y}_{t+1,k}=1$, then there is exactly one $i+1$ in the $k$-th column of $\burge(\bc_{t+1})$.
	Hence the bumping route of $b_t$ in $\left( \burge( \bc_{t+1} ) \leftarrow b_t \right)$ does not terminate in the $k$-th column. 
Instead, it bumps out the $i+1$ from the $k$-th column into the $(k+1)$-th column, 
thereby placing a new $i+1$ at the top of $(k+1)$-th column.
Consequently, one $i+1$ is removed from the $k$-th column while two $i+1$'s are added to the top of the $(k+1)$-th column of $\burge(\bc_{t+1})$. 
On the other hand, since ${\sf y}_{t+1,k}=1 $, immediately implies $\mc{F}(\bsfT_t)-\mc{F}(\bsfT_{t+1})=\bsfy_{k+1}(2)+\bsfy_k(-1)$. 
This shows that $\burge(\bc_{t}) = \tT(\bc_{t})$.
\smallskip

{\it Case 2}. Suppose that $a_t=i+1$ and $b_t=i+1$. In this case, we have
\begin{equation*} 
	\bsfT_{t+1}= \quad 
    \begin{tabular}{ccccccc}
         $\cdots$ & $n+1$ & $n$ & $\cdots$ & $1$ & $0$ & $-1$\\ \hline
         $\cdots$ & $0$ & $0$ & $\cdots$ & $0$ & $0$ & $0$\\ 
         $\cdots$ & $0$ & ${\sf y}_{t+1,n}$ & $\cdots$ & ${\sf y}_{t+1,1}$ & ${\sf y}_{t+1,0}$ & $0$
    \end{tabular},
\end{equation*}
and hence $\bsfT_t = \bsfT_{t+1} + \bsfy_0(2)$. 
Since ${\sf y}_{t+1,0}$ is even, it follows that $\widetilde{\sf Y}_{t+1,0}=0$.
Consequently, $\mc{F}(\bsfT_{t+1})-\mc{F}(\bsfT_{t})=\bsfy_1(2)$.
On the other hand, we have $\burge(\bc_t)=\left( \burge(\bc_{t+1}) \Bins  \binom{\,i+1\,}{\,i+1\,} \right)$.
Since all letters in $\burge(\bc_{t+1})$ are greater than or equal to $i+1$, the bumping route of $i+1$ in $\left( \burge(\bc_{t+1}) \leftarrow i+1 \right)$ terminates in the first column.
Therefore $\burge(\bc_t)$ is obtained from $\burge(\bc_{t+1})$ by adding two entrieds $i+1$ at the top of the first column of $\burge(\bc_{t+1})$. 
This shows that $\burge(\bc_{t}) = \tT(\bc_{t})$.
\smallskip

{\it Case 3}. Suppose that $a_t=i$ and $b_t>i+1$. In this case, we have
\begin{equation*} 
	\bsfT_{t+1}= \quad
    \begin{tabular}{ccccccccc}
         $\cdots$ & $n+1$ & $n$ & $\cdots$ & $k$ & $k-1$ & $\cdots$ & $0$ & $-1$ \\ \hline
         $\cdots$ & $0$ & ${\sf x}_{t+1,n}$ & $\cdots$ & ${\sf x}_{t+1,k}$ & $0$ & $\cdots$ & $0$ & $0$ \\ 
         $\cdots$ & $0$ & ${\sf y}_{t+1,n}$ & $\cdots$ & ${\sf y}_{t+1,k}$ & ${\sf y}_{t+1,k-1}$ & $\cdots$ & ${\sf y}_{t+1,0}$ & $0$
    \end{tabular},
\end{equation*}
and hence $\bsfT_t = \bsfT_{t+1} + \bsfx_k(1)$. 
Consider the following four subcases, illustrating how the bumping route behaves differently.

	{\it Case 3-1}. If $\widetilde{\sf Y}_{t+1,k}=\widetilde{\sf X}_{t+1,k}=0$, then $\mc{F}(\bsfT_{t})-\mc{F}(\bsfT_{t+1})=\bsfx_k(1)$.
	Since the bumping route of $b_t$ in $\left( \burge(\bc_{t+1}) \leftarrow b_t \right)$ terminates in the $k$-th column, $\burge(\bc_t)$ is obtained from $\burge(\bc_{t+1})$ by adding one $i$ to the $k$-th column.
	\smallskip
	
	{\it Case 3-2}. If $\widetilde{\sf Y}_{t+1,k}=0$ and $\widetilde{\sf X}_{t+1,k}=1$,
	then $\mc{F}(\bsfT_{t})-\mc{F}(\bsfT_{t+1})=\bsfx_{k+1}(2)+\bsfx_k(-1)$.
	In this case, the bumping route of $b_t$ in $\left( \burge(\bc_{t+1}) \leftarrow b_t \right)$ terminates in the $(k+1)$-th column, since there is one $i$ on top of the $k$-th column of $\burge(\bc_{t+1})$. 
	As a result, this $i$ is moved from the $k$-th column to the $(k+1)$-th column, and another $i$ is added to the $(k+1)$-th column.
	\smallskip
	
	{\it Case 3-3}. If $\widetilde{\sf Y}_{t+1,k}\geq 1$ and $\widetilde{\sf X}_{t+1,k+1}=0$, then $\mc{F}(\bsfT_{t})-\mc{F}(\bsfT_{t+1})=\bsfx_{k+1}(1)+\bsfy_{k+1}(1)+\bsfy_k(-1)$. 
	Here the bumping route of $b_t$ in $\left(\burge(\bc_{t+1}) \leftarrow b_t \right)$ terminates in the $(k+1)$-th column, since there is at least one $i+1$ on top of the $k$-th column of $\burge(\bc_{t+1})$. 
	Consequently one $i+1$ is moved from the $k$-th column to the $(k+1)$-th column, and one $i$ is added to the $(k+1)$-th column.
	\smallskip
	
	{\it Case 3-4}. If $\widetilde{\sf Y}_{t+1,k}\geq 1$ and $\widetilde{\sf X}_{t+1,k+1} \geq 1$, then $\mc{F}(\bsfT_{t})-\mc{F}(\bsfT_{t+1})=\bsfx_{k+2}(2)+\bsfx_{k+1}(-1)+\bsfy_{k+1}(1)+\bsfy_k(-1)$. 
	In this case, the bumping route of $b_t$ in $\left(\burge(\bc_{t+1}) \leftarrow b_t \right)$ terminates in the $(k+2)$-th column, since one $i+1$ in the $k$-th column of $\burge(\bc_{t+1})$ moves to the $(k+1)$-th column, 
	and this $i+1$ bumps out one $i$ from the top of the $(k+1)$-th column to the $(k+2)$-th column. 
	As a result, one $i+1$ is moved from the $k$-th column to the $(k+1)$-th column, 
	one $i$ is moved from the $(k+1)$-th column to the $(k+2)$-th column, 
	and another $i$ is added to the $(k+2)$-th column.
From {\it Case 3-1}--{\it Case 3-4}, we conlcude that $\burge(\bc_t) = \tT(\bc_t)$.
\medskip

{\it Case 4}. Suppose that $a_t=i$ and $b_t=i+1$. 
In this case, we have
\begin{equation*}
	\bsfT_{t+1}= \quad
    \begin{tabular}{ccccccc}
         $\cdots$ & $n+1$ & $n$ & $\cdots$ & $1$ & $0$ & $-1$\\ \hline
         $\cdots$ & $0$ & ${\sf x}_{t+1,n}$ & $\cdots$ & ${\sf x}_{t+1,1}$ & ${\sf x}_{t+1,0}$ & $0$ \\ 
         $\cdots$ & $0$ & ${\sf y}_{t+1,n}$ & $\cdots$ & ${\sf y}_{t+1,1}$ & ${\sf y}_{t+1,0}$ & ${\sf y}_{t+1,-1}$
    \end{tabular}
	,
\end{equation*}
where ${\sf y}_{t+1,-1} = 2 \qo{{\sf x}_{t+1,0}}$.
We consider two subcases.
\smallskip

	{\it Case 4-1}. If $\widetilde{\sf X}_{t+1,1} =0 $, then $\mc{F}(\bsfT_{t})-\mc{F}(\bsfT_{t+1})=\bsfx_{1}(1)+\bsfy_{1}(1)$.
	In this case, the bumping route of $i+1$ in $\left(\burge(\bc_{t+1}) \leftarrow i+1 \right)$ terminates in the first column. 
	Thus, $\burge(\bc_t)$ is obtained from $\burge(\bc_{t+1})$ by adding two $i+1$'s on top of the first column. 
	
	\smallskip
	
	{\it Case 4-2}. If $\widetilde{\sf X}_{t+1,1} =1 $, then $\mc{F}(\bsfT_{t})-\mc{F}(\bsfT_{t+1})=\bsfx_{2}(2)+\bsfx_{1}(-1)+\bsfy_{2}(1)$, which implies $\burge(\bc_{t}) = \tT(\bc_{t})$.
		In this case, the bumping route of $i+1$ in $\left( \burge(\bc_{t+1}) \leftarrow {i+1} \right) $ terminates in the second column,
		since there is one $i$ on top of the first column.
		Thus, $\burge(\bc_t)$ is obtained from $\burge(\bc_{t+1})$ by adding one $i+1$ to the first column, 
		bumping out $i$ from the first column to the second column, and then adding another $i$ to the second column.
From {\it Case 4-1} and {\it Case 4-2}, we conlcude that $\burge(\bc_t) = \tT(\bc_t)$.
\medskip

{\it Case 5}. Suppose that $a_t=i$ and $b_t=i$. Then $\mc{F}(\bsfT_{t})-\mc{F}(\bsfT_{t+1})=\bsfx_{1}(2)$.
As in {\it Case 2}, $\burge(\bc_t)$ is obtained from $\burge(\bc_{t+1})$ by adding two $i$'s on top of the first column. 
Hence $\burge(\bc_{t}) = \tT(\bc_{t})$.
\smallskip

By {\it Case 1}--{\it Case 5}, we have $\burge(\bc) = \tT(\bc)$ by induction on $s$. 
This completes the proof.
\qed

\begin{ex} 
{\em 
We continue Example \ref{ex: gluing i>m}.
Suppose $\g = \mf{d}_{2|4}$ and set $i=3$.
Let ${\bf c} \in \bM_{\bt_3}$ be given by $\bc_{\bti \setminus \btiII} = \bc( \ba, \bb) = \bc( \, 3333445 \,, \, 4566465 \,) $ and $\bc_{\btiII}=\bc_7=\bc( \,5 \,, \,5 \,)$.
Consider $\burge(\bc_2)= \left( \burge(\bc_{3}) \Bins \binom{\,3\,}{\,5\,} \right)$ where $\bc_3= \bc(\, 33445 \,, \, 66465 \,) $ and $\bc_2= \bc( \, 333445 \,, \,566465 \,) $.
We have
\begin{equation*}
\burge(\bc_3)=
\vcenter{%
\xymatrix@R=-0.82em @C=-0.82em{
	 & \boxed{3} & \boxed{4} \\
	 & \boxed{4} & \boxed{6} \\
	\boxed{3} & \boxed{5} & \boxed{6} \\
	\boxed{4} & \boxed{5} & \boxed{6}
	}
} \, , \quad \quad
\left(\burge(\bc_3) \leftarrow 5\right)=
\vcenter{%
\xymatrix@R=-0.82em @C=-0.82em{
	 & \boxed{4} & \boxed{5} \\
	\boxed{3} & \boxed{4} & \boxed{6} \\
	\boxed{3} & \boxed{5} & \boxed{6} \\
	\boxed{4} & \boxed{5} & \boxed{6}
	}
} \, , \quad \quad
\left( \burge(\bc_{3}) \Bins \binom{\,3\,}{\,5\,} \right)=
\vcenter{%
\xymatrix@R=-0.82em @C=-0.82em{
	\boxed{3} & \boxed{4} & \boxed{5} \\
	\boxed{3} & \boxed{4} & \boxed{6} \\
	\boxed{3} & \boxed{5} & \boxed{6} \\
	\boxed{4} & \boxed{5} & \boxed{6}
	}
}%
.
\end{equation*}
On the other hand, $\bsfT_t$ and $\mc{F}(\bsfT_t)$ for $t=2,3$ are given by
\begin{equation*}
	\bsfT_{3}= \quad 
    \begin{tabular}{c|cccccc}
        & $\cdots$ & $3$ & $2$ & $1$ & $0$ & $-1$ \\ \hline
        3 & $\cdots$ & $0$ & $1$ & $1$ & $0$ & $0$ \\ 
        4 & $\cdots$ & $0$ & $1$ & $0$ & $2$ & $0$
    \end{tabular} \,, \quad
	\bsfT_{2}= \quad 
    \begin{tabular}{c|cccccc}
        & $\cdots$ & $3$ & $2$ & $1$ & $0$ & $-1$ \\ \hline
        3 & $\cdots$ & $0$ & $1$ & $\bf{2}$ & $0$ & $0$ \\ 
        4 & $\cdots$ & $0$ & $1$ & $0$ & $2$ & $0$
    \end{tabular} \, ,
\end{equation*}
\begin{equation*}
	\mc{F}(\bsfT_{3})= \quad 
    \begin{tabular}{c|cccccc}
        & $\cdots$ & $3$ & $2$ & $1$ & $0$ & $-1$ \\ \hline
        3 & $\cdots$ & $1$ & $1$ & $0$ & $0$ & $0$ \\ 
        4 & $\cdots$ & $1$ & $1$ & $1$ & $0$ & $0$
    \end{tabular} \,, \quad
	\mc{F}(\bsfT_{2})= \quad 
    \begin{tabular}{c|cccccc}
        & $\cdots$ & $3$ & $2$ & $1$ & $0$ & $-1$ \\ \hline
        3 & $\cdots$ & $3$ & $0$ & $0$ & $0$ & $0$ \\ 
        4 & $\cdots$ & $1$ & $2$ & $0$ & $0$ & $0$
    \end{tabular} \,.
\end{equation*} \vskip 2mm
\noindent This corresponds to {\it Case 3-4} in the proof above.
We have $\mc{F}(\bsfT_{2})-\mc{F}(\bsfT_{3})=\bsfx_{3}(2)+\bsfx_{2}(-1)+\bsfy_{2}(1)+\bsfy_1(-1)$,
which shows the difference of the column indicies of $3$'s and $4$'s between $\burge(\bc_3)$ and $\burge(\bc_2)$.
}
\end{ex}

\subsubsection{Proof of Lemma \ref{lem: tfi on T-tab}} \label{subsubsec: proof of tfi on Q and Q bar for i > m}
	In this subsection, every two-rowed array $\bsfX=({\sf x}_k \, ; {\sf y}_k)_{k \in \Z}$ is assumed to have non-negative integer entries
	and finite support, that is, ${\sf x}_k={\sf y}_k=0$ for all but finitely many $k \in \Z$.
	The two-rowed array $\bsfT(\bc) $ defined in \eqref{eq: df of two array form} also satisfies this condition once we extend it by setting ${\sf x}_k={\sf y}_k=0$ for $k \leq -2$.
		
	Let us define the actions of $\te_i$ and $\tf_i$ on $\bsfX$, motivated by their actions on a tableau.
	We define a sequence of $+$ and $-$ as follows:
	$$
		\sigma_i(\bsfX) = \left(\, \cdots \, -^{\,{\sf y}_{k+1}}\, +^{\,{\sf x}_{k+1}}\, -^{\,{\sf y}_{k}}\, +^{\,{\sf x}_{k}}\, -^{\,{\sf y}_{k-1}}\, +^{\,{\sf x}_{k-1}}\, \cdots \,\right),
	$$
	where the superscript means the multiplicity of each $\pm$, which has finitely many $\pm$'s since ${\bsfX}$ has finite support. 
	We write $\pm^0 = \cdot$ by convention. 
	We call $\sigma_i(\bsfX)$ the $i$-signature (or signature simply) of $\bsfX$. 
	Write $\sigma_i(\bsfX) = (\xi_1 \,  \xi_2 \, \dots)$, where $\xi_s \in \left\{ -,\,+ \right\}$.
	We replace a pair $(\xi_r, \xi_t) = (+,-)$ by $(\,\cdot\,, \,\cdot\,)$, where $r < t$ and $\xi_s = \cdot$ for $ r < s < t$, and then repeat this process until we get a sequence $\ov{\sigma}_i(\bsfX)$ with no $-$ placed to the right of $+$. We call $\ov{\sigma}_i(\bsfX)$ the reduced $i$-signature of $\bsfX$. 
	Assume that the rightmost $-$ (resp. leftmost $+$) in $\ov{\sigma}_i(\bsfX)$ comes from $-^{{\sf y}_k}$ (resp. $+^{{\sf x}_k}$),
	then $\te_i \bsfX = \bsfX + \bsfx_k(1) - \bsfy_k(1)$ (resp. $\tf_i \bsfX = \bsfX - \bsfx_k(1) + \bsfy_k(1)$).
	In this way, we regard any two-rowed array of the form \eqref{eq:two-rowed array} with non-negative integers as an element of a $\U(\g_i)$-crystal.

	To prove Lemma \ref{lem: tfi on T-tab}, it suffices to show that 
\begin{equation} \label{eq: final}
	\mc{F}(\bsfT(\tf_i \bc)) = \tf_i \mc{F}(\bsfT(\bc))
\end{equation}
	for all $\bc \in \bM_\bti$. 

\begin{lem} \label{lem: sig1} 
For $\bc \in \bM_\bti$, we have
\begin{equation} \label{eq: sig1}
	\bsfT(\tf_i \bc) = \tf_i  \bsfT(\bc)
\end{equation} 
\end{lem}
\pf
If $\bc \in \bM_\ti$, then \eqref{eq: sig1} follows by direct calculation (see Remark \ref{rem: why (3) holds in general}).
If $\bc \in \bM_{\bti \setminus \ti}$, then \eqref{eq: sig1} holds by Lemma \ref{lem: tfi on P-tab and Q-tab}(2).
When $\tf_i$ acts on $\bc \in \bM_\bti$, we may regard $\bc$ as $\bc_{\ti} \ot \bc_{\bti \setminus \ti} $, viewed as an element of a $\U(\g_i)$-crystal. 
Recalling \eqref{eq: df of two array form}, we obtain
$$\bsfT(\bc) \cong \bsfT(\bc)_{\geq 1} \ot \bsfT(\bc)_{\leq 0} \cong \bsfT(\bc_{\bti \setminus \ti}) \ot \bsfT(\bc_{\ti})$$ 
as an element of a $\U(\g_i)$-crystal, 
where $\bsfX_{\geq 1}$ (resp.~$\bsfX_{\leq 0}$) denotes the subarray of $\bsfX$ consisting of the columns with positive (resp.~nonpositive) indices. 
This completes the proof.
\qed

The next lemma follows easily from Algorithm \ref{alg: Q bar i>m} (2).

\begin{lem} \label{lem: sig2}
For $k \in \Z $ and $a \in \Z_{\geq 0}$, let
\begin{equation*}
	\bsfX_k(a) = \bsfx_{k+1}(a) + \bsfy_k(a) =	
    \begin{tabular}{cccccc}
        $\cdots$ & $k+2$ & $k+1$ & $k$   & $k-1$ & $\cdots$ \\ \hline
        $\cdots$ & $0$ & $a$ & $0$ & $0$ & $\cdots$ \\ 
        $\cdots$ & $0$ & $0$ & $a$ & $0$ & $\cdots$
    \end{tabular} \quad .
\end{equation*}
We have
\begin{enumerate}[{\em (1)}]
	\item $\mc{F}(\bsfX_k(2a)) = \bsfX_{k+2}(2a)$,
	\item $\mc{F}(\bsfX + \bsfX_k(2a)) = \mc{F}(\bsfX) + \bsfX_{k+2}(2a)$,
	\item $\tf_i \left( \bsfX + \bsfX_k(2a) \right) = \tf_i \left( \bsfX \right) + \bsfX_k(2a)$,
\end{enumerate}
for any two-rowed array $\bsfX$ with non-negative integers.
\end{lem}
	For $\bsfX = ({\sf x}_k\,; {\sf y}_k)_{k \in \Z}$ with non-negative integer entries, 
	we say that $\bsfX$ is \emph{reduced} if $ \min \left\{\, {\sf x}_{k+1},\, \re{{\sf y}_k} \,\right\} \leq 1$ for all $k\in \Z$.
	Then every $\bsfX$ can be written as 
	$$\bsfX=\bsfX'+\sum_{k\in \Z} \bsfX_k(2a_k) $$
	for some $a_k \in \Z_{\geq 0} $ and a unique reduced two-rowed array $\bsfX'$, where the summation is finite since $\bsfX$ has finite support.
	Hence, by Lemma \ref{lem: sig2}, to prove \eqref{eq: final}, it suffices to consider reduced two-rowed arrays $\bsfX$ only.

\begin{rem} \label{rem: z is 01}
{\em 
Let $\bsfX = ({\sf x}_k\,; {\sf y}_k)_{k \in \Z}$ be a reduced two-rowed array. Then one can check that either $0 \leq {\sf x}_{k+1} \leq 1 $ or $0 \leq {\sf y}_k \leq 1 $ for all $k \in \Z$. 
Also we have ${\sf z}_k=\min \left\{\, {\sf x}_k,\, \re{{\sf y}_k} + 2\qo{{\sf y}_{k-1}} \,\right\} = 0$ or $1$ (recall Algorithm \ref{alg: Q bar i>m} (2)-(II));
If $0 \leq {\sf x}_k \leq 1 $, then clearly ${\sf z}_k \leq {\sf x}_k \leq 1 $.	
If $ 0\leq {\sf y}_{k-1} \leq 1 $, then $2\qo{{\sf y}_{k-1}}=0$ and ${\sf z}_k \leq \re{{\sf y}_k} + 2\qo{{\sf y}_{k-1}} \leq 1$.
These facts will be frequently used in the subsequent proofs.
}	
\end{rem}

The following lemma enables us to describe $\tf_i \mc{F}(\bsfX)$ from $\tf_i \bsfX$ for a reduced $\bsfX$.

\begin{lem} \label{lem: signature lemma} 
Let $\bsfX = ({\sf x}_k\,;{\sf y}_k)_{k \in \Z}$ be a reduced two-rowed array.
\begin{enumerate}[{\em (1)}]
	\item 
	If $\sigma_i(\mc{F}(\bsfX))$ has an adjacent pair $+-$ coming from $+^{\bsfx_k} -^{\bsfy_{k-1}}$,
	then $\sigma_i(\mc{F}(\bsfX))$ has either the same adjacent pair or $++--$, 
	where the inner pair $+-$ comes from $+^{\bsfx_{k+1}} -^{\bsfy_{k}}$, 
	and the outer pair $+-$ comes from $ +^{\bsfx_k} -^{\bsfy_{k-1}}$, both as seen in $\sigma_i(\bsfX)$.

	\item The $+$'s in $\sigma_i(\mc{F}(\bsfX))$ coming from $+^{\bsfx_k}$ in $\sigma_i(\bsfX)$ 
	cannot appear to the left of any $-$'s in $\sigma_i(\mc{F}(\bsfX))$ coming from $-^{\bsfy_l}$ in $\sigma_i(\bsfX)$ for all $l \ge k+1$.
	If such $+$'s appear to the left of the $-$'s coming from $-^{\bsfy_k}$, 
	then this $-$ must be unique and canceled by the adjacent $+$ arsing from $+^{{\sf x}_{k+1}}$.
	
	\item The $-$'s in $\sigma_i(\mc{F}(\bsfX))$ coming from $-^{\bsfy_{k-1}}$ in $\sigma_i(\bsfX)$ 
	cannot appear to the left of any $+$'s in $\sigma_i(\mc{F}(\bsfX))$ coming from $+^{\bsfx_l}$ in $\sigma_i(\bsfX)$ for all $l \ge k+1$.
	If such $-$'s appear to the left of the $+$'s coming from $+^{\bsfx_k}$, 
	then this $-$ must be unique and canceled by the adjacent $+$ coming from $+^{{\sf x}_{k}}$.

	\item We have $$\ov{\sigma}_i( \mc{F}( \bsfX ) ) = \ov{\sigma}_i( \bsfX ).$$
		Furthermore, the leftmost $+$ in $\ov{\sigma}_i(\mc{F}(\bsfX))$ comes from the leftmost $+$ in $\ov{\sigma}_i(\bsfX)$.
\end{enumerate}
\end{lem}
\begin{ex}\label{ex: signature lemma}
{\em 
	(1) Let $\bsfX=({\sf x}_k\,;{\sf y}_k)_{k \in \Z}$ be given by
	\begin{equation*}
    \begin{tabular}{cccccc}
        $\cdots$ & $2$ & $1$ & $0$   & $-1$ & $\cdots$ \\ \hline
        $\cdots$ & $0$ & $3$ & $2$ & $0$ & $\cdots$ \\ 
        $\cdots$ & $0$ & $2$ & $1$ & $0$ & $\cdots$
    \end{tabular}
	\ .
	\end{equation*}
	Then the $i$-signature of $\bsfX$ is $\sigma_i(\bsfX)=(-^2 +^3 -^1 +^2)$, and the reduced signature is $\ov{\sigma}_i(\bsfX)=(-^2 +^2 (+ -) +^2)=(-^2 +^4)$,
	where the canceled $(+-)$ comes from $+^{{\sf x}_1=3} $ and $-^{{\sf y}_0=1}$.
	Next, we consider $\mc{F}(\bsfX)$.
	Applying (I), we obtain
	\begin{equation*}
		\begin{tabular}{cccccc}
			$\cdots$ & $2$ & $1$ & $0$   & $-1$ & $\cdots$ \\ \hline
			$\cdots$ & $0$ & $3$ & $2$ & $0$ & $\cdots$ \\ 
			$\cdots$ & $2$ & $0$ & $1$ & $0$ & $\cdots$
		\end{tabular}
		\ .
	\end{equation*}
	Applying (II) then gives
	\begin{equation*}
		\begin{tabular}{cccccc}
			$\cdots$ & $2$ & $1$ & $0$   & $-1$ & $\cdots$ \\ \hline
			$\cdots$ & $0$ & $4$ & $1$ & $0$ & $\cdots$ \\ 
			$\cdots$ & $2$ & $1$ & $0$ & $0$ & $\cdots$
		\end{tabular}
		\ .
	\end{equation*}
	Finally, applying (III), 
	we find that $\mc{F}(\bsfX)=(\widetilde{\sf X}_k\,;\widetilde{\sf Y}_k)_{k \in \Z}$ is given by
	\begin{equation*}
		\begin{tabular}{cccccc}
			$\cdots$ & $2$ & $1$ & $0$   & $-1$ & $\cdots$ \\ \hline
			$\cdots$ & $4$ & $0$ & $1$ & $0$ & $\cdots$ \\ 
			$\cdots$ & $2$ & $1$ & $0$ & $0$ & $\cdots$
		\end{tabular}
		\ .
	\end{equation*}		
	We then compute $\sigma_i(\mc{F}(\bsfX))=(-^2 +^4 -^1 +^1)$ and $\ov{\sigma}_i(\mc{F}(\bsfX))=(-^2 +^3 (+ -) +^1)=(-^2 +^4) = \ov{\sigma}_i(\bsfX)$,
	showing that Lemma \ref{lem: signature lemma} (4) holds for this choice of $\bsfX$.	
	
	We observe that $+^{\widetilde{\sf X}_2=4}$ in $\sigma_i(\mc{F}(\bsfX))$ is made of $+^3$ from $+^{{\sf x}_1=3}$ and $+^1$ in $\sigma_i(\bsfX)$ from $+^{{\sf x}_0=2}$.
	Here, the canceled $+$ in $\ov{\sigma}_i(\mc{F}(\bsfX))$, which is one of the $+$'s in $+^{\widetilde{\sf X}_2=4}$, is interpreted as the one coming from $+^{{\sf x}_1}$ in $\sigma_i(\bsfX)$.
	This indicates that the canceled pair $(+-)$ in $ \ov{\sigma}_i(\mc{F}(\bsfX))$ originates from the same canceled pair in $\ov{\sigma}_i(\bsfX)$,
	thereby verifying Lemma \ref{lem: signature lemma} (1) holds for this example.
	
	On the other hand, there exists a $+$ in ${\sigma}_i(\mc{F}(\bsfX))$ coming from $+^{{\sf x}_1}$ in ${\sigma}_i(\bsfX)$
	which appears to the left of a $-$ in ${\sigma}_i(\mc{F}(\bsfX))$ coming from $-^{{\sf y}_1}$ in $\sigma_i(\bsfX)$.
	However, this $+$ can be regarded as canceled with $-$ in ${\sigma}_i(\mc{F}(\bsfX))$ coming from $-^{{\sf y}_0}$ in $\sigma_i(\bsfX)$,
	thereby satisfying Lemma \ref{lem: signature lemma} (2).
	Accordingly, this $+$ is canceled in $\ov{\sigma}_i(\mc{F}(\bsfX))$.

	Finally, one can similarly verify that Lemma \ref{lem: signature lemma} (3) also holds for this example.
	\smallskip
	
	(2) Let $\bsfX=({\sf x}_k\,;{\sf y}_k)_{k \in \Z}$ be given by
	\begin{equation*}
    \begin{tabular}{cccccc}
        $\cdots$ & $2$ & $1$ & $0$   & $-1$ & $\cdots$ \\ \hline
        $\cdots$ & $0$ & $1$ & $1$ & $0$ & $\cdots$ \\ 
        $\cdots$ & $0$ & $0$ & $1$ & $1$ & $\cdots$
    \end{tabular}
	\ .
	\end{equation*}
	Then $\mc{F}(\bsfX)=(\widetilde{\sf X}_k\,;\widetilde{\sf Y}_k)_{k \in \Z}$ is given by
	\begin{equation*}
		\begin{tabular}{cccccc}
			$\cdots$ & $2$ & $1$ & $0$   & $-1$ & $\cdots$ \\ \hline
			$\cdots$ & $2$ & $0$ & $0$ & $0$ & $\cdots$ \\ 
			$\cdots$ & $0$ & $1$ & $0$ & $1$ & $\cdots$
		\end{tabular}
		\ .
	\end{equation*}
	We then compute $\sigma_i(\bsfX)=(+-+-)$, $\sigma_i(\mc{F}(\bsfX))=(++--)$, and $\ov{\sigma_i}(\bsfX)=\ov{\sigma_i}(\mc{F}(\bsfX))=(\cdot)$.
	Here, we interpret the second $+$ in $\sigma_i(\mc{F}(\bsfX))$ as coming from $+^{{\sf x}_1}$ in $\sigma_i(\bsfX)$.
	Then the inner pair $+-$ in $\sigma_i(\mc{F}(\bsfX))=(++--)$ comes from $+^{{\sf x}_1}$ and $-^{{\sf y}_0}$, 
	while the outer pair $+-$ comes from $+^{{\sf x}_0}$ and $-^{{\sf y}_{-1}}$.
	This provides an example of the second case described in (1) of Lemma \ref{lem: signature lemma}.
	}
\end{ex}
\pf
We will prove (1), (2) and (4) since the proof of (3) is essentially identical to that of (2).
\smallskip

(1)	Assume that there exists an adjacent pair $+-$ coming from $+^{\bsfx_k}$ and $-^{\bsfy_{k-1}}$ in $\sigma_i(\bsfX)$.
Since $\bsfX$ is reduced, we have $\mathrm{min}({\sf x}_k, {\sf y}_{k-1})=1$.
Hence, there are two possible cases: either ${\sf x}_k \geq {\sf y}_{k-1} = 1$ or ${\sf x}_k = 1 < {\sf y}_{k-1}$.
\smallskip

{\it Case 1}. Suppose ${\sf x}_k \geq {\sf y}_{k-1} = 1$.
Note that under $\mc{F}$, each $+$ or $-$ in $\sigma_i(\bsfX)$ may be shifted by at most two positions.
In particular, the $+$ (resp. $-$) coming from $+^{{\sf x}_k}$ (resp. $-^{{\sf y}_{k-1}}$) may be moved to $+^{\widetilde{\sf X}_k}$, $+^{\widetilde{\sf X}_{k+1}}$, or $+^{\widetilde{\sf X}_{k+2}}$ 
(resp. $-^{\widetilde{\sf Y}_{k-1}}$, $-^{\widetilde{\sf Y}_{k}}$, or $-^{\widetilde{\sf Y}_{k+1}}$) in $ \sigma_i(\mc{F}(\bsfX))$.
We first focus on the $-=-^{{\sf y}_{k-1}}$ in $ \sigma_i(\mc{F}(\bsfX))$.
By Algorithm \ref{alg: Q bar i>m} (2), this $-$ can be moved by at most one position, since ${\sf y}_{k-1} = 1$.

Let us consider two subcases:

{\it Case 1-1}. Suppose that $-=-^{{\sf y}_{k-1}}$ is moved to $-^{\widetilde{\sf Y}_{k}}$ in $ \sigma_i(\mc{F}(\bsfX))$ by $\mc{F}$.
This implies ${\sf z}_{k-1}=1=\mathrm{min}\{{\sf x}_{k-1},1+2\qo{{\sf y}_{k-2}}\}$, and hence ${\sf x}_{k-1}\geq 1 $.
We now examine the positions of the $+$'s in $ \sigma_i(\mc{F}(\bsfX))$ that comes from $+^{{\sf x}_{k}}$ .
If ${\sf X}_k={\sf x}_k-{\sf z}_k+1 \geq 2$, then at least one $+$ from $+^{{\sf x}_k}$ is moved to $+^{\widetilde{{\sf X}_{k+1}}}$, so (1) holds.
Otherwise, suppose ${\sf x}_k-{\sf z}_k=0$. 
Then ${\sf x}_k={\sf z}_k=1$ and $\re{{\sf y}_k}={\sf z}_k=1$, since $\qo{{\sf y}_{k-1}}=0$.
If $\re{{\sf X}_{k+1}}=\re{{\sf x}_{k+1}-{\sf z}_{k+1}+1}=1$, then one $+$ from $+^{{\sf x}_k}$ remains at $+^{\widetilde{{\sf X}_{k+1}}}$, and again (1) holds.
Otherwise, suppose $\re{{\sf x}_{k+1}-{\sf z}_{k+1}+1}=0$. 
Then ${\sf x}_{k+1}-{\sf z}_{k+1}\geq 1$ and ${\sf y}_{k+1}-{\sf z}_{k+1}=0$.
Thus, in $\sigma_i(\mc{F}(\bsfX))$, we obtain:
\begin{itemize}
	\item $+^{\widetilde{{\sf X}_{k+2}}}$ contains at least two $+$'s, each coming from $+^{{\sf x}_{k}}$ and $+^{{\sf x}_{k+1}}$,
	\item $\widetilde{{\sf Y}_{k+1}}=1$, originating from ${\sf y}_k=1$,
	\item $\widetilde{{\sf X}_{k+1}}=0$,
	\item $\widetilde{{\sf Y}_{k}}=1$, originating from ${\sf y}_{k-1}=1$.
\end{itemize}
This is precisely the second case of (1), namely ($++--$), verifying (1) in this subcase.

{\it Case 1-2}. Suppose instead that $-=-^{{\sf y}_{k-1}}$ is not moved by $\mc{F}$.
	This implies ${\sf z}_{k-1}=0=\mathrm{min}({\sf x}_{k-1},1+2\qo{{\sf y}_{k-2}})$, hence ${\sf x}_{k-1}=0 $.
	If $\re{{\sf x}_k-{\sf z}_k}=1$, then one $+$ from $+^{{\sf x}_k}$ remains at $+^{\widetilde{{\sf X}_{k}}}$, verifying (1) in this subcase.
	Otherwise, suppose $\re{{\sf x}_k-{\sf z}_k}=0$. 
	Then ${\sf x}_k={\sf z}_k=1$ and $\re{{\sf y}_k}-{\sf z}_k=0$. 
	In $\sigma_i(\mc{F}(\bsfX))$ we then have:
\begin{itemize}
	\item $+^{\widetilde{{\sf X}_{k+1}}}$ contains at least one $+$ coming from $+^{{\sf x}_{k}}$,
	\item $\widetilde{{\sf Y}_{k}}=0$, since $\re{{\sf y}_k}-{\sf z}_k=0$,
	\item $\widetilde{{\sf X}_{k}}=0$, since $\re{{\sf x}_k-{\sf z}_k}=0$ and ${\sf X}_{k-1}={\sf z}_{k-2} \leq 1$,
	\item $-^{\widetilde{{\sf Y}_{k-1}}}$ contains at least one $-$ coming from $-^{{\sf y}_{k-1}}$.
\end{itemize}
Thus, $+=+^{{\sf x}_{k}}$ and $-=-^{{\sf y}_{k-1}}$ remain adjacent in $\sigma_i(\mc{F}(\bsfX))$, verifying (1).
\smallskip

{\it Case 2}. Suppose ${\sf x}_k = 1 < {\sf y}_{k-1}$.
We focus on the $+$ coming from $+^{{\sf x}_{k}}$ in $ \sigma_i(\mc{F}(\bsfX))$.
Since $\qo{{\sf y}_{k-1}}\geq 1$, we have ${\sf z}_k=\mathrm{{\sf x}_k,\re{{\sf y}_k}+2\qo{{\sf y}_{k-1}}}=1$.
Therefore, this $+$ is moved to either $+^{\widetilde{\sf X}_{k+2}}$ or $+^{\widetilde{\sf X}_{k+1}}$ by $\mc{F}$ under $\mc{F}$.

Suppose first that such $+$ in $+^{{\sf x}_{k}}$ is moved to $+^{\widetilde{\sf X}_{k+2}}$.
This means ${\sf z}_k=1$, so the $+$ is moved to $+^{\widetilde{\sf X}_{k+1}}$ by (II),
and since $\re{{\sf x}_{k+1}-{\sf z}_{k+1}+1}=0$, it is further moved to $+^{\widetilde{\sf X}_{k+2}}$ by (III).
From $\re{{\sf x}_{k+1}-{\sf z}_{k+1}+1}=0$ we deduce ${\sf x}_{k+1}-{\sf z}_{k+1} \geq 1$, which implies ${\re{{\sf y}_{k+1}}}+2\qo{{\sf y}_{k}}-{\sf z}_{k+1}=0$, hence $\qo{{\sf y}_{k}}=0$.
Thus, in $\sigma_i(\mc{F}(\bsfX))$ we have the following:
\begin{itemize}
	\item $+^{\widetilde{{\sf X}_{k+2}}}$ contains at least two $+$'s, each coming from $+^{{\sf x}_{k}}$ and $+^{{\sf x}_{k+1}}$,
	\item $\widetilde{{\sf Y}_{k+1}}=1={\sf z}_k$, originating from either ${\sf y}_k$ or ${\sf y}_{k-1}$,
	\item $\widetilde{{\sf X}_{k+1}}=0$, since $\re{{\sf x}_{k+1}-{\sf z}_{k+1}+1}=0$ and ${\sf X}_k={\sf z}_{k-1}\leq 1$,
	\item $\widetilde{{\sf Y}_{k}}=\re{{\sf y}_k}+2\qo{{\sf y}_{k-1}}-1+{\sf z}_{k-1}$.
\end{itemize}
	If $\widetilde{{\sf Y}_{k+1}}=1$ originates from ${\sf y}_{k-1}$, then (1) holds.
	Otherwise, suppose it originates from ${\sf y}_{k}$, so that $\re{{\sf y}_k}=1$.
	Then $-^{\widetilde{{\sf Y}_{k+1}}}$ coming from $-^{{\sf y}_k}$ is canceled with one $+$ in $+^{\widetilde{{\sf X}_{k+2}}}$ coming from $+^{{\sf x}_{k+1}}$.
	Since $ 2\qo{{\sf y}_{k-1}} \geq 2$, we have $\widetilde{{\sf Y}_{k}}\geq 1$, so $-^{\widetilde{{\sf Y}_{k}}}$ contains at least one $-$ from $-^{{\sf y}_{k-1}}$. 
	Hence (1) holds.
	\smallskip
	
	Next suppose that such $+$ in $+^{{\sf x}_{k}}$ is moved to $+^{\widetilde{\sf X}_{k+1}}$ under $\mc{F}$.
	Then for this $+$ is moved by (II), we must have ${\sf z}_k=1$ and $\re{{\sf x}_{k+1}-{\sf z}_{k+1}}=0$.
	Since $\widetilde{{\sf Y}_k}=\re{{\sf y}_k} + 2 \qo{{\sf y}_{k-1}} - 1 + {\sf z}_{k-1} $ and $ 2\qo{{\sf y}_{k-1}} \geq 2$,
	it follows that $-^{\widetilde{{\sf Y}_k}} $ contains at least one $-$ from $-^{{\sf y}_{k-1}}$. 
	Hence (1) holds.
\smallskip

(2) Let us refer to Algorithm \ref{alg: Q bar i>m} (2)-(I),(II),(III) simply as (I),(II),(III).
Let us focus on the location of $+^{{\sf x}_k}$ under (I),(II) and (III).
According to Algorithm \ref{alg: Q bar i>m} (2), $+^{{\sf x}_k}$ can be moved at most two positions to the left, 
so it cannot be positioned to the left of $-^{{\sf y}_l}$ for any $l \geq k+2$.
Thus it suffices to compare its position with $-^{{\sf y}_{k+1}}$ and $-^{{\sf y}_k}$.
If all $+$'s of $+^{{\sf x}_k}$ remain to the left of $+^{{\widetilde{\sf X}_{k}}}$ in $\sigma_i(\mc{F}(\bsfX))$, then there is nothing to prove.
Otherwise, we consider the following two cases.
\medskip

{\it Case 1}. There is a $+$ of $+^{{\sf x}_k}$ in $\sigma_i(\bsfX)$ that is moved to $+^{{\widetilde{\sf X}_{k+2}}}$ in $\sigma_i(\mc{F}(\bsfX))$.

First, at least one $+$ of $+^{{\sf x}_k}$ is moved by (II), so that ${\sf z}_k=1$ by Remark \ref{rem: z is 01}.
This $+$ is then moved again by (III),
which implies ${\sf X}_{k+1} = {\sf x}_{k+1} - {\sf z}_{k+1} + {\sf z}_k \ge 2$ and $\re{{\sf X}_k} =0$.
Since ${\sf x}_{k+1} - {\sf z}_{k+1} \ge 1 $, we have $\re{{\sf y}_{k+1}}+2\qo{{\sf y}_k} -{\sf z}_{k+1}=0 $ by the definition of ${\sf z}_{k+1}$.

From this, every $-^{{\sf y}_{k+1}}$ in $\sigma_i(\bsfX)$ is moved to either $-^{{\widetilde{\sf Y}_{k+2}}}$ or $-^{{\widetilde{\sf Y}_{k+3}}}$ in $\sigma_i(\mc{F}(\bsfX))$.
Hence, no $+$ coming from $+^{{\sf x}_k}$ appear to the left of any $-$ coming from in $\sigma_i(\bsfX)$.

Next, consider the position of $-$'s coming from $-^{\bsfy_{k}}$.
We have $2 \qo{{\sf y}_k}=0$ since $\re{{\sf y}_{k+1}}+2\qo{{\sf y}_k} -{\sf z}_{k+1}=0 $ and ${\sf z}_{k+1} \leq 1$.
Assume ${\sf y}_k=1$, since if ${\sf y}_k=0$ there is nothing to prove.
This $-=-^{{\sf y}_k}$ is moved to $-^{{\widetilde{\sf Y}_{k+1}}}$ by (II) since ${\sf z}_k=1$,
so there exists a $+$ coming from $+^{{\sf x}_{k}}$ that appears to the left of this $-=-^{{\sf y}_{k}}$ in $\sigma_i(\mc{F}(\bsfX))$.

We then show that this $-$ is canceled by $+$ coming from $+^{{\sf x}_{k+1}}$.
Since ${\sf x}_{k+1} - {\sf z}_{k+1} + {\sf z}_k \ge 2$ we have ${\sf x}_{k+1}\geq 1$, 
so there exists an adjacent pair $+-$ coming from $+^{{\sf x}_k} -^{{\sf y}_{k-1}}$ in $\sigma_i(\bsfX)$.
The $+$ from $+^{{\sf x}_{k+1}}$ ends at $+^{\widetilde{\sf X}_{k+2}}$ in $\sigma_i(\mc{F}(\bsfX))$,
since ${\sf X}_{k+1} = {\sf x}_{k+1} - {\sf z}_{k+1} + {\sf z}_k \ge 2$.
Therefore the $-$ coming from $ -^{{\sf y}_k}$ is canceled by adjacent $+$ coming from $ +^{{\sf x}_{k+1}}$.
\medskip

{\it Case 2}. There is a $+$ from $+^{{\sf x}_k}$ in $\sigma_i(\bsfX)$ that is moved to $+^{{\widetilde{\sf X}_{k+1}}}$ in $\sigma_i(\mc{F}(\bsfX))$, 
and no $+$ is moved to $+^{{\widetilde{\sf X}_{k+1}}}$ in $\sigma_i(\mc{F}(\bsfX))$.
To derive a contradiction, assume that a $-$ coming from $-^{{\sf y}_k}$ remains at $-^{\widetilde{{\sf Y}_{k}}}$ in $\sigma_i(\mc{F}(\bsfX))$.
Then $\re{{\sf y}_k}=1$ and ${\sf z}_k=0$.
By the definition of ${\sf z}_k$, we have ${\sf x}_k=0$, so there is nothing to prove.

\smallskip

 	(4)
	Let $\sigma_i'$ be the sequence obtained from $\sigma_i(\bsfX)$ 
	by canceling only the adjacent pairs $+-$ coming from $+^{{\sf x}_k}  -^{{\sf y}_{k-1}}$, for all $k \in \Z$.
	Next, let $\sigma_i''$ be the sequence obtained from $\sigma_i(\mc{F}(\bsfX))$
	by canceling only the pairs $+-$ induced by the adjacent pairs $+-$ in $\sigma_i(\bsfX)$,
	which is possible by (1).
	By (2) and (3), we have $\sigma_i'=\sigma_i''$, since any uncanceled $+$'s in $\sigma_i''$ come from $+^{{\sf x}_k}$ 
	located to the left of $+^{{\sf x}_{k-1}}$ and $-^{{\sf y}_{k-1}}$, 
	and to the right of $+^{{\sf x}_{k+1}}$ and $-^{{\sf y}_{k}}$.
	Therefore, the reduced sequence $\ov{\sigma}_i(\bsfX)$ and $\ov{\sigma}_i(\mc{F}(\bsfX))$ coincide. This proves (4).
\qed
\medskip

Now we are ready to prove \eqref{eq: final}. 
\smallskip

For $\bc \in \bM_\bti$, we have $\mc{F}\left( \tf_i \left( \bsfT(\bc) \right)  \right) = \tf_i\, \mc{F}\left( \bsfT(\bc) \right)$.
By Lemma \ref{lem: sig1} and \ref{lem: sig2}, it suffices to verify that $\mc{F}(\tf_i(\bsfX))=\tf_i \mc{F}(\bsfX)$ for any reduced two-rowed array $\bsfX$. 
Assume that $\tf_i \bsfX = \bsfX - \bsfx_k(1) + \bsfy_k(1)$. 
In other words, the left-most uncanceled $+$ in $\ov{\sigma_i}(\bsfX)$ comes from $+^{{\sf x}_k}$ in $\sigma_i(\bsfX)$.
Then ${\sf x}_{k} > {\sf y}_{k-1}$, since at least one $+$ coming from $+^{{\sf x}_k}$ is not canceled, so ${\sf y}_{k-1} \leq 1$.
Similarly, ${\sf x}_{k+1} \leq {\sf y}_{k}$, since the left-most uncanceled $+$ comes from $+^{{\sf x}_k}$, hence ${\sf x}_{k+1} \leq 1$.

Now, depending on the number of positions an uncanceled $+$ coming from $+^{{\sf x}_k}$ is moved under $\mc{F}$, we consider the following three cases.

\medskip
{\it Case 1}.
Suppose that there is a $+$ of $+^{{\sf x}_{k}}$ moved to $+^{\widetilde{\sf X}_{k+2}}$ in $ \sigma_i(\mc{F}(\bsfX))$ by $\mc{F}$.
Then ${\sf z}_k = 1$ and $\re{{\sf x}_{k+1}-{\sf z}_{k+1}+1}=0$.
Since $\qo{{\sf y}_{k-1}}=0$ and ${\sf z}_k=1$, we have $\re{{\sf y}_k}=1$.
Moreover, as ${\sf x}_{k+1}-{\sf z}_{k+1}\geq 1$ and ${\sf x}_{k+1} \leq 1$, we deduce ${\sf x}_{k+1}=1$, ${\sf z}_{k+1}=0$, and $\re{{\sf y}_{k+1}}=\qo{{\sf y}_k}=0$.
Hence, 
\begin{itemize}
	\item $\widetilde{{\sf X}_{k+2}}=\re{({\sf x}_{k+2} - {\sf z}_{k+2})}+2$, where $+2$ is the sum of ${\sf x}_{k+1}=1$ and $1$ from ${\sf x}_{k}$,
	\item $\widetilde{{\sf Y}_{k+1}}=1$, originating from ${\sf y}_k=1$,
	\item $\widetilde{{\sf X}_{k+1}}=2\qo{({\sf x}_k-1+{\sf z}_{k-1})}$,
	\item $\widetilde{{\sf Y}_{k}}={\sf z}_{k-1}$.
\end{itemize}

We claim that the $+$ in $+^{\widetilde{{\sf X}_{k+2}}}$ coming from $+^{{\sf x}_k}$ is not canceled.
Indeed, $\re{({\sf x}_{k+2} - {\sf z}_{k+2})}=0$ by Lemma \ref{lem: signature lemma}(4).
The $+$ in $+^{\widetilde{{\sf X}_{k+2}}}$ coming from $+^{{\sf x}_{k+1}}$ is canceled with $-$ in $-^{\widetilde{{\sf Y}_{k+1}}}$. 
Consider two subcases:
If ${\sf y}_{k-1}=0$, the claim follows by Lemma \ref{lem: signature lemma}(4).
If ${\sf y}_{k-1}=1$, then ${\sf z}_{k-1}=1$ and ${\sf x}_k-1+{\sf z}_{k-1} \geq 2$, 
so $-^{\widetilde{{\sf Y}_{k}}}$ is canceled by $+^{\widetilde{{\sf X}_{k+1}}}$, 
and the claim holds by Lemma \ref{lem: signature lemma}(4).

Thus the claim is verified.

Summarizing, in {\it Case 1}, we have
\begin{equation*}
	\bsfX= \quad
	\begin{tabular}{cccccc}
		$\cdots$ & $k+2$ & $k+1$ & $k$   & $k-1$ & $\cdots$ \\ \hline
		$\cdots$ & ${\sf x}_{k+2}$ & $1$ & ${\sf x}_{k}$ & ${\sf x}_{k-1}$ & $\cdots$ \\ 
		$\cdots$ & ${\sf y}_{k+2}$ & ${\sf y}_{k+1}$ & $1$ & ${\sf y}_{k-1}$ & $\cdots$
	\end{tabular}
	\ ,
\end{equation*}
and 
\begin{equation*}
	\mc{F}(\bsfX)= \quad
	\begin{tabular}{cccccc}
	$\cdots$ & $k+2$ & $k+1$ & $k$   & $k-1$ & $\cdots$ \\ \hline
	$\cdots$ & $2$ & $2\qo{({\sf x}_k-1+{\sf z}_{k-1})}$ & $\widetilde{{\sf X}_k}$ & $\widetilde{{\sf X}_{k-1}}$ & $\cdots$ \\ 
	$\cdots$ & $\widetilde{{\sf Y}_{k+2}}$ & $1$ & ${\sf z}_{k-1}$ & $\widetilde{{\sf Y}_{k-1}}$ & $\cdots$
	\end{tabular}
	\ ,
\end{equation*}
with $\widetilde{{\sf X}_k}=\re{{\sf x}_k-{\sf z}_k+{\sf z}_{k-1}} + 2\qo{({\sf x}_{k-1}-{\sf z}_{k-1}+{\sf z}_{k-2})}$ and $\widetilde{{\sf Y}_{k}}=\re{{\sf y}_{k}}+2\qo{{\sf y}_{k-1}}-{\sf z}_{k}+{\sf z}_{k-1}$.
Then, by the claim, 
\begin{equation*}
	\tf_i \mc{F}(\bsfX) = \quad
	\begin{tabular}{cccccc}
		$\cdots$ & $k+2$ & $k+1$ & $k$   & $k-1$ & $\cdots$ \\ \hline
		$\cdots$ & $1$ & $2\qo{({\sf x}_k-1+{\sf z}_{k-1})}$ & $\widetilde{{\sf X}_k}$ & $\widetilde{{\sf X}_{k-1}}$ & $\cdots$ \\ 
		$\cdots$ & $\widetilde{{\sf Y}_{k+2}}+1$ & $1$ & ${\sf z}_{k-1}$ & $\widetilde{{\sf Y}_{k-1}}$ & $\cdots$
	\end{tabular}
	\ .
\end{equation*}
On the other hand, by assumption
\begin{equation*}
	\tf_i \bsfX= \quad
	\begin{tabular}{cccccc}
		$\cdots$ & $k+2$ & $k+1$ & $k$   & $k-1$ & $\cdots$ \\ \hline
		$\cdots$ & ${\sf x}_{k+2}$ & $1$ & ${\sf x}_{k}-1$ & ${\sf x}_{k-1}$ & $\cdots$ \\ 
		$\cdots$ & ${\sf y}_{k+2}$ & ${\sf y}_{k+1}$ & $2$ & ${\sf y}_{k-1}$ & $\cdots$
	\end{tabular}
	\ .
\end{equation*}
Hence we conclude $\tf_i \mc{F}(\bsfX) = \mc{F}(\tf_i \bsfX)$.
\medskip

{\it Case 2}.
Suppose that the $+$ of $+^{{\sf x}_{k}}$ is maximally moved to $+^{\widetilde{\sf X}_{k+1}}$ in $ \sigma_i(\mc{F}(\bsfX))$ by $\mc{F}$.
Then there are two possibilities: the $+$ is moved by (II) or by (III) as follows.
\smallskip

{\it Case 2-1}. Suppose that the $+$ is moved by (II).
Then ${\sf z}_k=1$ and $\re{({\sf x}_{k+1}-{\sf z}_{k+1}+1)}=1$.
Since ${\sf z}_k=1$ and $\qo{{\sf y}_{k-1}}=0$, we have $\re{{\sf y}_k}=1$.

We claim that there exists a $+$ in $+^{\widetilde{\sf X}_{k+1}}$ coming from $+^{{\sf x}_{k}}$ which is not canceled.
By Lemma \ref{lem: signature lemma}, these $+$'s can only be canceled with $-$'s coming from $-^{{\sf y}_{k-1}}$.
If ${\sf y}_{k-1}=1$ and ${\sf z}_{k-1}=1$, then one $+$ is canceled with one $-$,
but since ${\sf x}_k \geq 2$, ${\sf X}_k={\sf x}_k-{\sf z}_k+{\sf z}_{k-1}\geq 2$, 
so $+^{\widetilde{\sf X}_{k+1}}$ contains at least two $+$'s coming from $+^{{\sf x}_k}$.
Hence, the claim holds, which implies $\tf_i \mc{F}(\bsfX)= \mc{F}(\tf_i \bsfX)$.
	\smallskip
	
{\it Case 2-2}. Suppose that the $+$ is moved by (III). 
	Then ${\sf z}_k=0$ and ${\sf x}_{k}+{\sf z}_{k-1}\geq 2$.
	We have $\re{{\sf y}_{k}}=0$ since ${\sf z}_k=0$.

	We claim that ${\sf x}_k \geq 2$.
	Indeed, if ${\sf x}_k=1 $, then ${\sf y}_{k-1}=0$.
	Since $+^{{\sf x}_k}$ is not canceled in $\ov{\sigma_i}(\bsfX)$,
	we have ${\sf x}_{k-1} \geq {\sf y}_{k-2}$, so ${\sf y}_{k-2} \leq 1$,
	and then ${\sf z}_{k-1}=\mathrm{{\sf x}_{k-1},\re{{\sf y}_{k-1}}+2\qo{{\sf y}_{k-2}}}=0$, 
	contradicting ${\sf x}_{k}+{\sf z}_{k-1}\geq 2$. 
	Hence ${\sf x}_k \geq 2$, and $+^{\widetilde{\sf X}_{k+1}}$ contains an uncanceled $+$ from $+^{{\sf x}_k}$,
	giving $\tf_i \mc{F}(\bsfX)= \mc{F}(\tf_i \bsfX)$.
\medskip

{\it Case 3}.
Suppose that none of $+$'s in $+^{{\sf x}_k}$ is moved under $\mc{F}$. 
Then ${\sf z}_k=0$ and ${\sf x}_k+{\sf z}_{k-1} \leq 1$.
Consequently, we have $\re{{\sf y}_{k}}=0$.
Since ${\sf x}_k+{\sf z}_{k-1} \leq 1$, we have ${\sf x}_k=1$, ${\sf z}_{k-1}=0$, and ${\sf y}_{k-1}=0$.
Hence, $\tf_i \mc{F}(\bsfX)= \mc{F}(\tf_i \bsfX)$.
\qed

\begin{ex}
{\em
We continue Example \ref{ex: signature lemma}(1). Then
\begin{equation*}
	\bsfX=
    \begin{tabular}{c|cccccc}
    	& $\cdots$ & $2$ & $1$ & $0$ & $-1$  & $\cdots$ \\ \hline 
        $i$ & $\cdots$ & $0$ & $3$ & $2$ & $0$  & $\cdots$\\ 
        $i+1$ & $\cdots$ & $0$ & $2$ & $1$ & $0$  & $\cdots$
    \end{tabular}
    \quad
    \leadsto
    \quad
    \mc{F}(\bsfX)=
	\begin{tabular}{c|cccccc}
				& $\cdots$ & $2$ & $1$ & $0$ & $-1$  & $\cdots$ \\ \hline
		$i$     & $\cdots$ & $4$ & $0$ & $1$ & $0$  & $\cdots$ \\ 
		$i+1$   & $\cdots$ & $2$ & $1$ & $0$ & $0$  & $\cdots$
	\end{tabular}
\end{equation*}	
so that
\begin{equation*}
	\tf_i \mc{F}(\bsfX)=
	\begin{tabular}{c|cccccc}
		& $\cdots$ & $2$ & $1$ & $0$ & $-1$  & $\cdots$\\ \hline
		$i$ & $\cdots$ & $3$ & $0$ & $1$ & $0$ & $\cdots$ \\ 
		$i+1$ & $\cdots$ & $3$ & $1$ & $0$ & $0$ & $\cdots$
	\end{tabular}
	\ .
\end{equation*}		
On the other hand, we have 
\begin{equation*}
	\tf_i \bsfX=
    \begin{tabular}{c|cccccc}
    	& $\cdots$ & $2$ & $1$ & $0$ & $-1$ & $\cdots$ \\ \hline 
        $i$ & $\cdots$ & $0$ & $2$ & $2$ & $0$ & $\cdots$  \\ 
        $i+1$ & $\cdots$ & $0$ & $3$ & $1$ & $0$  & $\cdots$
    \end{tabular}
	\quad
	\leadsto
	\quad
	\mc{F}(\tf_i \bsfX)=
	\begin{tabular}{c|cccccc}
		& $\cdots$ & $2$ & $1$ & $0$ & $-1$ & $\cdots$ \\ \hline
		$i$ & $\cdots$ & $3$ & $0$ & $1$ & $0$ & $\cdots$  \\ 
		$i+1$ & $\cdots$ & $3$ & $1$ & $0$ & $0$ & $\cdots$
	\end{tabular}
\end{equation*}
where $\mc{F}(\tf_i \bsfX)$ coincides with $\tf_i \mc{F}(\bsfX)$.
This is an example of {\it Case 2} (b) in the above proof.
}
\end{ex}

\end{document}